\documentclass[12pt, twoside]{article}
\usepackage{amsmath,amsthm,amssymb,enumerate,times}

\usepackage{
verbatim,
graphicx,
epstopdf,
mathtools,
hyperref
}

\theoremstyle{definition}
\newtheorem{thm}{Theorem}[section]

\newtheorem{rem}[thm]{Remark}

\numberwithin{equation}{section}

\newcommand{\subjclass}[1]{\bigskip\noindent\emph{2010 Mathematics Subject Classification:}\enspace#1}
\newcommand{\keywords}[1]{\noindent\emph{Keywords:}\enspace#1}

\def\div{\nabla \cdot }
\def\vec[#1,#2]{\left(\begin{array}{c}{#1}\\{#2}\end{array}\right)}
\def\vecc[#1,#2,#3]{\left(\begin{array}{c}{#1}\\{#2}\\{#3}\end{array}\right)}
\def\mat[#1,#2,#3,#4]{\left(\begin{array}{cc}{#1}&{#2}\\{#3}&{#4}\end{array}\right)}
\def\matt[#1,#2,#3,#4,#5,#6,#7,#8,#9]{\left(\begin{array}{ccc}{#1}&{#2}&{#3}\\{#4}&{#5}&{#6}\\{#7}&{#8}&{#9}\\\end{array}\right)}
\def\bbb[#1]{\boldsymbol{#1}}
\def\mmm[#1]{\mathcal{#1}}
\def\R{\mathbb{R}} 

\def\N{\mathbb{N}} 
\def\dd{ \,\text{d}}

\def\t{\theta}
\def\p{\partial}
\def\md[#1]{\partial_t^{\bullet ({#1})}}
\def\nd[#1]{\partial_t^{\circ ({#1})}}

\def\O[#1]{\Omega^{({#1})}}
\def\G[#1,#2]{\Gamma^{({#1},{#2})}}
\def\T[#1,#2,#3]{T^{({#1},{#2},{#3})}}
\def\Textern[#1,#2]{T^{({#1},{#2})}_{\Omega}}
\def\Q[#1,#2,#3,#4]{Q^{({#1},{#2},{#3},{#4})}}
\def\Qextern[#1,#2,#3]{Q^{({#1},{#2},{#3})}_{\Omega}}
\def\n[#1,#2]{\bbb[\nu]^{({#1},{#2})}}
\def\nextern{\bbb[\nu]_{\Omega}}
\def\m[#1,#2,#3]{\bbb[\mu]^{({#1},{#2},{#3})}}
\def\mextern[#1,#2]{\bbb[\mu]^{({#1},{#2})}_{\Omega}}
\def\tang[#1,#2]{\bbb[\tau]^{({#1},{#2})}}
\def\curv[#1,#2]{\bbb[\kappa]^{({#1},{#2})}}
\def\scurv[#1,#2]{\kappa^{({#1},{#2})}}
\def\u[#1,#2]{u^{({#1},{#2})}}
\def\ndu[#1,#2]{\partial_t^{\circ (u^{({#1},{#2})})}}
\def\uu[#1,#2,#3]{\bbb[u]^{({#1},{#2},{#3})}}
\def\nabg[#1,#2]{\nabla_{\Gamma^{({#1},{#2})}}} 
\def\nabd[#1,#2]{\nabla_{\Gamma^{({#1},{#2})}} \cdot}
\def\cho[#1]{\chi_{\Omega^{({#1})}}}
\def\chg[#1,#2]{\chi_{\Gamma^{({#1},{#2})}}}
\def\d[#1]{\delta_{{#1}}}

\def\r[#1]{\overline{\rho}^{({#1})}}
\def\v{\bbb[v]}
\def\vel[#1]{\bbb[v]^{(i)}}
\def\et[#1]{{\eta}^{({#1})}}
\def\Ti{\bbb[T]^{(i)}}
\def\dimstr{\bbb[T]}
\def\c[#1]{c^{({#1})}}
\def\cc[#1,#2]{c^{({#1},{#2})}}
\def\Jc[#1]{\bbb[j]_c^{({#1})}}
\def\Jcc[#1,#2]{\bbb[j]_c^{({#1},{#2})}}
\def\Mc[#1]{M_c^{({#1})}}
\def\Mcc[#1,#2]{M_c^{({#1},{#2})}}
\def\g[#1]{g_{#1}}
\def\gg[#1,#2]{\gamma_{{#1},{#2}}}
\def\sig[#1,#2]{\sigma_{{#1},{#2}}}
\def\sigt[#1,#2]{\tilde{\sigma}_{{#1},{#2}}}
\def\lamt[#1]{\tilde{\lambda}_{{#1}}}

\def\e{\varepsilon}

\def\einvs{\tfrac{1}{\varepsilon}}
\def\rh[#1]{{\rho^{({#1})}_\e}}
\def\vv{{\bbb[v]_\e}}
\def\vve[#1]{{\bbb[v]_{{#1}}}}
\def\mdvve[#1]{\partial_t^{\bullet (\bbb[v]_{{#1}})}}
\def\vvE[#1]{{\bbb[V]_{{#1}}}}
\def\vvT[#1]{{\hat{\bbb[V]}_{{#1}}}}

\def\pp{{p_\e}}
\def\ppe[#1]{{p_{{#1}}}}
\def\ppE[#1]{{P_{{#1}}}}
\def\ppT[#1]{{\hat{P}_{{#1}}}}
\def\tpe{{\tilde{p}_\e}}
\def\qq{{q_\e}}
\def\qqe[#1]{{q_{{#1}}}}
\def\qqE[#1]{{Q_{{#1}}}}
\def\qqT[#1]{{\hat{Q}_{{#1}}}}
\def\Jq{\bbb[j]_{\qq}}
\def\Jqe[#1]{\bbb[j]_{q,#1}}
\def\JqE[#1]{\bbb[J]_{q,#1}}
\def\JqT[#1]{\hat{\bbb[J]}_{q,#1}}
\def\JJc[#1]{\bbb[j]_{c,\e}^{({#1})}} 
\def\JJcc[#1,#2]{\bbb[j]_{c,\e}^{({#1},{#2})}}
\def\cgg[#1,#2]{\check{\gamma}_{{#1},{#2}}}
\def\a[#1,#2]{a_{{#1},{#2}}}
\def\w[#1,#2]{\overline{w}_{{#1},{#2}}}
\def\del[#1,#2]{\delta_{{#1},{#2}}}
\def\mba{{\mu_\e}}
\def\mbae[#1]{{\mu_{{#1}}}}
\def\mbaE[#1]{{M_{{#1}}}}
\def\mbaT[#1]{{\hat{M}_{{#1}}}}
\def\mb[#1]{{\mu^{({#1})}_\e}}
\def\mbe[#1,#2]{{\mu^{({#1})}_{{#2}}}}
\def\mbE[#1,#2]{{M^{({#1})}_{{#2}}}}
\def\mbT[#1,#2]{{\hat{M}^{({#1})}_{{#2}}}}

\def\Mb[#1]{M^{({#1})}}
\def\L[#1,#2]{\mmm[L]^{({#1},{#2})}}
\def\pha{{\varphi_\e}}
\def\npha{{\nabla\varphi_\e}}
\def\phae[#1]{{\varphi_{{#1}}}}
\def\phaE[#1]{{\Phi_{{#1}}}}
\def\phaT[#1]{{\hat{\Phi}_{{#1}}}}
\def\ph[#1]{{\varphi^{({#1})}_\e}} 
\def\phe[#1,#2]{{\varphi^{({#1})}_{{#2}}}} 
\def\phE[#1,#2]{{\Phi^{({#1})}_{{#2}}}} 
\def\phT[#1,#2]{{\hat{\Phi}^{({#1})}_{{#2}}}} 
\def\Jpa{\bbb[j]_\pha}
\def\Jp[#1]{\bbb[j]_\pha^{({#1})}}
\def\Jpe[#1,#2]{\bbb[j]_{\varphi,{#2}}^{({#1})}}
\def\JpE[#1,#2]{\bbb[J]_{\varphi,{#2}}^{({#1})}}
\def\JpT[#1,#2]{\hat{\bbb[J]}_{\varphi,{#2}}^{({#1})}}
\def\Jbar{\overline{\bbb[j]}_{\e}}
\def\be[#1]{\bbb[e]_{{#1}}}
\def\mplus{p}
\def\mminus{n}
\def\param[#1,#2]{\pi^{({#1},{#2})}}
\def\mIn{n}
\def\mIp{p}
\def\mIm{m}
\def\mIm{r}
\def\triple[#1,#2,#3]{\boldsymbol{\theta}^{({#1},{#2},{#3})}}
\def\paramtri[#1,#2,#3]{\boldsymbol{\pi}^{({#1},{#2},{#3})}}
\def\Y[#1,#2,#3]{Y^{({#1},{#2},{#3})}}
\def\projY[#1,#2,#3]{\boldsymbol{P}_{Y^{({#1},{#2},{#3})}}}
\def\tangtri[#1,#2,#3]{\boldsymbol{\tau}^{({#1},{#2},{#3})}}
\def\quadruple[#1,#2,#3,#4]{\boldsymbol{\eta}^{({#1},{#2},{#3},{#4})}}

\begin{document}



\title{Phase field modelling of surfactants in multi-phase flow}

\author{
Oliver R. A. Dunbar\thanks{o.dunbar.1@warwick.ac.uk}\\
Mathematics Institute, University of Warwick, \\
Coventry CV4 7AL, United Kingdom\\
Kei Fong Lam\thanks{kei-fong.lam@math.cuhk.edu.hk}\\
Department of Mathematics,\\
The Chinese University of Hong Kong,\\
Shatin, N.T., Hong Kong\\
Bj\"orn Stinner\thanks{bjorn.stinner@warwick.ac.uk, corresponding author}\\
Mathematics Institute and Centre for Scientific Computing, \\
University of Warwick, Coventry CV4 7AL, United Kingdom
}

\date{\today}

\maketitle


\begin{abstract}
A diffuse interface model for surfactants in multi-phase flow with three or more fluids is derived. A system of Cahn-Hilliard equations is coupled with a Navier-Stokes system and an advection-diffusion equation for the surfactant ensuring thermodynamic consistency. By an asymptotic analysis the model can be related to a moving boundary problem in the sharp interface limit, which is derived from first principles. Results from numerical simulations support the theoretical findings. The main novelties are centred around the conditions in the triple junctions where three fluids meet. Specifically the case of local chemical equilibrium with respect to the surfactant is considered, which allows for interfacial surfactant flow through the triple junctions. 
 

\end{abstract}

\section{Introduction}
\label{sec:intro}

Surfactants (surface active agents) are chemicals that, when dissolved in a system of multiple immiscible fluids, tend to form layers at the fluid-fluid interfaces and thus reduce the surface tension. Such manipulation is exploited in nature and industry, and we refer to \cite{Mye05,Tad05} for overviews of the vast applicability of surfactants and to \cite{BruDelHarETAL15,KarCraMat11,KraSchKowTre15,KumBhaKulKum15,SeiFouFroDef03} for specific applications involving more than two phases. 

Several approaches to such problems based on the representation of the interfaces by hypersurfaces (here called {\em sharp interface models}) are available, among which we mention interface tracking methods \cite{Hym84, KhaTor11, LaiTseHua08, MurTry08, MurTry14, TryBunEsmetal01}, volume-of-fluid methods \cite{AlkBot09, HirNic81, JamLow04, NicHirHot81}, and ALE methods \cite{BarGarNue15, GanTob09, YanJam09}, see also the books \cite{BotReu17, GroReu11}. In general, the fluid-fluid interfaces undergo changes of topology, which may manifest as the breakup of droplets, pinching, coalescence, or cusp formation or tip-streaming driven by Marangoni forces. To overcome the analytical and numerical complications associated to such events one can turn to interface-capturing methods such as level-set methods \cite{AdaSet95,SusSmeOsh94, XuHuaLaiLi14, XuLiLowZha06}, or \emph{diffuse interface} approaches, which comprise the {\em phase field methodology}. 

In this work we address the phase field modelling of surfactant dynamics in multi-phase flow with more than two fluids. The classical description of fluid-fluid interfaces with hypersurfaces is replaced by one with thin transition layers of a thickness that scales with a small parameter $\e$. Within these thin layers, some form of microscopic mixture of the macroscopically immiscible fluids is allowed. One then introduces \emph{order parameters} or \emph{phase field variables} that serve to distinguish between the bulk phases, where the phase fields are close to constants, and the interfacial layers, across which the phase fields change values quickly but smoothly. 

The notion of diffuse interfaces dates at least back to van der Waals \cite{Row79}. In \cite{HohHal77}, {\em model H} couples a Cahn--Hilliard equation with a Navier--Stokes system. Subsequent efforts have been directed to extend this type of model with regards to non-matched densities \cite{LowTru98}, divergence-free mixture velocities \cite{DinSpeShu07}, thermodynamic consistency \cite{AbeGarGru12,GurPolVin96}, and flows with more than two fluids \cite{BanNue17, BoyLap06, Don14, Don15, Kim09, Kim12, KimLow05}. Regarding the inclusion of surfactants we refer to \cite{EngDoQAmbTor13, LiKim12, LiuZha10, TeiSonLowVoi11, vdSvdG06, YunLiKim14}, all of which are restricted to two fluids. 

Our phase field model builds up on \cite{GarLamSti14}, where surfactants in two-phase flow are studied within a free energy framework. The focus of this study is on the {\em instantaneous adsorption} regime when the adsorption-desorption process between interfacial surfactant and bulk surfactant in the adjacent sublayers occurs at a much faster timescale relative to other diffusive or convective processes in the system. The relation between interfacial and bulk surfactant is commonly described by {\em isotherms} \cite{EasDal00}. Within a free energy framework this local chemical equilibrium condition can be expressed as an equality of the chemical potentials of the (surfactant dependent) interfacial and bulk free energies \cite{DiaAnd96,RowWid13}. 

In the case of more than two fluids the fluid-fluid interfaces can meet at {\em triple junctions}, which are points (in two spatial dimensions, $d=2$) or lines ($d=3$), the latter possibly forming quadruple points if four or more fluids are present. Mass flux of interfacial surfactant through the triple junctions is of relevance in applications \cite{KarCraMat11}. We here make the assumption that no surfactant mass is associated with the triple junction and that the assumption of local chemical equilibrium at the interfaces extends to the triple junction. In more mathematical terms, the net flux into a triple junction from the adjacent interfaces is zero and the interfacial surfactant chemical potentials match up at the triple junction.

Under some convexity assumptions on the bulk and surface free energies the local chemical equilibrium assumption enables us to introduce a single continuous (chemical) potential field in which the balance laws for the bulk and interfacial surfactants can be expressed, as can the surfactant dependent surface tensions. These coupled equations can be formulated in a distributional form in the context of sharp interface models. Here, we can follow the lines of \cite{Alt09}, which covers the two-phase case and can be extended to account for triple junctions. This form allows for an incorporation into a phase field model with thin layers representing the interfaces. It is achieved by smoothing the distributions associated with the bulk domains and interfaces in terms of the phase field variables. We refer to \cite{GarLamSti14,LiLowRaeVoi09,RaeVoi06,TeiLXLowWanVoi09,TeiSonLowVoi11} for the ideas and to \cite{AbeLamSti15,BurElvSch17,EllSti09} for rigorous analytical investigations in the two-phase case. 

One of the challenges in the case of multiple fluids is to choose suitable smoothing functions such that (a) the interfacial surfactant equation for each specific fluid-fluid interface is consistently approximated and (b) the conditions in the triple junctions are consistently approximated, as the diffuse interface thickness converges to zero. A key ingredient to solve this problem are phase field potentials that avoid {\em third phase contributions} at interfaces, i.e., in the interfacial layer between two phases only the phase field variable associated with these phases are present \cite{BoyLap06,BoyMin14,BreMas17,GarNesSto99B,Sti05}. This allows for a precise localisation of fields and functional dependencies that are supposed to be present at a specific interface or a triple junction only. 

It is also desirable that the smoothing leads to a model with good structural properties such as thermodynamic consistency and a solenoidal velocity field, which are beneficial for numerical approximations. As in \cite{GarLamSti14} our model for the multi-phase flow is based on \cite{AbeGarGru12}. This approach is extended to account for multiple fluids. As in \cite{GurPolVin96} it assumes that, within the interfacial layers where the fluids mix, inertia and kinetic energy due to the motion of the constituent fluids relative to the gross motion of the mixture is negligible. Often, the mass-averaged velocity is chosen to define the gross motion \cite{GuoLinLowWis17,Kim09,Kim12,KimLow05,LowTru98}. But taking the volume-averaged velocity as in \cite{BanNue17,DinSpeShu07,Don14,Don15} leads to a divergence-free velocity field. Moreover, we can ensure that the {\em calibration} of our phase field model is convenient in the following sense: Parameters (fluid densities, viscosities, surfactant diffusivities, etc) and relations in the sharp interface model (dependence of the surface tension on the surfactant, etc) directly reappear in the phase field model, no adjustment or rescaling is required. 

In Section \ref{sec:SIM} we derive the sharp interface model that we aim to approximate with the phase field methodology. In particular, the conditions in the triple junctions are motivated within a free energy framework that is discussed in detail. A summary of the governing equations can be found in Section \ref{sec:sumSIM}. The phase field model is derived in Section \ref{sec:DIM} and follows a similar procedure by postulating balance equations and free energies, and then closing the equations accounting for the instantaneous adsorption assumption mentioned above and ensuring thermodynamic consistency. A summary of the model is contained in Section \ref{sec:sumDIM}. An asymptotic analysis based on matching suitable expansions in the small interfacial thickness parameter $\e$ is presented in Section \ref{sec:asymp}. As one of the main novelties we show that the conditions for the surfactant in the triple junction indeed are obtained in the sharp interface limit. The conditions at the fluid-fluid interfaces have been analysed in \cite{GarLamSti14}, to which the present multi-phase case arguably reduces if the third phase contributions, which were mentioned above, can be avoided. However, the details are presented as they are required for the analysis around the triple junctions. We have performed some numerical simulations on a qualitative level in order to validate and support the theoretical results of the asymptotic analysis, see Section \ref{sec:numsim}.

\section{Sharp interface model}
\label{sec:SIM}

In deriving the free boundary problem, which we intend to approximate with a phase field model, we extend \cite{GarLamSti14} by accounting for multiple phases. This implies that the conditions at points where several phases meet have to be discussed. A general study of balance equations in three phase systems including bulk, surface, and triple line fields is presented in \cite{BotPru16}. We re-state some of the theory in order to introduce our notation and to define and briefly discuss our specific closing conditions.

\subsection{Setting}
\label{sec:SIMsetting}

\begin{figure}[ht]
 \centering
 \includegraphics[width=0.9\textwidth]{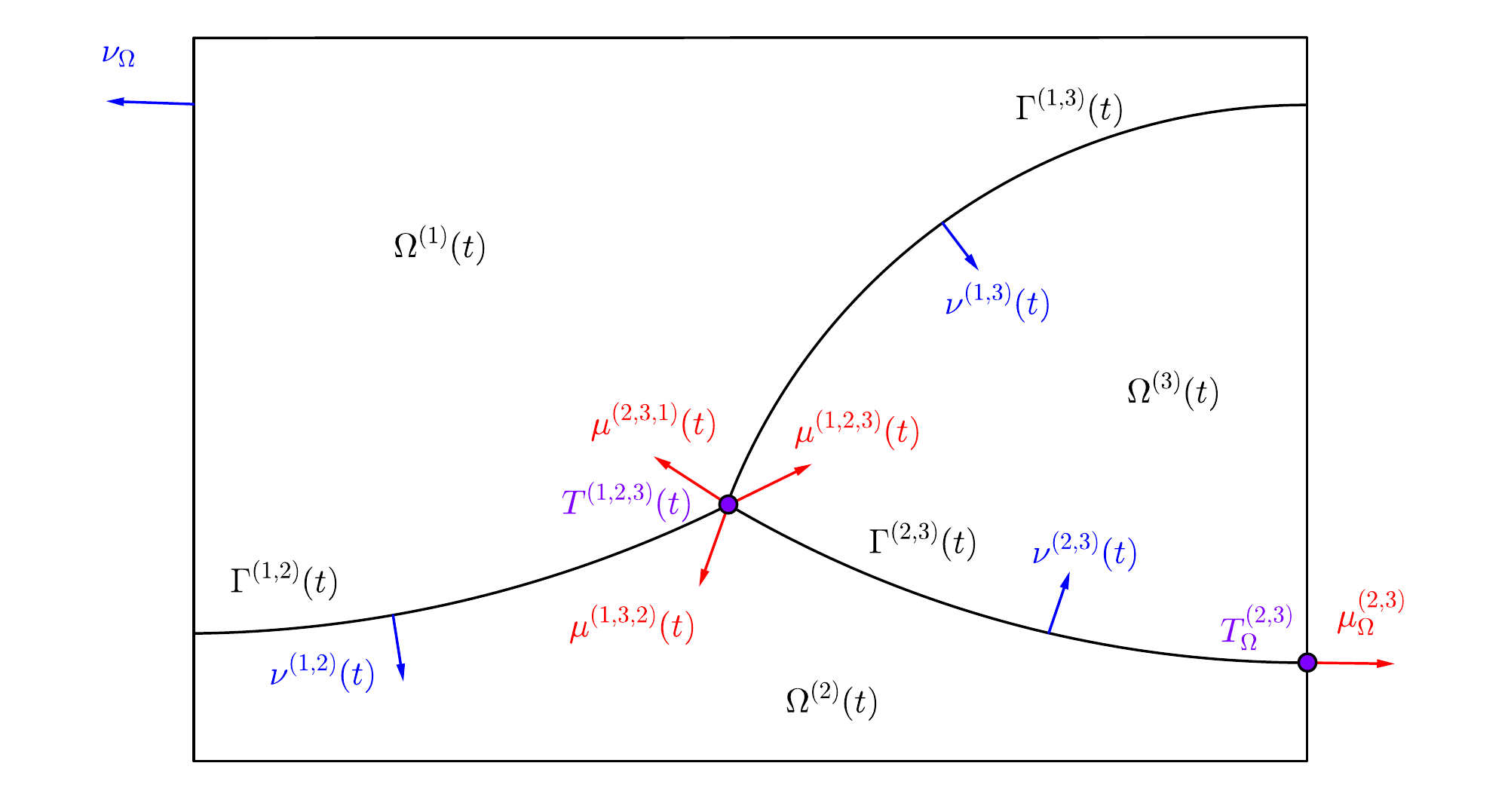}
 \caption{\em Illustration of the setting described in Section \ref{sec:SIMsetting}.} \label{fig:subdivision}
\end{figure}

Let $\Omega \subset \R^d$, $d \in \{1,2,3\}$, be a bounded domain and $I=[0,T)$, $T \in (0,\infty]$ be a time interval. Assume that $\Omega$ is partitioned by moving hypersurfaces $\G[i,j]$ into $M$ time dependent open subdomains $\O[k]$, $i,j,k \in \{ 1, \dots, M \}$ (for brevity, here and in the following we neglect noting the time dependence unless it is required). Intersections of three hypersurfaces are denoted by $\T[i,j,k]$ and form triple points ($d=2$) or form triple lines ($d=3$) ending at quadruple points $\Q[i,j,k,l]$, $i,j,k,l \in \{ 1, \dots, M \}$. For simplicity, with regards to $\T[i,j,k]$ we will only talk about \emph{triple junctions} in the following. Similarly, on the external boundary $\p \Omega $ there are triple points or lines $\Textern[i,j]$ with quadruple points $\Qextern[i,j,k]$ if $d=3$. The unit normal on $\G[i,j]$ pointing out of $\O[i]$ into $\O[j]$ is denoted by $\n[i,j]$ and by $\nextern$ on $\p \Omega$. For the conormal of $\G[i,j]$ at $\T[i,j,k]$ pointing into $\O[k]$ we write $\m[i,j,k]$, and we write $\mextern[i,j]$ for it on $\p \Omega$. Figure \ref{fig:subdivision} is a sketch of a configuration as we have it in mind. 

The whole configuration is transported by a continuous velocity field $\v : [0,T) \times \Omega \to \R^d$. This implies that 
\begin{equation} \label{eq:intvelcont}
 [ \v ]_i^j = 0, \quad \u[i,j] = \v \cdot \n[i,j] \qquad \text{on } \G[i,j],
\end{equation}
where $[ \cdot ]_i^j = (\cdot)^{(j)} - (\cdot)^{(i)}$ stands for the jump from domain $\O[i]$ into $\O[j]$ across the interface $\G[i,j]$ and $\u[i,j]$ is the normal velocity of $\G[i,j]$ in direction $\n[i,j]$, and furthermore
\begin{equation} \label{eq:triplevelcont}
 \uu[i,j,k] = \bbb[P]_{(\T[i,j,k])^\perp} \v \qquad \text{at } \T[i,j,k],
\end{equation}
where $\uu[i,j,k]$ is the normal velocity of $\T[i,j,k]$ and $\bbb[P]_{(\T[i,j,k])^\perp}$ is the projection to the space normal to $\T[i,j,k]$.

The surface derivative and divergences along the hypersurface $\G[i,j]$ are denoted by $\nabg[i,j]$ and $\nabd[i,j]$, respectively. Further notation concerns the material derivative for a field $w : [0,T) \times \Omega \to \R$,
\begin{equation} \label{eq:defmatvel}
 \md[\v] w \coloneqq \p_t w + \v \cdot \nabla w,
\end{equation}
where we remark that this operator is well-defined for fields restricted to a hypersurface $\G[i,j]$. For such fields $w$ we also consider the normal time derivative,
\begin{equation} \label{eq:defntd}
 \ndu[i,j] w \coloneqq \p_t w + \u[i,j] \n[i,j] \cdot \nabla w,
\end{equation}
and note that
\begin{equation} \label{eq:identtd}
 \md[\v] w = \ndu[i,j] w + \v \cdot \nabg[i,j] w.
\end{equation}
Some identities such as a transport identity on evolving surfaces and integration by parts formula on surfaces are stated in the Appendix. Finally, by $\curv[i,j]$ we denote the mean curvature vector of $\G[i,j]$.

\subsection{Balance equations}

$M\in \N$ represents the number of fluids, which are assumed to be immiscible, incompressible, and Newtonian. Correspondingly, for each fluid there is a subdomain $\O[i]$ indicating the regions occupied by the fluid. Denoting by $\r[i]$ and $\et[i]$ the mass density and viscosity of fluid $i \in \{ 1, \dots, M \}$, the mass and linear momentum balances in $\O[i]$ read
\begin{align}
 \div \v &=0, \label{eq:massbal}\\
 \p_t(\r[i] \v) + \div(\r[i] \v \otimes \v) &= \div \Ti, \label{eq:mombal} \\
 \Ti &= -p \bbb[I] + 2\et[i] D(\v), \label{eq:deftensor}
\end{align}
with the rate of deformation tensor $D(\v) = \frac{1}{2}((\nabla \v) +(\nabla \v)^\top)$, pressure $p$ and the identity tensor $\bbb[I]$.

For simplicity we only consider a single surfactant. Its bulk and surface mass density in each subdomain $\O[i]$ and hypersurface $\G[i,j]$ are denoted by $\c[i] : \O[i] \to \R$ and $\cc[i,j] : \G[i,j] \to \R$, respectively. The effect on the total mass density is assumed to be small. Thus, only mass balances are considered and any effects on the momentum are neglected. Following the derivation in \cite{GarLamSti14} the surfactant mass balance equations read
\begin{align}
 \md[\v] \c[i] &= -\div \Jc[i], & \text{in } \O[i], \label{eq:sagbulk} \\
 \md[\v]\cc[i,j] + \cc[i,j] \nabg[i,j] \cdot \v &= -\nabd[i,j] \Jcc[i,j] + q_{AD}^{(i,j)} & \text{on } \G[i,j], \label{eq:sagsurf}
\end{align}
where $\Jc[\cdot]$ and $\Jcc[\cdot,\cdot]$ are associated bulk and surface diffusive fluxes, and with the adsorption-desorption flux
\[
 q_{AD}^{(i,j)} = \Jc[i] \cdot \n[i,j] + \Jc[j] \cdot \n[j,i] = [ \Jc[\cdot] ]_j^i \cdot \n[i,j].
\] 

Let now $V(t) \subset \Omega$ be an arbitrary material test volume with external unit normal $\bbb[\nu]_V$. Using the transport identities \eqref{eq:Rey1} and \eqref{eq:Rey2} and the above identities \eqref{eq:sagbulk} and \eqref{eq:sagsurf} one can then derive that
\begin{multline}
 \frac{\dd}{\dd t} \Big{(} \sum_{i} \int_{V \cap \O[i]} \c[i] + \sum_{i<j} \int_{V \cap \G[i,j]} \cc[i,j] \Big{)} \\
 = -\sum_i \int_{\p V \cap \O[i]} \Jc[i] \cdot \bbb[\nu]_V - \sum_{i<j} \int_{\p V \cap \G[i,j]} \Jcc[i,j] \cdot \bbb[\mu]^{(i,j)}_V \\
 + \sum_{i<j<k} \int_{V \cap \T[i,j,k]} \Jcc[i,j] \cdot \m[i,j,k] + \Jcc[j,k] \cdot \m[j,k,i] + \Jcc[k,i] \cdot \m[k,i,j], \label{eq:motitripcond}
\end{multline}
where $\bbb[\mu]^{(i,j)}_V$ is the external conormal of $V \cap \G[i,j]$ on $\p V \cap \G[i,j]$. We omit the details of the calculation as similar techniques are presented in Section \ref{sec:SIMfreeenergy} within a more extensive calculation for the free energy.

We now make the assumption that no surfactant mass is stored at the triple points if $d=2$ nor at triple lines or quadruple points if $d=3$. For the diffusive surface fluxes this means that
\begin{equation}
\Jcc[i,j] \cdot \m[i,j,k] + \Jcc[j,k] \cdot \m[j,k,i] + \Jcc[k,i] \cdot \m[k,i,j] =0 \quad \text{at } \T[i,j,k]. \label{eq:sagtripmasscond} 
\end{equation}
The last line of \eqref{eq:motitripcond} then vanishes and the remainder of this identity thus reads that the (instantaneous) change of surfactant mass in the material volume $V$ (left-hand side) is given by the surfactant mass flux across $\p V$ (right-hand side).

\subsection{Free energy}
\label{sec:SIMfreeenergy}

In order to close the balance equations and relate the fluxes to the conserved fields we consider an energetic framework. With regards to the surfactant we postulate bulk free energies $\g[i](\c[i])$ and surface free energies $\gg[i,j](\cc[i,j])$ that are strictly convex, i.e. $\g[i]''>0$ and $\gg[i,j]''>0$. The total free energy including the kinetic energy is then
\begin{equation} \label{eq:SIMenerg}
 E \coloneqq \sum_{i=1}^{M} \int_{\O[i]} \left(\frac{\r[i]}{2} |\v|^2 + \g[i](\c[i])\right) + \sum_{\substack{i,j=1 \\ i < j}}^{M} \int_{\G[i,j]} \gg[i,j](\cc[i,j]).
\end{equation}
Related to the surface free energy we define the {\em surface tensions}:
\begin{equation} \label{eq:sig}
\sig[i,j](\cc[i,j]) \coloneqq \gg[i,j](\cc[i,j]) - \cc[i,j]\gg[i,j]'(\cc[i,j]).
\end{equation}

Let $V(t) \subset \Omega$ be an arbitrary material test volume with external unit normal $\bbb[\nu]_V$. Then thanks to the transport identities \eqref{eq:Rey1} and \eqref{eq:Rey2} and the incompressibility of the fluids \eqref{eq:massbal}
\begin{align}
 & \frac{\dd}{\dd t} \left(\sum_i \int_{V \cap \O[i]} \left(\frac{\r[i]}{2}|\v|^2 + \g[i](\c[i])\right) + \sum_{i < j}\int_{V\cap\G[i,j]} \gg[i,j](\cc[i,j]) \right) \nonumber \\
=& \sum_i \int_{V\cap\O[i]} (\r[i] \v \cdot \md[\v] \v + \g[i]'(\c[i]) \md[\v]\c[i]) + \sum_{i < j} \int_{V \cap \G[i,j]} (\gg[i,j]' \md[\v] \cc[i,j] + \gg[i,j] \nabg[i,j] \cdot \v), \nonumber \displaybreak[0] \\
\intertext{inserting the balance laws \eqref{eq:mombal}, \eqref{eq:sagbulk}, and \eqref{eq:sagsurf}, this is }
=& \sum_i \int_{V\cap\O[i]} (\v\cdot (\div\Ti)+ \g[i]'(-\div \Jc[i])) \nonumber \\
 & + \sum_{i < j} \int_{V\cap\G[i,j]} \gg[i,j]' \big{(} -\nabg[i,j] \cdot \Jcc[i,j] + [\Jc[\cdot]]_j^i\cdot\n[i,j] \big{)} + (-\cc[i,j] \gg[i,j]' + \gg[i,j]) \nabg[i,j] \cdot \v, \nonumber \displaybreak[0] \\
\intertext{using \eqref{eq:sig}, applying \eqref{div} where we note that $\Jcc[i,j] \cdot \curv[i,j] = 0$ as the flux is tangential, and using the symmetry of $\Ti$,}
=& \sum_i \Big{(} \int_{V \cap \O[i]} ( -\nabla \v \colon \Ti + \nabla \g[i]' \cdot \Jc[i] ) + \int_{\p V \cap \O[i]} ( \Ti \v - \g[i]' \Jc[i] ) \cdot \bbb[\nu]_V \Big{)} \nonumber \\
 & + \sum_i \sum_{j \neq i} \int_{V \cap \G[i,j]} ( \Ti \v - \g[i]' \Jc[i] ) \cdot \n[i,j] \nonumber \displaybreak[0] \\
 & + \sum_{i<j} \int_{V \cap \G[i,j]} \big{(} \nabg[i,j] \gg[i,j]' \cdot \Jcc[i,j] + [\gg[i,j]' \Jc[\cdot]]_j^i \cdot \n[i,j] - \nabg[i,j]\sig[i,j] \cdot \v - \sig[i,j] \curv[i,j] \cdot \v \big{)} \nonumber \\
 & + \sum_{i<j} \Big{(} \int_{\p V \cap \G[i,j]} (\gg[i,j]' \Jcc[i,j] + \sig[i,j] \v) \cdot \bbb[\mu]^{(i,j)}_V + \sum_{k \neq i,j} \int_{V \cap \T[i,j,k]} (-\gg[i,j]' \Jcc[i,j] + \sig[i,j] \v) \cdot \m[i,j,k] \Big{)}, \nonumber \displaybreak[0] \\
\intertext{with the external conormal $\bbb[\mu]^{(i,j)}_V$ of $\G[i,j]$ on $\p V$, rewriting the double sums we obtain}
=& \sum_i \int_{V \cap \O[i]} ( -D(\v) \colon \Ti + \nabla \g[i]' \cdot \Jc[i] ) + \sum_{i<j} \int_{V \cap \G[i,j]} \nabg[i,j] \gg[i,j]' \cdot \Jcc[i,j] \label{eq:dissterm1} \\
 & + \sum_{i<j} \int_{V \cap \G[i,j]} \big{[} (\gg[i,j]' - \g[(\cdot)]') \Jc[\cdot] \big{]}_j^i \cdot \n[i,j] \label{eq:dissterm2} \displaybreak[0] \\
 & + \sum_{i<j} \int_{V \cap \G[i,j]} \big{(} [\bbb[T]^{(\cdot)}]_j^i \n[i,j] - \nabg[i,j]\sig[i,j] - \sig[i,j] \curv[i,j] \big{)} \cdot \v \label{eq:forcebalterm1} \displaybreak[0] \\
 & + \sum_{i<j<k} \int_{V \cap \T[i,j,k]} \big{(} \sig[i,j] \m[i,j,k] + \sig[j,k] \m[j,k,i] + \sig[k,i] \m[k,i,j] \big{)} \cdot \v \label{eq:forcebalterm2} \displaybreak[0] \\
 & - \sum_{i<j<k} \int_{V \cap \T[i,j,k]} \big{(} \gg[i,j]' \Jcc[i,j] \cdot \m[i,j,k] + \gg[j,k]' \Jcc[j,k] \cdot \m[j,k,i] + \gg[k,i]' \Jcc[k,i] \cdot \m[k,i,j] \big{)} \label{eq:tripleterm} \displaybreak[0] \\
 & + \sum_i \int_{\p V \cap \O[i]} \big{(} ( \Ti \bbb[\nu]_V ) \cdot \v - \g[i]' \Jc[i] \cdot \bbb[\nu]_V \big{)} \label{eq:workterms1} \\
 & + \sum_{i<j} \int_{\p V \cap \G[i,j]} \big{(} -\gg[i,j]' \Jcc[i,j] \cdot \bbb[\mu]^{(i,j)}_V + \sig[i,j] \bbb[\mu]^{(i,j)}_V \cdot  \v \big{)}. \label{eq:workterms2}
\end{align}

\subsection{Instantaneous adsorption}
\label{subsec:intad}

We assume that the adsorption-desorption dynamics of the surfactant at the interfaces is fast and therefore may be considered as instantaneous at the time scale of the interface and fluid flow dynamics. These local equilibrium conditions result in relations between the surfactant densities in the sublayers close to interfaces with the interfacial densities, which are known as \emph{isotherms} \cite{EasDal00}. In terms of the chemical potentials $\g[i]'$ and $\gg[i,j]'$ these conditions read
\begin{equation}
\g[i]'(\c[i]) = \g[j]'(\c[j]) = \gg[i,j]'(\cc[i,j]) \quad \text{on } \G[i,j]. \label{eq:equilibsurf}
\end{equation}
In addition, we also assume a local chemical equilibrium at the triple junctions:
\begin{equation}
\gg[i,j]'(\cc[i,j]) = \gg[j,k]'(\cc[j,k]) = \gg[k,i]'(\cc[k,i]) \quad \text{at } T^{(i,j,k)}. \label{eq:equilibtriple}
\end{equation}
Thus, the chemical potential 
\begin{align*}
q \coloneqq 
\begin{cases}
\g[i]'(\c[i]) & \quad \mbox{in } \O[i], \\
\gg[i,j]'(\cc[i,j]) & \quad \mbox{on } \G[i,j]
\end{cases}
\end{align*} 
is continuous in $\Omega$. 

Recall that, by assumption, the free energies are convex as functions of the mass densities. Hence, they can be locally inverted so that we can express the surfactant bulk and surface mass densities in terms of $q$:
\begin{equation}\label{q}
 \cc[i,j](q) = (\gg[i,j]')^{-1}(q), \quad  \c[i](q) = (\g[i]')^{-1}(q).
\end{equation}
We can then also express the surface tension as a function of $q$:
\begin{equation} \label{eq:defsigt}
 \sigt[i,j](q) \coloneqq \sig[i,j](\cc[i,j](q)) = \gg[i,j](\cc[i,j](q)) - q \, \cc[i,j](q).
\end{equation}

We note that by \eqref{eq:equilibsurf} the term \eqref{eq:dissterm2} vanishes. Similarly, the condition \eqref{eq:equilibtriple} together with \eqref{eq:sagtripmasscond} ensures that \eqref{eq:tripleterm} vanishes.

\subsection{Further constitutive assumptions}
\label{sec:SIMconstitassum}
The terms in \eqref{eq:dissterm1} motivate us to define the surfactant fluxes by
\begin{align}
\Jc[i] &\coloneqq -\Mc[i] \nabla \g[i]'(\c[i]) = -\Mc[i] \nabla q & \quad & \mbox{in } \O[i], \label{eq:sagbulkflux} \\
\Jcc[i,j] &\coloneqq -\Mcc[i,j] \nabg[i,j] \gg[i,j]'(\cc[i,j]) = -\Mcc[i,j] \nabg[i,j] q & \quad & \mbox{on } \G[i,j], \label{eq:sagsurfflux}
\end{align}
with nonnegative mobilities $\Mc[i]$ and $\Mcc[i,j]$ that may be functions of the $\c[i]$ and the $\c[i,j]$, respectively, but are assumed to be constants for simplicity. 

At the interfaces $\G[i,j]$ we assume the force balances
\begin{equation} \label{eq:surfforcebal}
 [\bbb[T]^{(\cdot)}]_j^i \n[i,j] = \sig[i,j](\cc[i,j]) \curv[i,j] + \nabg[i,j] \sig[i,j](\cc[i,j]) = \sigt[i,j](q) \curv[i,j] + \nabg[i,j] \sigt[i,j](q),
\end{equation}
which mean that the stresses exerted by the fluids adjacent to the interfaces are counterbalanced by intrinsic forces, namely the surface tension forces $\sigt[i,j] \curv[i,j]$ and the \emph{Marangoni forces} $\nabg[i,j] \sigt[i,j]$. 

In the triple points or lines we assume the following balances of capillary forces:
\begin{equation} \label{eq:tripleforcebal}
 \sigt[i,j](q) \m[i,j,k] + \sigt[j,k](q) \m[j,k,i] + \sigt[k,i](q) \m[k,i,j] = 0.
\end{equation}
This triple junction condition is also known as \emph{Young's law}, see \cite{GarNesSto99A} for a discussion in the context of general anisotropic surface energies. In particular, it determines the angles at which the three phases meet at the triple junction. In the case $d=3$ the condition \eqref{eq:tripleforcebal} also fully determines the configuration and angles at the quadruple junctions $\Q[i,j,k,l]$, see \cite{BroGarSto98}, Section 3, for a discussion.

Condition \eqref{eq:tripleforcebal} is a local mechanical equilibrium condition, which may not always hold true. Indeed, wetting or spreading phenomena are of great relevance in many applications. The \emph{wetting} or \emph{spreading coefficients} \cite{HarFel22}
\begin{equation} \label{eq:defWettcoeff}
 \tilde{S}^{(i,j,k)}(q) \coloneqq \sigt[i,j](q) - \big{(} \sigt[i,k](q) + \sigt[j,k](q) \big{)},
\end{equation}
may be positive so that a thin layer of fluid $k$ between fluids $i$ and $j$ is energetically favourable to an $i$-$j$ interface. The condition \eqref{eq:tripleforcebal} then cannot be satisfied but other closing conditions, for instance, involving precursor films have to be postulated \cite{PopOshDieCaz12}. We will not cover the spreading case in the free boundary problem and the subsequent asymptotic analysis but note that some phase field models are able to deal with it \cite{BoyLap06}. 

Accounting for all constitutive assumptions \eqref{eq:equilibsurf}, \eqref{eq:equilibtriple}, \eqref{eq:sagbulkflux}, \eqref{eq:sagsurfflux}, \eqref{eq:surfforcebal}, and \eqref{eq:tripleforcebal} we obtain from \eqref{eq:dissterm1}--\eqref{eq:workterms2} that
\begin{align}
 & \frac{\dd}{\dd t} \left(\sum_i \int_{V \cap \O[i]} \left(\frac{\r[i]}{2}|\v|^2 + \g[i](\c[i](q))\right) + \sum_{i < j} \int_{V\cap\G[i,j]} \gg[i,j](\cc[i,j](q)) \right) \label{eq:changeenergy} \\
=& - \sum_i \int_{V \cap \O[i]} ( 2 \et[i] | D(\v) |^2 + \Mc[i] | \nabla q |^2 ) - \sum_{i<j} \int_{V \cap \G[i,j]} \Mcc[i,j] | \nabg[i,j] q |^2 \label{eq:dissterms} \displaybreak[0] \\
 & + \sum_i \int_{\p V \cap \O[i]} ( \Ti \v ) \cdot \bbb[\nu]_V + \sum_{i<j} \int_{\p V \cap \G[i,j]} \sig[i,j] \v \cdot \bbb[\mu]^{(i,j)}_V \label{eq:workterms} \\
 & - \sum_i \int_{\p V \cap \O[i]} q \Jc[i] \cdot \bbb[\nu]_V - \sum_{i<j} \int_{\p V \cap \G[i,j]} q \Jcc[i,j] \cdot \bbb[\mu]^{(i,j)}_V. \label{eq:energyfluxterms}
\end{align}
The terms in \eqref{eq:dissterms} are dissipative contributions to the change of energy. The terms in \eqref{eq:workterms} represent the working done on $V$ by the external fluid, and \eqref{eq:energyfluxterms} lists the loss (or gain) of energy due to the surfactant mass fluxes across $\p V$.

\subsection{Boundary conditions}
\label{sec:SIMBC}

The terms \eqref{eq:workterms} and \eqref{eq:energyfluxterms} also motivate \emph{natural} boundary conditions in the sense that if $V$ is replaced by $\Omega$ in \eqref{eq:changeenergy}--\eqref{eq:energyfluxterms} then all the terms in \eqref{eq:workterms} and \eqref{eq:energyfluxterms} vanish. 

With regards to the velocity, we consider the impenetrable boundary condition
\begin{equation} \label{eq:BCvel1}
 \v \cdot \nextern = 0 \quad \mbox{on } \p \Omega,
\end{equation}
i.e., the velocity is tangential on $\p \Omega$. We obtain a stress-free boundary condition by imposing the condition
\begin{equation} \label{eq:BCvel2}
 \bbb[P]_{\p \Omega} D(\v) = 0 \quad \mbox{on } \p \Omega,
\end{equation}
where $\bbb[P]_{\p \Omega} = \bbb[I] - \nextern \otimes \nextern \in \R^{d \times d}$ is the projection of $\R^d$ to the tangential space at each point of $\p \Omega$. For the interfaces $\G[i,j]$ we impose the condition
\begin{equation} \label{eq:BCangle}
 \bbb[P]_{\p \Omega} \mextern[i,j] = 0 \quad \mbox{on } \p \Omega \cap \G[i,j],
\end{equation}
which means that the interfaces intersect with $\p \Omega$ with a $90^\circ$ angle. No-flux conditions for both the bulk and the surface surfactant are natural conditions, too, and thanks to \eqref{eq:BCvel1} reduce to
\begin{align}
 \Jc[i] \cdot \nextern &= 0 \quad \mbox{on } \p \Omega \cap \p \O[i], \label{eq:BCsagbulk} \\
 \Jcc[i,j] \cdot \mextern[i,j] &= 0 \quad \mbox{on } \p \Omega \cap \p \G[i,j]. \label{eq:BCsagsurf}
\end{align}

Depending on the application, other boundary conditions may be of relevance. Instead of \eqref{eq:BCvel1} and \eqref{eq:BCvel2} one could consider a Dirichlet condition for $\v$. In the case of in-flow or out-flow \eqref{eq:BCsagbulk} then should read $(\c[i] \v + \Jc[i]) \cdot \nextern = 0$ and similarly for \eqref{eq:BCsagsurf}. For the surfactant a Dirichlet or a Robin condition may be of interest, too, and can easily be expressed in terms of the $\c[i]$ and the $\cc[i,j]$ or, equivalently, in terms of $q$. In all of these cases of non-natural boundary conditions the terms in \eqref{eq:workterms} and \eqref{eq:energyfluxterms} with $V(t) = \Omega$ will not vanish any more, in general, but may be interpreted as working performed by the boundary condition.

\subsection{Summary of the sharp interface model}
\label{sec:sumSIM}

Let us summarise the equations governing the evolution of the multi-phase flow with surfactant. The problem consists in finding a continuous velocity field $\v$, a pressure $p$ and a continuous chemical potential $q$ such that in the domains $\O[i]$ 
\begin{align}
 \div \v &=0, \label{eq:massbal2}\\
 \p_t(\r[i] \v) + \div(\r[i] \v \otimes \v) &= \div \big{(} -p \bbb[I] + 2\et[i] D(\v) \big{)}, \label{eq:mombal2} \\
 \md[\v] \c[i](q) &= \div \big{(} \Mc[i] \nabla q \big{)}, \label{eq:sagbulk2}
\end{align}
on the interfaces $\Gamma^{(i,j)}$
\begin{align}
 \u[i,j] =& \, \v \cdot \n[i,j], \label{eq:intvelcont2} \\
 [ -p \bbb[I] + 2\et[\cdot] D(\v) ]_j^i \n[i,j] =& \, \sigt[i,j](q) \curv[i,j] + \nabg[i,j] \sigt[i,j](q), \label{eq:surfforcebal2} \\
 \md[\v] \cc[i,j](q) + \cc[i,j](q) \nabd[i,j] \v =& \, \nabd[i,j] \big{(} \Mcc[i,j] \nabg[i,j] q \big{)} + \big{[} \Mc[\cdot] \nabla q \big{]}_i^j \cdot \n[i,j], \label{eq:sagsurf2}
\end{align}
and at the triple junctions $T^{(i,j,k)}$
\begin{align}
 \uu[i,j,k] =& \, \bbb[P]_{(\T[i,j,k])^\perp} \v, \label{eq:triplevelcont2} \\
 0 =& \, \Mcc[i,j] \nabg[i,j] q \cdot \m[i,j,k] + \Mcc[j,k] \nabg[j,k] q  \cdot \m[j,k,i] + \Mcc[k,i] \nabg[k,i] q \cdot \m[k,i,j], \label{eq:sagtripmasscond2} \\
 0 =& \, \sigt[i,j](q) \m[i,j,k] + \sigt[j,k](q) \m[j,k,i] + \sigt[k,i](q) \m[k,i,j]. \label{eq:tripleforcebal2}
 \end{align}
These equations then are completed with suitable initial conditions and boundary conditions as discussed in Section \ref{sec:SIMBC}. 

Observe that thanks to \eqref{eq:identtd} the surface surfactant equation \eqref{eq:sagsurf2} can also be written in the following form, which is more convenient for the asymptotic analysis:
\begin{equation} \label{eq:sagsurf2B}
 \ndu[i,j] \cc[i,j](q) + \nabd[i,j] \big{(} \cc[i,j](q) \v \big{)} = \nabd[i,j] \big{(} \Mcc[i,j] \nabg[i,j] q \big{)} + \big{[} \Mc[\cdot] \nabla q \big{]}_i^j \cdot \n[i,j].
\end{equation}

The phase field approach to the surfactant equations will be based on the following distributional form, which can be derived following the lines of \cite{Alt09}:
\begin{equation} \label{distall}
\md[\v] \Big( \sum_{i} \cho[i] \c[i](q) + \sum_{i<j} \delta_{\G[i,j]} \cc[i,j](q) \Big) = - \div \Big( \sum_{i} \cho[i] \Jc[i] + \sum_{i<j} \delta_{\G[i,j]} \Jcc[i,j] \Big).
\end{equation}
Here, $\delta_{\G[i,j]}$ and $\cho[i]$ are the distributions associated with the $\G[i,j]$ and the $\O[i]$, respectively, i.e., 
\[
 \langle \delta_{\G[i,j]}, \phi \rangle = \int_{\G[i,j]} \phi, \qquad \langle \cho[i], \phi \rangle = \int_{\O[i]} \phi, \qquad \forall \phi \in C^\infty_0((0,T) \times \Omega).
\]

\section{Diffuse interface model}
\label{sec:DIM}

The objective is now to derive a phase field model to approximate the free boundary problem that was presented in Section \ref{sec:SIM}. As in \cite{GarLamSti14} we postulate abstract balance equations for phase field variables, mass, momentum, and surfactant and close them within an energetic framework. We postulate a suitable free energy density that approximates the free energy of the sharp interface model. The phase field model for multi-phase flow is based on \cite{AbeGarGru12}, which is extended to multiple phases. 

\subsection{Phase field approach and balance equations}
\label{sec:DIMbaleq}

We begin by introducing a small length scale $\e>0$, the {\em interfacial thickness parameter}, that characterises the length scales of interfacial layers between the different fluids or, more precisely, the different {\em phases} of a fluid domain. It is a fundamental parameter of the approximation, thus we shall use it as an index for all newly defined variables depending on $\e$.
As usual in phase field approaches to multi-phase problems we introduce one phase field variable for each phase (here, the immiscible fluids) that serves to model its presence. Denoting by $\rh[i]$ the mass density of fluid $i$ we define the phase field variables by 
\begin{equation} \label{eq:defphi}
 \ph[i] \coloneqq \frac{\rh[i]}{\r[i]}, \quad i=1, \dots, M.
\end{equation}
As the fluids are immiscible one will expect that $\rh[i] \approx \r[i]$ in the domain of fluid $i$ and $\rh[i] \approx 0$ in the other domains. Only in the thin layers between the fluid domains the fluids are allowed to mix and $\ph[i]$ may take values between zero and one. We assume that there is no excess volume of mixing in these layers so that\footnote{In a small control volume $V$, the masses of the fluids are given by $M^{(i)} = \rh[i] V$. No excess volume of mixing means that $V$ coincides with the sum of the volumes $V^{(i)} = M^{(i)} / \r[i]$ occupied by the same masses of pure fluids, $V = \sum_i V^{(i)}$. Dividing this identity by $V$ yields \eqref{eq:phiconstraint}.}
\begin{equation} \label{eq:phiconstraint}
 \sum_{i=1}^M \ph[i] = 1.
\end{equation}
Introducing the \emph{Gibb's Simplex}
\[
 \Sigma^M \coloneqq \Big{\{} u = (u_1,\dots,u_M)\in \R^M\ :\ \sum_{i=1}^M u_i = 1,\text{ where } 0 \leq u_i \leq 1  \Big{\}},
\]
as well as
\[
 T\Sigma^M \coloneqq \Big{\{} u = (u_1,\dots,u_M)\in \R^M\ :\ \sum_{i=1}^M u_i = 0 \Big{\}},
\]
which can be naturally identified with the tangent space on $\Sigma^M$ at each point, we thus have that $\pha = (\ph[1], \dots, \ph[M]) \in \Sigma^M$. Note that the corners of the Gibb's simplex correspond to the pure fluids as at those points one of the phase field variables equals one and all the others are zero. We write $\be[k] = (\hat{\delta}_{k,l})_{l=1}^M$, $k=1, \dots, M$ for these corners, where $\hat{\delta}_{k,l}$ stands for the Kronecker symbol. For later use we also introduce $\bbb[1] = (1, \dots, 1) \in \R^M$ and note that vectors $u \in T\Sigma^M$ are characterised by $u \cdot \bbb[1] = 0$. 

Denoting by $\vel[i]$ the velocity of mass particles of fluid $i$ the mass balances for the fluids read
\begin{equation} \label{eq:massbalspecies}
 \p_t \rh[i] + \div ( \rh[i] \vel[i] ) = 0.
\end{equation}
In order to describe the motion of the fluid mixture we resort to the \emph{volume averaged velocity} by
\[
 \vv \coloneqq \sum_{i=1}^M \ph[i] \vel[i],
\]
which is solenoidal: Using \eqref{eq:phiconstraint}, \eqref{eq:defphi}, and \eqref{eq:massbalspecies}
\begin{equation} \label{eq:velsolenoidal}
 \div \vv = \p_t \Big{(} \sum_{i=1}^M \ph[i] \Big{)} + \div \Big{(} \sum_{i=1}^M \ph[i] \vel[i] \Big{)} = \sum_{i=1}^M \frac{1}{\r[i]} \big{(} \p_t \rh[i] + \div (\rh[i] \vel[i]) \big{)} = 0.
\end{equation}
As in the previous section (see \eqref{eq:defmatvel}) we define the  \emph{material derivative}
\[
  \md[\vv] w \coloneqq \p_t w + \vv \cdot \nabla w, 
\]
with respect to the velocity field $\vv$. The mass balances \eqref{eq:massbalspecies} yield that
\begin{align} 
 \md[\vv] \ph[i] + \ph[i] \div \vv 
 &= - \div \Jp[i], \label{eq:massbalphi} \\
 \Jp[i] &= \ph[i] (\vel[i] - \vv). \label{eq:identJphi}
\end{align}
Note that, thanks to \eqref{eq:velsolenoidal}, the total mass density
\[
 \rho_\e \coloneqq \sum_{i=1}^M \ph[i] \rh[i],
\]
satisfies the equation
\begin{equation}
 \md[\vv] \rho_\e + \rho_\e \div \vv = -\div \Jbar \ \ \ \text{ with } \ \Jbar = \sum_{i=1}^M \r[i] \Jp[i]. \label{eq:Jbar}
\end{equation}

We now assume that the inertia and the kinetic energy, which are due to the motion of the fluids relative to the gross motion given in terms of $\vv$, is negligible. Thus, rather than formulating momentum balances for the individual velocities $\vel[i]$ we will formulate the conservation of (linear) momentum in terms of $\vv$ and, within an energetic framework presented further below, make assumptions on the fluxes $\Jp[i]$. With a stress tensor $\dimstr_\e$ yet to be determined we postulate
\begin{equation} \label{eq:dimmombal}
 \md[\vv] (\rho_\e \vv) + \rho_\e \vv \, \div \vv 
 = \div \dimstr_\e.
\end{equation}

In order to approximate the (distributional form of the) surfactant equation \eqref{distall} we need to approximate the distributions $\delta_{\G[i,j]}$ and $\cho[i]$ with the help of the phase field variables. Denote by $\del[i,j](\pha,\npha)$ an approximation to $\delta_{\G[i,j]}$, which will be picked later on (see \eqref{eq:def_delta}), and let
\begin{equation} \label{eq:defxi}
 \xi_i(\ph[i]) \coloneqq 
 \begin{cases}
    0 & \quad \text{ if } \ph[i] \leq 0, \\
    1 & \quad \text{ if } \ph[i] \geq 1, \\
    (\ph[i])^2 (3 -2 (\ph[i])) & \quad \text{ else, }
 \end{cases}
\end{equation}
denote an approximation of the characteristic function $\cho[i]$. Recalling that we are studying the case of instantaneous sorption at the phase interfaces, we consider the following regularisation of the surfactant mass balance equation \eqref{distall} for a variable $\qq$:
\begin{multline}
 \md[\vv] \Big( \sum_{i} \xi_i(\ph[i])\c[i](\qq) + \sum_{i<j}\del[i,j](\pha,\nabla\pha) \cc[i,j](\qq) \Big) \\
 + \Big( \sum_{i} \xi_i(\ph[i])\c[i](\qq) + \sum_{i<j}\del[i,j](\pha,\nabla\pha) \cc[i,j](\qq) \Big) \div \vv \\
 + \div \Big (\sum_{i} \xi_i(\ph[i]) \JJc[i] + \sum_{i<j} \del[i,j](\pha,\nabla\pha) \JJcc[i,j] \Big) = 0 \label{regsurfbal}
\end{multline}
with fluxes $\JJc[i]$ and $\JJcc[i,j]$ to be determined later on.  The variable $\qq$ is a diffuse interface approximation of the continuous chemical potential $q$ in the sharp interface model.  In particular, we have an analogous relation $\qq = \g[i]'(\c[i](\qq)) = \gg[i,j]'(\cc[i,j](\qq))$ to \eqref{q}.

\begin{rem}
Here are a few remarks on the above generalisation of Model C in \cite{GarLamSti14}, which is based on the two-phase flow model by \cite{AbeGarGru12} to multiple phases and surfactant fields:
\begin{itemize}
 \item In practice, the hard constraint $\ph[i] \in [0,1]$ often is dropped in favour of a soft one, i.e., values outside of the interval are permitted but energetically expensive. 
 \item We could have dropped the terms with $\div \vv$ in \eqref{eq:massbalphi}, \eqref{eq:dimmombal}, and \eqref{regsurfbal} thanks to \eqref{eq:velsolenoidal}. However, keeping them we get a better idea of pressure contributions to the stress tensor from the thermodynamic analysis below. In particular, we can identify terms associated with the interface that are scaling with $\e^{-1}$, which is beneficial for the subsequent asymptotic analysis. 
 \item Instead of the mass density ratio one could pick different fields for the order parameters $\ph[i]$ such as the $\rh[i]$ or the mass concentrations $\rh[i] / \rho_\e$, see \cite{AbeGarGru12} for a discussion. The essential requirement is that the mass densities $\rh[i]$ and the total mass density $\rho_\e$ can be expressed in terms of the $\ph[i]$.
 \item The expectation is that the phase field variables $\ph[i]$ converge to the $\cho[i]$ as the interfacial thickness converges to zero. The above choice of $\xi_i$ is a $C^1$ function of $\ph[i]$ and satisfies $\xi_i'(p) = 0$ if $p \in \{ 0, 1 \}$, which will enable to recover the sharp interface model as we will see in the asymptotic analysis.
\end{itemize}
\end{rem}

\subsection{Free energy}
\label{sec:choiceFE}

The significance of the small parameter $\e$ is how it features in a Ginzburg-Landau type energy for the phase field variables that serves to approximate the surface energies of the various possible interfaces. Let $\check a : \Sigma^M \times (T \Sigma^M)^d \to [0,\infty)$ be a gradient potential, which is positive ($\check a(\phi,X) >0$ whenever $X \neq 0$) and even and two-homogeneous in the second argument ($\check a(\phi,\eta X) = \eta^2 a(\phi, X)$ for all $\eta \geq 0$), and let $\check w : \Sigma^M \to [0,\infty]$ be a multi-well potential satisfying $\check w(\phi) =0$ if and only if $\phi$ is one of the corners of $\Sigma^M$. Under some more regularity and technical assumptions on $\check a$ and $\check w$, which we skip for brevity, it is shown in \cite{BelBraRie05} that, as $\e \to 0$, 
\[
\int_\Omega \Big{(} \e\check  a(\pha,\npha) + \frac{1}{\e} \check w(\pha) \Big{)} \quad
\to \quad \sum_{i < j } \int_{\G[i,j]} \cgg[i,j](\n[i,j]),
\]
in the sense of a $\Gamma$-limit. The relation between the potential and the surface energies is given by the minimisation problems (see \cite{Ste91,GarNesSto98})
\begin{multline*}
 \cgg[i,j] (\n[i,j]) = \inf_p \Big{\{} 2 \int_{-1}^1 \sqrt{\check w(p)\check  a(p,p'\otimes \n[i,j])} dy \, \Big{|} \\
p: [-1,1] \to \Sigma^M \mbox{ Lipschitz }, \, p(-1) = \be[i], \, p(1) = \be[j] \Big{\}},
\end{multline*}
where $\be[i], \be[j] \in \R^M$ the corners of the Gibb's simplex corresponding to the fluids $i$ and $j$. Note that this formula even holds for some anisotropic surface energies but we here only consider isotropic surface energies. 

For na\"{i}ve choices of $\check a$ and $\check w$, minimisers lie in the interior of $\Sigma^M$ rather than along the edge that connects $e_i$ with $e_j$. In numerical simulations so-called \emph{third phase contributions} then can be observed within the thin interfacial layers \cite{GarNesSto99B}. While they may be considered unphysical the main issue is that they make the recovery of given surface energies $ \cgg[i,j]$ difficult, see \cite{Sti05} for an outline of the problem. But suitable potentials avoid those interfacial third phase contributions (or satisfy the \emph{consistency principle} introduced in \cite{BoyMin14} of reducing to a two-phase system given suitable initial and boundary data). These potentials also enable the approximation of given surface energies $ \cgg[i,j]$, see \cite{GarNesSto99A,GarNesSto99B,BoyLap06,BoyMin14}. During the asymptotic analysis in Section \ref{sec:in_sol} the impact of the choice of such suitable potentials will be clarified. We build up on these works to approximate the energy \eqref{eq:SIMenerg} and consider an energy of the form
\begin{equation} \label{eq:DIMenerg}
 E_\e \coloneqq \int_\Omega e_\e, \quad e_\e \coloneqq \frac{\rho_\e}{2} |\vv|^2 + f(\qq,\pha) + \frac{1}{\e} w(\qq,\pha) + \e a(\qq,\pha,\npha),
\end{equation}
with the contributions
\begin{align}
 a(\qq,\pha,\npha) &\coloneqq \sum_{\substack{i,j=1,\dots,M \\ i < j}} \gg[i,j] (\cc[i,j](\qq)) \a[i,j](\pha,\npha), \label{eq:potw} \\
 w(\qq,\pha) &\coloneqq \sum_{\substack{i,j=1,\dots,M \\ i < j}} \gg[i,j] (\cc[i,j](\qq)) \w[i,j](\pha), \label{eq:pota} \\
 f(\qq,\pha) &\coloneqq \sum_{i=1,\dots,M} \xi_i(\ph[i]) \g[i] (\c[i](\qq)). \nonumber 
\end{align}
See \cite{GarNesSto99A,BoyLap06} for possible choices of the $\a[i,j]$ and the $\w[i,j]$. 

As with the sharp interface model we wish for thermodynamic consistency in the sense of the dissipation of the energy being non-negative. We thus have to ensure that
\[
\md[\vv] e_\e + e_\e \div \vv + \div \bbb[j]_{e_\e} \leq 0,
\]
where the free energy density $e_\e$ is defined in \eqref{eq:DIMenerg} and its flux $\bbb[j]_{e_\e}$ will be defined below. We recall the identities $\sigt[i,j] = \gg[i,j](\cc[i,j](\qq)) - \qq \cc[i,j](\qq)$ from \eqref{eq:defsigt} and, for brevity, define an analogous field for the bulk by
\begin{equation} \label{eq:deflamt}
 \lamt[k](\qq) \coloneqq \g[k](\c[k](\qq)) - \qq \c[k](\qq).
\end{equation}

Using the identities \eqref{eq:dimmombal} and \eqref{eq:Jbar}, a straightforward calculation shows that 
\begin{align}
\md[\vv] \big( \frac{\rho_\e}{2}|\vv|^2 \big) 
&= \vv \cdot \md[\vv] (\rho_\e\vv) - \frac{|\vv|^2}{2} \md[\vv] \rho_\e \nonumber \\
&= \div \Big{(} (\dimstr_\e^\perp + (\vv\otimes\Jbar)^\perp ) \vv - \frac{|\vv|^2}{2} \Jbar \Big{)} \nonumber \\
& \quad - (\dimstr_\e + \vv\otimes\Jbar) \colon \nabla \vv - \rho_\e \frac{|\vv|^2}{2} (\div \vv). \label{eq:mdenergy1}
\end{align}
For the other energy contribution we recall the definition of $\qq$ in \eqref{q} and obtain that 
\begin{align}
\md[\vv] \big{(} f & \, + \frac{1}{\e} w + \e a \big{)} \nonumber \\
=& \, \sum_i \qq (\md[\vv] \c[i]) \xi_i + \g[i] \xi_i' \md[\vv] \ph[i] \nonumber \\
 & \, + \sum_{i<j} \qq (\md[\vv] \cc[i,j]) \e \a[i,j] \nonumber \displaybreak[0] \\
 & \, + \sum_{i<j} \gg[i,j] \sum_k \e \big{(} \p_{\ph[k]} \a[i,j] \md[\vv] \ph[k] + \p_{\nabla \ph[k]} \a[i,j] \cdot \md[\vv] (\nabla \ph[k]) \big{)} \nonumber \\
 & \, + \sum_{i<j} \qq (\md[\vv] \cc[i,j]) \einvs \w[i,j] + \gg[i,j] \sum_k \einvs \p_{\ph[k]} \w[i,j] \md[\vv] \ph[k] \nonumber \displaybreak[0] \\
=& \, \sum_i \qq \md[\vv] (\c[i] \xi_i) + (\g[i] - \c[i] \qq ) \xi_i' \md[\vv] \ph[i] \nonumber \\
 & \, + \sum_{i<j} \big{(} \gg[i,j] - \qq \cc[i,j] \big{)} \sum_k \einvs \p_{\ph[k]} \w[i,j] \md[\vv] \ph[k] \displaybreak[0] \nonumber \\
 & \, + \sum_{i<j} \big{(} \gg[i,j] - \qq \cc[i,j] \big{)} \sum_k \e \big{(} \p_{\ph[k]} \a[i,j] \md[\vv] \ph[k] + \p_{\nabla \ph[k]} \a[i,j] \cdot \md[\vv] (\nabla \ph[k]) \big{)} \nonumber \\
 & \, + \sum_{i<j} \qq \md[\vv] \big{(} \cc[i,j] (\einvs \w[i,j] + \e \a[i,j]) \big{)}. \label{eq:mdenergy2}
\end{align}
Using \eqref{eq:defsigt} and the identity
\[
 \md[\vv](\nabla\ph[k]) = \nabla \md[\vv]\ph[k] - (\nabla\vv)^\perp \nabla\ph[k],
\]
and setting 
\begin{equation} \label{eq:def_delta}
 \del[i,j](\pha,\nabla\pha) \coloneqq \e \a[i,j](\pha,\nabla\pha) + \frac{1}{\e} \w[i,j](\pha),
\end{equation}
we obtain
\begin{align*}
  & \, \big{(} \gg[i,j] - \qq \cc[i,j] \big{)} \e \p_{\nabla \ph[k]} \a[i,j] \cdot \md[\vv] (\nabla \ph[k]) \big{)} \\
 =& \, \sigt[i,j] \p_{\nabla \ph[k]} \del[i,j] \cdot \big{(} \nabla \md[\vv] \ph[k] - (\nabla \vv)^\perp \nabla \ph[k] \big{)} \\ 
 =& \, \div \big{(} \sigt[i,j] \p_{\nabla \ph[k]} \del[i,j] \md[\vv] \ph[k] \big{)} - \div \big{(} \sigt[i,j] \p_{\nabla \ph[k]} \del[i,j] \big{)} \md[\vv] \ph[k] \\
 & \, - \sigt[i,j] \nabla \ph[k] \otimes \p_{\nabla \ph[k]} \del[i,j] : \nabla \vv.
\end{align*}
Therefore, continuing with \eqref{eq:mdenergy2} and using \eqref{eq:defsigt} and \eqref{eq:deflamt}
\begin{align}
\md[\vv] \big{(} f & \, + \frac{1}{\e} w + \e a \big{)} \nonumber \\
=& \, \qq \, \md[\vv] \Big{(} \sum_i \xi_i \c[i] + \sum_{i<j} \del[i,j] \cc[i,j] \Big{)} \nonumber \\
 & \, + \sum_i \lamt[i] \xi_i' \md[\vv] \ph[i] \displaybreak[0] \nonumber \\
 & \, + \sum_{i<j} \sum_k \Big{(} \sigt[i,j] \p_{\ph[k]} \del[i,j] - \div \big{(} \sigt[i,j] \p_{\nabla \ph[k]} \del[i,j] \big{)} \Big{)} \md[\vv] \ph[k]  \nonumber \\
 & \, + \sum_{i<j} \sum_k \div \big{(} \sigt[i,j] \p_{\nabla \ph[k]} \del[i,j] \md[\vv] \ph[k] \big{)} - \sigt[i,j] \nabla \ph[k] \otimes \p_{\nabla \ph[k]} \del[i,j] : \nabla \vv, \displaybreak[0] \nonumber \\
\intertext{that, when inserting the balance equations \eqref{regsurfbal} and \eqref{eq:massbalphi}, yields}
=& \, -\qq \, \div \Big{(} \sum_{i} \xi_i \JJc[i] + \sum_{i<j} \del[i,j] \JJcc[i,j] \Big{)} \nonumber \\
 & \, - \qq \Big{(} \sum_i \xi_i \c[i] + \sum_{i<j} \del[i,j] \cc[i,j] \Big{)} \div \vv \nonumber \\
 & \, - \sum_i \lamt[i] \xi_i' \, (\div \Jp[i] + \ph[i] \div \vv) \displaybreak[0] \nonumber \\
 & \, - \sum_{i<j} \sum_k \Big{(} \sigt[i,j] \p_{\ph[k]} \del[i,j] - \div \big{(} \sigt[i,j] \p_{\nabla \ph[k]} \del[i,j] \big{)} \Big{)} (\div \Jp[k] + \ph[k] \div \vv) \nonumber \\
 & \, + \sum_{i<j} \sum_k \div \big{(} \sigt[i,j] \p_{\nabla \ph[k]} \del[i,j] \md[\vv] \ph[k] \big{)} - \sigt[i,j] \nabla \ph[k] \otimes \p_{\nabla \ph[k]} \del[i,j] : \nabla \vv \displaybreak[0] \nonumber \\
=& \, \div \Big{[} -\qq \big{(} \sum_{i} \xi_i \JJc[i] + \sum_{i<j} \del[i,j] \JJcc[i,j] \big{)} \nonumber \\
 & \phantom{\div \Big{(}} - \sum_k \Big{(} \lamt[k] \xi_k' + \sum_{i<j} \big{(} \sigt[i,j] \p_{\ph[k]} \del[i,j] - \div \big{(} \sigt[i,j] \p_{\nabla \ph[k]} \del[i,j] \big{)} \big{)} \Big{)} \Jp[k] \nonumber \\
 & \phantom{\div \Big{(}} + \sum_k \sum_{i<j} \big{(} \sigt[i,j] \p_{\nabla \ph[k]} \del[i,j] \md[\vv] \ph[k] \big{)} \Big{]} \displaybreak[0] \nonumber \\
 & \, + \sum_{i} \xi_i \nabla \qq \cdot \JJc[i] + \sum_{i<j} \del[i,j] \nabla \qq \cdot \JJcc[i,j] \displaybreak[0] \nonumber \\
 & \, + \sum_k \nabla \Big{(} \lamt[k] \xi_k' + \sum_{i<j} \big{(} \sigt[i,j] \p_{\ph[k]} \del[i,j] - \div \big{(} \sigt[i,j] \p_{\nabla \ph[k]} \del[i,j] \big{)} \big{)} \Big{)} \cdot \Jp[k] \displaybreak[0] \nonumber \\
 & \, - \Big{[} \sum_i \xi_i \qq \c[i] + \sum_{i<j} \del[i,j] \qq \cc[i,j] \Big{]} \div \vv \displaybreak[0] \nonumber \\
 & \, - \Big{[} \sum_k \Big{(} \lamt[k] \ph[k] \xi_k' + \ph[k] \sum_{i<j} \big{(} \sigt[i,j] \p_{\ph[k]} \del[i,j] - \div \big{(} \sigt[i,j] \p_{\nabla \ph[k]} \del[i,j] \big{)} \big{)} \Big{)} \Big{]} \div \vv \nonumber \\
 & \, - \sum_k \sum_{i<j} \sigt[i,j] \nabla \ph[k] \otimes \p_{\nabla \ph[k]} \del[i,j] : \nabla \vv. \label{eq:mdenergy3}
\end{align}
Defining
\begin{align*}
 \bbb[j]_{e_\e} 
\coloneqq & \, - (\dimstr_\e^\perp + (\vv\otimes\Jbar)^\perp ) \vv + \frac{|\vv|^2}{2} \Jbar \\
 & \, + \qq \big{(} \sum_{i} \xi_i \JJc[i] + \sum_{i<j} \del[i,j] \JJcc[i,j] \big{)} \displaybreak[0] \nonumber \\
 & \, + \sum_k \Big{(} \lamt[k] \xi_k' + \sum_{i<j} \big{(} \sigt[i,j] \p_{\ph[k]} \del[i,j] - \div \big{(} \sigt[i,j] \p_{\nabla \ph[k]} \del[i,j] \big{)} \big{)} \Big{)} \Jp[k] \nonumber \\
 & \, - \sum_k \sum_{i<j} \big{(} \sigt[i,j] \p_{\nabla \ph[k]} \del[i,j] \md[\vv] \ph[k] \big{)},
\end{align*}
we obtain that from \eqref{eq:mdenergy1} and \eqref{eq:mdenergy3} that
\begin{align}
 \md[\vv] e_\e & \, + e_\e \div \vv + \div \bbb[j]_{e_\e} \nonumber \\
=& \, \sum_{i} \xi_i \nabla \qq \cdot \JJc[i] + \sum_{i<j} \del[i,j] \nabla \qq \cdot \JJcc[i,j] \nonumber \\
 & \, + \sum_k \nabla \Big{(} \lamt[k] \xi_k' + \sum_{i<j} \big{(} \sigt[i,j] \p_{\ph[k]} \del[i,j] - \div \big{(} \sigt[i,j] \p_{\nabla \ph[k]} \del[i,j] \big{)} \big{)} \Big{)} \cdot \Jp[k] \displaybreak[0] \nonumber \\
 & \, + \Big{[} \sum_i \xi_i \lamt[i] + \sum_{i<j} \del[i,j] \sigt[i,j] \Big{]} \div \vv \displaybreak[0] \nonumber \\
 & \, - \Big{[} \sum_k \Big{(} \lamt[k] \ph[k] \xi_k' + \ph[k] \sum_{i<j} \big{(} \sigt[i,j] \p_{\ph[k]} \del[i,j] - \div \big{(} \sigt[i,j] \p_{\nabla \ph[k]} \del[i,j] \big{)} \big{)} \Big{)} \Big{]} \div \vv \displaybreak[0] \nonumber \\
 & \, - \Big{(} \dimstr_\e + \vv\otimes\Jbar + \sum_k \sum_{i<j} \sigt[i,j] \nabla \ph[k] \otimes \p_{\nabla \ph[k]} \del[i,j] \Big{)} \colon \nabla \vv. 
 \label{eq:mdenergy4}
\end{align}

\subsection{Constitutive assumptions and boundary conditions}
\label{sec:DIMass}

The calculations resulting in \eqref{eq:mdenergy4} motivate to make the following assumptions that ensure non-negative energy dissipation:
\begin{align*}
 \JJc[i] \, &\coloneqq -\Mc[i]\nabla \qq, \\
 \JJcc[i,j] \, &\coloneqq -\Mcc[i,j]\nabla \qq,
\end{align*}
with the mobilities $\Mc[i]$ and $\Mcc[i,j]$ as in \eqref{eq:sagbulkflux}, \eqref{eq:sagsurfflux},
\begin{align*}
 \Jp[k] \, &\coloneqq -\sum_{l=1}^M \L[k,l] \nabla\mb[l], \quad \mbox{where } \\
 \mb[l] \, &\coloneqq \lamt[l] \xi_l' + \sum_{i<j} \big{(} \sigt[i,j] \p_{\ph[l]} \del[i,j] - \div \big{(} \sigt[i,j] \p_{\nabla \ph[l]} \del[i,j] \big{)} \big{)},
\end{align*}
with mobilities $\mathcal{L}^{(k,l)}$ that may depend on $\pha$ and $\qq$, form a symmetric positive semi-definite matrix, and satisfy 
\begin{equation} \label{eq:propL}
 \sum_{k=1}^M \L[k,l](\pha,\qq) = 0 \quad \forall \pha \in \Sigma^M, \qq \in \R,
\end{equation}
which ensures that \eqref{eq:phiconstraint} is fulfilled during the evolution, and finally
\begin{align*}
 \dimstr_\e \coloneqq \, & -\tpe \bbb[I] + 2\eta(\pha) D(\vv) - \vv\otimes\Jbar \\
 &- \sum_k \sum_{i<j} \big{(} \sigt[i,j] \nabla \ph[k] \otimes \p_{\nabla \ph[k]} \del[i,j] \big{)} + \Big{(} \sum_k \big{(} \xi_k \lamt[k] - \mb[k] \ph[k] \big{)} + \sum_{i<j} \del[i,j] \sigt[i,j] \Big{)} \bbb[I],
\end{align*}
with a pressure $\tpe$, and where $\eta(\pha)$ is a non-negative smooth interpolation function between the viscosities of the pure fluids, i.e., $\eta(\ph[1], \dots, \ph[M]) = \et[i]$ if $\ph[i]=1$ (and then $\ph[j] = 0$ for $j \neq i$ by \eqref{eq:phiconstraint}). We can absorb some of terms multiplying the identity tensor $\bbb[I]$ into the pressure but keep those terms that, in the interfacial regions, are required to identify the terms to leading order in $\e$. Setting
\[
 \pp \coloneqq \tpe - \sum_k \big{(} \mb[k] \ph[k] - \xi_k \lamt[k] \big{)}, 
\]
we obtain 
\begin{align}
 \dimstr_\e = & \, -\pp \bbb[I] + 2\eta(\pha) D(\vv) - \vv\otimes\Jbar \nonumber \\
 & + \sum_{i<j} \sigt[i,j] \Big{(} \del[i,j] - \sum_k \nabla \ph[k] \otimes \p_{\nabla \ph[k]} \del[i,j] \Big{)} \bbb[I], \label{eq:DIMdefT}
\end{align}
where we also recall the definition of $\Jbar$ from \eqref{eq:Jbar}.

Natural boundary conditions on $\p \Omega$ arise from assuming a closed system so there are no mass and energy fluxes into or out of the domain. For the fluid flow they read as in the sharp interface model \eqref{eq:BCvel1}, \eqref{eq:BCvel2}, 
\begin{align}
 0 &= \vv \cdot \nextern, \label{eq:PFMBCnoflux} \\
 0 &= \bbb[P]_{\p \Omega} D(\vv) \quad \mbox{with } \bbb[P]_{\p \Omega} = \bbb[I] - \nextern \otimes \nextern. \label{eq:PFMBCnoforce}
\end{align}
For the phase fields natural conditions are 
\begin{align*}
 0 &= \nabla \mb[l] \cdot \nextern, \\
 0 &= \p_{\nabla \ph[k]} \del[i,j] \cdot \nextern,
\end{align*}
for all $k,l = 1, \dots, M$, where the first one ensures a no-flux condition for the $\Jp[k]$ and, as we shall see in the asymptotic analysis, the second one is related to angles between the interface $\G[i,j]$ and the external boundary $\p \Omega$. In order to guarantee a no-flux boundary condition for the surfactant mass one may assume that
\begin{equation} \label{eq:PFMBCqnoflux}
 0 = \nabla \qq \cdot \nextern.
\end{equation}

\subsection{Summary of the diffuse interface model}
\label{sec:sumDIM}

Summarising the phase field equations we have a Cahn-Hilliard type system for the phase fields of the form
\begin{align}
\md[\vv] \ph[k] &= - \div \Jp[k], \label{eq:PFMphi1} \\
\Jp[k] &= - \sum_{l} \L[k,l] \nabla \mb[l], \label{eq:PFMphi2} \displaybreak[0] \\
\mb[l] &= \lamt[l] \xi_l' + \sum_{i<j} \big{(} \sigt[i,j] \p_{\ph[l]} \del[i,j] - \div \big{(} \sigt[i,j] \p_{\nabla \ph[l]} \del[i,j] \big{)} \big{)}, \label{eq:PFMmu}
\end{align}
for $k,l=1,\dots,M$. It is coupled to an equation for the surfactant
\begin{align}
 \md[\vv] \Big( \sum_{i} \xi_i \c[i](\qq) + \sum_{i<j} \del[i,j] \cc[i,j](\qq) \Big) = - \div \Jq, \label{eq:PFMsag1} \\
 \Jq = - \Big{(} \sum_{i} \xi_i \Mc[i] \nabla \qq + \sum_{i<j} \del[i,j] \Mcc[i,j] \nabla \qq \Big{)} , \label{eq:PFMsag2} 
\end{align}
while the fluid flow is subject to the Navier-Stokes system
\begin{align}
 \div\vv &= 0, \label{eq:PFMmassbal} \\ 
 \md[\vv] (\rho_\e\vv) &= \div \Big( -p\bbb[I] + 2\eta(\pha) D(\vv) - \vv \otimes \sum_{k} \r[k] \Jp[k] \Big{)} \nonumber \\
 &\quad + \div \Big{(} \sum_{i<j} \sigt[i,j] \Big{(} \del[i,j] - \sum_{k} \nabla \ph[k] \otimes \p_{\nabla \ph[k]} \del[i,j] \Big{)} \bbb[I] \Big{)}.  \label{eq:PFMmombal}
\end{align}

For completion of the problem, boundary conditions as discussed in Section \ref{sec:DIMass} and suitable initial conditions have to be imposed.

We may reform the capilliary forcing in the Navier-Stokes system. Starting with
\begin{align*}
    &\sum_{k}\mb[k]\nabla\ph[k] \\
    &=\sum_k \sum_{i<j} \Big(- \div(\sigt[i,j]\p_{\nabla \ph[k]} \del[i,j])\nabla\ph[k] + \sigt[i,j] \p_{ \ph[k]} \del[i,j]\nabla\ph[k]\Big) + \lambda_k \xi_k'\nabla\ph[k] \displaybreak[0] \\
    &=\sum_{i<j} \Big(- \sum_k \div(\sigt[i,j] \p_{\nabla \ph[k]} \del[i,j])\nabla\ph[k] + \sum_k \sigt[i,j] \p_{ \ph[k]} \del[i,j]\nabla\ph[k]\Big) + \sum_k\lambda_k\nabla \xi_k \displaybreak[0] \\
    &=\sum_{i<j}  \Big(\div(- \sigt[i,j]\sum_k\nabla \ph[k]\otimes\p_{\nabla \ph[k]} \del[i,j] ) + \sigt[i,j]\nabla \del[i,j] \Big) + \sum_k\lambda_k\nabla \xi_k \\
    &= \div(\sum_{i<j}\sigt[i,j]( \del[i,j]\bbb[I] - \sum_k \nabla \ph[k]\otimes\p_{\nabla \ph[k]} \del[i,j]))- \sum_{i<j} \del[i,j]\nabla\sigt[i,j] + \sum_k\lambda_k\nabla \xi_k,
\end{align*}
and rearranging we find that
\[
 \div(\sum_{i<j}\sigt[i,j]( \del[i,j]\bbb[I] - \sum_k \p_{\nabla \ph[k]} \del[i,j]\otimes \nabla \ph[k])) = \sum_{k}\mb[k]\nabla\ph[k] + \sum_{i<j} \del[i,j]\nabla\sigt[i,j] - \sum_k\lambda_k\nabla \xi_k, 
\]
which can be substituted into \eqref{eq:PFMmombal}.

\subsection{Specific Example}
\label{sec:BLM}

For some numerical simulations, the results of which are presented in Section \ref{sec:numsim}, we pick the model in \cite{BoyLap06, BoyMin14} for $M=3$ phases with further choice of the mobility matrix \eqref{eq:PFMphi2}. More precisely, $w$ is of the form \eqref{eq:potw} with 
\begin{multline} \label{eq:BLMw}
\w[i,j](\pha) = 12\Big((\ph[i])^2(\ph[j])^2 + \sum_{k\neq i,j} \big( \ph[j]\ph[k](\ph[i])^2+\ph[i]\ph[k](\ph[j])^2-\ph[i]\ph[j](\ph[k])^2 \big) \Big)\\
+ 4 \Lambda \sum_{k\neq i,j} (\ph[i])^2(\ph[j])^2(\ph[k])^2,
\end{multline}
where the sixth order polynomial with a sufficiently large $\Lambda > 0$ serves to prevent the leaking of third phase contributions between two other phases outside of the triple junction regions (see the discussion in Section \ref{sec:choiceFE}). The gradient potential $a$ is of the form \eqref{eq:pota} with  
\begin{equation}\label{eq:BLMa}
 \a[i,j](\nabla\pha) = \frac{3}{8} \big( |\nabla\ph[i]|^2 + |\nabla\ph[j]|^2 - \sum_{k\neq i,j} |\nabla \ph[k]|^2 \big).
\end{equation}
For the mobility matrix $\mmm[L]$ in \eqref{eq:PFMphi2} we choose a $\qq$ dependent matrix defined as follows:
\begin{equation}\label{eq:BLML}
 \mmm[L]^{(k,l)}(\qq) = \begin{cases}
                    -\frac{M_c \bar{S}(\qq)}{3S_k(\qq)S_l(\qq)}, & \text{ for } l\neq k, \\
                    \sum_{i\neq l}\frac{M_c \bar{S}(\qq)}{3S_i(\qq)S_l(\qq)}, & \text{ for } k=l.
                   \end{cases}
\end{equation}
Here, $S_k(\qq) = \sigt[i,k](\qq) + \sigt[j,k](\qq)-\sigt[i,j](\qq) = - \tilde{S}^{(i,j,k)}(\qq)$ (see \eqref{eq:defWettcoeff} for the wetting coefficients), and their harmonic average is $\bar{S}= \sum_{i=1}^3 \frac{3}{S_i(\qq)}$. Finally, we take a constant mobility parameter $M_c$. 

In the absense of fluid flow we obtain the following Cahn-Hilliard system: For $i=1,2,3$
\begin{align}\label{eq:BLMequation}
\p_t \ph[i] &=\div\Big(\frac{M_c}{S_i(\qq)} \nabla \mb[i]\Big),\\
\mb[i] &= -\frac{3}{4} \e S_i(\qq) \Delta \ph[i] + \frac{4\bar{S}}{\e}\mathcal{D}_i w(\qq,\mathbf{\pha}), \label{eq:BLMequation2}
\intertext{where }
\mathcal{D}_i w(\qq,\mathbf{\pha}) &= \sum_{j\neq i} \frac{1}{S_j(\qq)} \Big( \p_{\ph[i]} w(\qq,\mathbf{\pha}) - \p_{\ph[j]} w(\qq,\mathbf{\pha}) \Big). \nonumber
\end{align}

\section{Asymptotic Analysis}
\label{sec:asymp}

By matching suitable asymptotic expansions of solutions we show in this section that the formal asymptotic limit of the phase field model presented in Section \ref{sec:sumDIM} is the free boundary problem presented in Section \ref{sec:sumSIM}. The situation in the phases and along the interface layers reduces to the two-phase case. Its asymptotic analysis is presented in \cite{GarLamSti14} in great detail. We still present many details, for the notation is quite different under reformulation with multiple phase fields and, more importantly, we subsequently will require some of the findings to deal with the triple junctions. For the latter, techniques presented in \cite{BroGarSto98,BroRei93,GarNesSto98} are used and further developed to treat the surfactant equation. The case $d=2$, in which the triple junctions are points, is investigated first. Building up on this, triple lines and quadruple points in the case $d=3$ are then considered. 

\subsection{Setting and assumptions}
\label{sec:asympass}

Let $\{ \pha, \Jpa, \mba, \vv, \pp, \qq, \Jq \}_{\e>0}$ denote a family of solutions to \eqref{eq:PFMphi1}--\eqref{eq:PFMmombal}. We make some \emph{assumptions} on the model and the solutions, which are sketched here and further detailed and clarified during the following analysis:

\renewcommand{\theenumi}{A\arabic{enumi}}
\begin{enumerate}
 \item \label{ass:AA1} 
 We are interested in the solution regime where interfacial layers of thickness $\sim \e$ have emerged between the domains in which the phase field is close to one of the minimisers of the multi-well potential $w(\qq,\pha)$. That is, these \emph{phases} are where $\pha \approx \be[i]$ for some $i \in \{ 1, \dots, M \}$ and thus, notionally, the domain is occupied by fluid $i$. 
 \item \label{ass:AA2} 
 The potentials $\a[i,j]$ and $\w[i,j]$ are such that no third-phase contributions appear along the interface layers. See Section \ref{sec:choiceFE} before \eqref{eq:DIMenerg} for a brief discussion and references. The clear meaning of this assumption and its consequences are discussed around equation \eqref{eq:inprofile} below.
 \item \label{ass:AA3} 
 The potentials $\a[i,j]$ and $\w[i,j]$ furthermore are such that the equation for $\qq$, \eqref{eq:PFMsag1} with \eqref{eq:PFMsag2}, is non-degenerate in the triple junctions where any three interfacial layers meet. In particular, close to the triple junction the vector of the phase fields is away from the corners of the Gibb's simplex. In the case $d=3$ this property is assumed to extend to the quadruple points. Around equations \eqref{eq:Q0const} and \eqref{eq:triH2} this assumption is discussed and exploited. 
 \item \label{ass:AA4} The mobilities of the phase fields are of the form
 \[
  \L[k,l](\pha,\qq) = \L[k,l]_0(\pha,\qq) + \e \L[k,l]_1(\pha,\qq),
 \]
 where both the $\L[k,l]_0$ and the $\L[k,l]_1$ form symmetric matrices satisfying \eqref{eq:propL}. Moreover, $\L[k,l](\be[j],\qq) = 0$ for all $j \in \{1, \dots, M\}$ and $\qq \in \R$ but if $\check{\varphi}_\e \in \Sigma^M \backslash \{ \be[j] \}_j$ then the kernel of $\{ \L[k,l]_0(\check{\varphi}_\e,\qq) \}_{k,l}$ is the span of $\bbb[1] = (1, \dots, 1) \in \R^M$. In turn, the matrix $\{ \L[k,l]_1(\pha,\qq) \}_{k,l}$ is non-degenerate for all $(\pha,\qq) \in \Sigma^M \times \R$ in the sense that its kernel is only the span of $\bbb[1]$. 
\end{enumerate}

\subsection{Outer expansions and solutions}
\label{sec:out_expsol}

In points $(x,t)$ in the phases away from the interface layers we consider expansions of the form
\[
 \zeta_\e(x,t) = \zeta_0(x,t) + \e \zeta_1(x,t) + \e^2 \zeta_2(x,t) + \dots,
\]
for all fields $\ph[k]$, $\mb[l]$, $\vv$, $\pp$, and $q_\e$, and also for the fluxes $\Jp[k]$. The flux $\Jq$ contains a term scaling with $\e^{-1}$ whence we assume that it can be expanded in the form
\[
 \Jq = \e^{-1} \Jqe[-1] + \e^{0} \Jqe[0] + \dots\ .
\]
These expansions are plugged into the phase field equations \eqref{eq:PFMphi1}--\eqref{eq:PFMmombal} and all non-linearities are Taylor-expanded.

From \eqref{eq:PFMmassbal} we obtain to leading order $0$ that
\begin{equation} \label{eq:outR1}
 \div \vve[0] = 0.
\end{equation}
Equation \eqref{eq:PFMmu} yields to leading order $-1$ that
\[
 0 = \sum_{i<j} \sigt[i,j](\qqe[0]) \p_{\pha} \w[i,j] (\phae[0]).
\]
As we are in a phase by assumption this implies that $\phae[0]$ is one of the corners of the Gibb's simplex, $\phae[0] = \be[m]$ for some $m \in \{ 1, \dots, M \}$. To the next order $0$ we obtain that
\begin{equation} \label{eq:outH1}
 \mbae[0] = \sum_{i<j} \sigt[i,j](\qqe[0]) \p_{\pha \pha} \w[i,j] (\phae[0]) \phae[1],
\end{equation}
where we used that $\xi_k'(\phe[k,0]) = 0$ (thanks to \eqref{eq:defxi}). 

Considering \eqref{eq:PFMphi1}, \eqref{eq:PFMphi2} to leading order $0$ yields 
\[ 
 0 = - \div \Jpe[k,0], \quad \Jpe[k,0] = - \sum_l \L[k,l]_0(\phae[0],\qqe[0]) \nabla \mbe[l,0]. 
\]
But as $\phae[0] = \be[m]$ we have that $\L[k,l]_0(\phae[0],\qqe[0]) = 0$ so that 
\begin{equation} \label{eq:outH2}
 \Jpe[k,0] = 0.
\end{equation}
Moreover, $\p_{q} \L[k,l]_0(\phae[0],\qqe[0]) = 0$ so that \eqref{eq:PFMphi1}, \eqref{eq:PFMphi2} to the next order read 
\begin{align*}
 & \p_t \phe[k,1] + \vve[0] \cdot \nabla \phe[k,1] = - \nabla \cdot \Jpe[k,1], \\
 & \Jpe[k,1] = - \sum_l (\p_{\pha} \L[k,l]_0(\phae[0],\qqe[0]) \cdot \phae[1] + \L[k,l]_1(\phae[0],\qqe[0])) \nabla \mbe[l,0].
\end{align*}
Inserting \eqref{eq:outH1} this becomes a parabolic problem for $\phae[1]$ that allows for the solution $\phae[1] = 0$. Whether this is the unique solution will depend on the boundary conditions both on the external boundary of the domain as well as the free boundaries. However, we do not need any specific knowledge of these solutions for our asymptotic analysis.

As $\w[i,j](\phae[0]) = 0$ there are no terms to order $-1$ in the momentum equation \eqref{eq:PFMmombal}. To order $0$ it yields that
\begin{equation} \label{eq:outR2}
 \p_t (\r[m] \vve[0]) + (\vve[0] \cdot \nabla) (\r[m] \vve[0]) = \div \big{(} - \ppe[0] \bbb[I] + 2 \et[m] D (\vve[0]) \big{)},
\end{equation}
where we used that $\p_{\pha} \w[i,j] (\phae[0]) = 0$ and \eqref{eq:outH2}. 

Finally, recalling \eqref{eq:def_delta}, using that $\w[i,j](\phae[0]) = 0$ and $\p_{\pha} \w[i,j](\phae[0]) = 0$ in \eqref{eq:PFMsag2}, and using \eqref{eq:defxi} we see that
\begin{align}
 \Jqe[-1] & = - \sum_{i<j} \Mcc[i,j] \w[i,j](\phae[0]) \nabla \qqe[0] = 0, \label{eq:outH4} \\
 \Jqe[0] & = - \Big{(} \sum_i \Mc[i] \xi_i(\phe[i,0]) \nabla \qqe[0] + \sum_{i<j} \Mcc[i,j] \big{(} \p_{\pha} \w[i,j](\phae[0]) \cdot \phae[1] \nabla \qqe[0] + \w[i,j](\phae[0]) \nabla \qqe[1] \big{)} \Big{)} \nonumber \\
 & = - \Mc[m] \nabla \qqe[0]. \label{eq:outH5} 
\end{align}

The same arguments apply to the left-hand side of \eqref{eq:PFMsag1} so that, to order $0$, it reads
\[
 \p_t \c[m](\qqe[0]) + \vve[0] \cdot \nabla \c[m](\qqe[0]) = - \div \Jqe[0] = \div \big{(} \Mc[m] \nabla \qqe[0] \big{)}.
\]
With this equation and \eqref{eq:outR1} and \eqref{eq:outR2} we have recovered the bulk equations \eqref{eq:massbal2}--\eqref{eq:sagbulk2} of the sharp interface model.

\subsection{Inner expansions and matching conditions}
\label{sec:in_expmat}

Consider now an interfacial layer between two domains where $\phae[0] \approx \be[\mminus]$ and $\phae[0] \approx \be[\mplus]$, respectively, for two phase indices $\mminus < \mplus$. For simplicity, we restrict the analysis to the two-dimensional case, $d=2$. However, the final results consisting of \eqref{eq:inR1}, \eqref{eq:inR2}, and \eqref{eq:inR3} can also be retrieved in the higher dimensional case by following exactly the line of argument below. We refer to \cite{RaeVoi06} for the techniques that are required to do so. 

We use the limiting curve of the layer, which belongs to $\G[\mminus,\mplus]$, in order to introduce new coordinates. By $s$ we denote a tangential coordinate along $\G[\mminus,\mplus](t)$ such that, for $t$ given, an arc-length parametrisation is obtained, which is denoted by $\param[\mminus,\mplus](s,t)$. Then $\tang[\mminus,\mplus](\param[\mminus,\mplus](s,t),t) = \p_s \param[\mminus,\mplus](s,t)$ is a unit tangent vector field to $\G[\mminus,\mplus]$. We assume the orientation of $s$ to be such that in $x = \param[\mminus,\mplus](s,t)$
\begin{equation} \label{eq:Frenet}
 \frac{d}{ds} \tang[\mminus,\mplus](x,t) = \curv[\mminus,\mplus](x,t) =: \scurv[\mminus,\mplus](x,t) \n[\mminus,\mplus](x,t),
\end{equation} 
where we introduced the scalar mean curvature $\scurv[\mminus,\mplus]$ of $\G[\mminus,\mplus]$. Then
\[
 \frac{d}{ds} \n[\mminus,\mplus](x,t) = - \scurv[\mminus,\mplus](x,t) \tang[\mminus,\mplus](x,t).
\] 
For any surface resident field $r(t) : \G[\mminus,\mplus](t) \to \R$, written as $R(s,t) = r(x,t)$, $x = \param[\mminus,\mplus](s,t)$, in these new coordinates, we note the following identity (for instance, see \cite{StiDiss06} for a derivation):
\begin{equation}
 \label{eq:identnd}
 \p_t R(s,t) - \p_t \param[\mminus,\mplus](s,t) \p_s R(s,t) = \ndu[\mminus,\mplus] r(x,t), 
\end{equation}
where we recall the notation \eqref{eq:defntd} for the normal time derivative. A further coordinate in direction $\n[\mminus,\mplus]$ is denoted by $z$, which is the signed distance to $\G[\mminus,\mplus](t)$ divided by $\e$, i.e., positive on the side of $\O[\mplus](t)$ and negative on the side of $\O[\mminus](t)$. 

As before, expansions of the solutions fields are plugged into the equations of the phase field model. But this time the expansions are of the form
\begin{align*}
 \zeta_\e(x,t) & = Z_0(s,z,t) + \e Z_1(s,z,t) + \e^2 Z_2(s,z,t) + \dots, \\
 \Jp[k](x,t) & = \e^{-1} \JpE[k,-1](s,z,t) + \e^{0} \JpE[k,0](s,z,t) + \e^{1} \JpE[k,1](s,z,t) + \dots, \\
 \Jq(x,t) & = \e^{-2} \JqE[-2](s,z,t) + \e^{-1} \JqE[-1](s,z,t) + \e^{0} \JqE[0](s,z,t) + \dots,
\end{align*}
for inner variables $Z \in \{ \Phi, M, Q, \bbb[V], P\}$ corresponding to $\{\pha, \mba, \qq, \vv, \pp\}$ 
in points $(x,t)$ close to $\G[\mminus,\mplus](t)$ where the distance function, which is required to define the coordinate $z$, is well-defined. The tangential coordinate $s$ for such a point $x$ is such that $\param[\mminus,\mplus](s,t)$ is the closest point to $x$ on $\G[\mminus,\mplus]$. The differential operators read as follows in the new coordinates \cite{GarLamSti14}:
\begin{align*}
\p_{t} \zeta(x,t) & = - \e^{-1} \u[\mminus,\mplus] \p_{z} Z(s,z,t) + \p_t Z(s,z,t) - \p_t \param[\mminus,\mplus] \p_s Z(s,z,t) + \mmm[O](\e), \\
\nabla \zeta(x,t) & = \e^{-1} \p_{z} Z(s,z,t) \n[\mminus,\mplus]  + (1 + \e \scurv[\mminus,\mplus]) \p_s Z(s,z,t) \tang[\mminus,\mplus] + \mmm[O](\e^2). 
\end{align*}
Here and in the following all interface resident fields such as $\u[\mminus,\mplus]$, $\n[\mminus,\mplus]$, and $\tang[\mminus,\mplus]$ are evaluated in $\param[\mminus,\mplus](s,t) \in \G[\mminus,\mplus](t)$. 

Requiring inner and outer expansions to match leads to the following \emph{matching conditions} \cite{GarSti06}: 
As $z \to \pm \infty$,
\begin{align}
Z_{0}(s,z,t) & \sim \zeta_{0}^{\pm}, \label{MC0} \\
\p_{z} Z_{0}(s,z,t) & \sim 0, \label{MC1} \\
\p_{z} Z_{1}(s,z,t) & \sim \nabla \zeta_{0}^{\pm} \cdot \n[\mminus,\mplus], \label{MC2} 
\end{align}
\begin{align}
\JpE[k,-1](s,z,t) & \sim 0, & \quad \JqE[-2](s,z,t) & \sim 0, \label{MC4} \\
\JpE[k,0](s,z,t) & \sim (\Jpe[k,0])^{\pm}, & \quad \JqE[-1](s,z,t) & \sim \Jqe[-1]^{\pm}, \label{MC5} \\
\JpE[k,1](s,z,t) & \sim (\Jpe[k,1])^{\pm} + z \nabla (\Jpe[k,0])^{\pm} \n[\mminus,\mplus], & \quad \JqE[0](s,z,t) & \sim \Jqe[0]^{\pm} + z \nabla \Jqe[-1]^{\pm} \n[\mminus,\mplus], \label{MC6}
\end{align}
where $(\cdot)^{\pm}$ denotes the limit $\lim_{\delta \searrow 0} (\cdot)(x \pm \delta \n[\mminus,\mplus])$ in  $x = \param[\mminus,\mplus](s,t) \in \G[\mminus,\mplus](t)$.

\subsection{Inner solutions}
\label{sec:in_sol}

The surfactant equation \eqref{eq:PFMsag1}, \eqref{eq:PFMsag2} to leading order $-3$ reads
\begin{align}
 0 &= - \n[\mminus,\mplus] \cdot \p_z \JqE[-2], \nonumber \\
 \JqE[-2] & = - \sum_{i<j} \Mcc[i,j] \big{(} \a[i,j](\phaE[0], \p_z \phaE[0] \otimes \n[\mminus,\mplus]) + \w[i,j](\phaE[0]) \big{)} \p_z \qqE[0] \n[\mminus,\mplus]. \label{eq:inH0}
\end{align}
Integrating with respect to $z$ from $-\infty$ to a variable denoted by $z$ again and using the matching condition \eqref{MC4} we conclude that $\p_z \qqE[0] = 0$ so that also all fields depending on $\qqE[0]$ such as $\sigt[i,j](\qqE[0])$ are constant across the interface layer to leading order. In particular, 
\begin{equation} \label{eq:inQ0}
\big{[} \qqe[0] \big{]}_{\mminus}^{\mplus} = 0, \quad \qqe[0]^\pm(\param[\mminus,\mplus](s,t),t) = \qqE[0](s,t)
\end{equation}
thanks to the matching condition \eqref{MC0}. Equation \eqref{eq:PFMmu} to order $-1$ then becomes
\begin{multline} \label{eq:inH1}
 0 = \sum_{i<j} \sigt[i,j](\qqE[0]) \Big{(} - \n[\mminus,\mplus] \cdot \frac{d}{dz} \big{(} \p_{\nabla \ph[l]} \a[i,j](\phaE[0], \p_z \phaE[0] \otimes \n[\mminus,\mplus]) \big{)} \\
 + \p_{\ph[l]} \a[i,j](\phaE[0], \p_z \phaE[0] \otimes \n[\mminus,\mplus]) + \p_{\ph[l]} \w[i,j](\phaE[0]) \Big{)}.
\end{multline}
This second order ODE in $z$ is supplied with the boundary conditions $\phaE[0] \sim \be[\mplus], \be[\mminus]$ and $\p_z \phaE[0] \sim 0$ as $z \to \pm \infty$, which are due to the matching conditions \eqref{MC0}, \eqref{MC1}. By Assumption \ref{ass:AA2} on the potentials $\a[i,j]$ and $\w[i,j]$ there are no third phase contributions, i.e., the leading order solution $\phaE[0]$ is such that $\phE[k,0] = 0$ if $k \not\in \{ \mplus, \mminus \}$. In fact, with choices as in \cite{GarNesSto99A,Sti05,BoyLap06}, for a wide range of surface energies $\gg[i,j]$ and related tensions $\sigt[i,j]$ the solution only depends on $z$ and is of the form
\begin{equation} \label{eq:inprofile}
 \phaE[0](z) = \chi(z) \be[\mplus] + (1-\chi(z)) \be[\mminus],
\end{equation}
with some monotone function $\chi : \R \to [0,1]$ (the \emph{transition profile}) satisfying 
\[
 \chi(0) = \frac{1}{2}, \quad \lim_{z \to \infty} \chi(z) = 1, \quad \lim_{z \to -\infty} \chi(z) = 0.
\]
The potentials are also such that for $i<j$ then $\a[i,j](\phaE[0], \p_z \phaE[0] \otimes \n[\mminus,\mplus]) = 0$ and $\w[i,j](\phaE[0]) = 0$ if $(i,j) \neq (\mminus,\mplus)$. Hence, $\phaE[0]$ satisfies 
\begin{multline} \label{eq:inH2}
 0 = \sigt[\mminus,\mplus](\qqE[0]) \Big{(} \p_{\ph[l]} \a[\mminus,\mplus](\phaE[0], \p_z \phaE[0] \otimes \n[\mminus,\mplus]) + \p_{\ph[l]} \w[\mminus,\mplus](\phaE[0]) \\
 - \frac{d}{dz} \big{(} \p_{\nabla \ph[l]} \a[\mminus,\mplus](\phaE[0], \p_z \phaE[0] \otimes \n[\mminus,\mplus]) \big{)} \cdot \n[\mminus,\mplus] \Big{)}.
\end{multline}
To avoid tracking of dimensionless calibration constants it is also convenient (and possible) to assume that the potentials are normalised in the sense that
\begin{equation} \label{eq:inH3}
 \int_{-\infty}^{\infty} \a[\mminus,\mplus](\phaE[0], \p_z \phaE[0] \otimes \n[\mminus,\mplus]) + \w[\mminus,\mplus](\phaE[0]) = 1.
\end{equation}
Multiplying \eqref{eq:inH2} with $\p_z \phE[l,0]$, summing over $l$, integrating with respect to $z$ from $-\infty$ to a variable denoted by $z$ again, using the two-homogeneity of the $\a[i,j]$ and using that $\qqE[0]$ is independent of $z$ we obtain (see \cite{GarNesSto98} for details on the calculation)
\begin{equation} \label{eq:equi}
 \mbox{\emph{equipartition of energy}:} \quad \sigt[\mminus,\mplus](\qqE[0]) \a[\mminus,\mplus](\phaE[0],\p_z \phaE[0] \otimes \n[\mminus,\mplus]) = \sigt[\mminus,\mplus](\qqE[0]) \w[\mminus,\mplus](\phaE[0]).
\end{equation}
In the following, for brevity, $\a[\mminus,\mplus]$ and its derivatives are evaluated at $(\phaE[0], \p_z \phaE[0] \otimes \n[\mminus,\mplus])$ and $\w[\mminus,\mplus]$ and its derivatives at $\phaE[0]$. 

Equation \eqref{eq:PFMmassbal} yields to leading order $-1$ that
\begin{equation} \label{eq:inH2b}
 \p_z \vvE[0] \cdot \n[\mminus,\mplus] = 0.
\end{equation}
Considering \eqref{eq:PFMphi1}, \eqref{eq:PFMphi2} to order $-2$ we obtain that 
\[
 0 = -\n[\mminus,\mplus] \cdot \p_z \JpE[k,-1], \quad \JpE[k,-1] = -\sum_l \L[k,l]_0(\phaE[0],\qqE[0]) \p_z \mbE[l,0] \n[\mminus,\mplus].
\]
After integrating the first identity with respect to $z$ from $-\infty$ to a variable denoted by $z$ again and using the matching condition \eqref{MC4} we conclude that 
\begin{equation} \label{eq:inH3b}
 \JpE[k,-1] = 0.
\end{equation}
As $\phaE[0]$ is no corner of the Gibb's simplex, the kernel of $\{ \L[k,l]_0(\phaE[0],\qqE[0]) \}_{k,l}$ is the span of $\bbb[1]$ by assumption, hence there is a scalar function $\psi(s,z,t)$ such that
\begin{equation} \label{eq:inH3c}
 \p_z \mbaE[0] = \psi(s,z,t) \bbb[1].
\end{equation}
Using \eqref{eq:inH3b}, the momentum equation \eqref{eq:PFMmombal} to order $-2$ becomes
\begin{align}
 0 =& \, \n[\mminus,\mplus] \cdot \frac{d}{dz} \big{(} \eta(\phaE[0]) (\p_z \vvE[0] \otimes \n[\mminus,\mplus] + \n[\mminus,\mplus] \otimes \p_z \vvE[0]) \big{)} \nonumber \\
 & + \sigt[\mminus,\mplus](\qqE[0]) \frac{d}{dz} \big{(} \a[\mminus,\mplus] + \w[\mminus,\mplus] \big{)} \n[\mminus,\mplus] \nonumber \\
 & - \sigt[\mminus,\mplus](\qqE[0]) \n[\mminus,\mplus] \cdot \frac{d}{dz} \Big{(} \sum_k \p_z \phE[k,0] \n[\mminus,\mplus] \otimes \p_{\nabla \ph[k]} \a[\mminus,\mplus] \Big{)}. \label{eq:inH4}
\end{align}
Note that by the two-homogeneity of the $\a[i,j]$ in the second argument
\begin{align*}
 & \, \n[\mminus,\mplus] \cdot \Big{(} \sum_k \p_z \phE[k,0] \n[\mminus,\mplus] \otimes \p_{\nabla \ph[k]} \a[\mminus,\mplus] \Big{)} \\
=& \, \sum_k \big{(} \p_{\nabla \ph[k]} \a[\mminus,\mplus] \cdot (\p_z \phE[k,0] \n[\mminus,\mplus]) \big{)} \n[\mminus,\mplus] \\
=& \, \big{(} \p_{\nabla \pha} \a[\mminus,\mplus] : (\p_z \phaE[0] \otimes \n[\mminus,\mplus]) \big{)} \n[\mminus,\mplus] = \big{(} 2 \a[\mminus,\mplus] \big{)} \n[\mminus,\mplus].
\end{align*}
Thus, thanks to the equipartition of energy \eqref{eq:equi} the terms involving $\sigt[\mminus,\mplus](\qqE[0])$ in \eqref{eq:inH4} cancel out and we obtain that 
\begin{equation} \label{eq:inH4b}
 0 = \eta(\phaE[0]) \p_z \vvE[0].
\end{equation}
Hence, $\vvE[0]$ is independent of $z$, and with the matching conditions we conclude that the velocity is continuous to leading order across the interface, 
\begin{equation} \label{eq:inV0}
 \big{[} \vve[0] \big{]}_{\mminus}^{\mplus} = 0, \quad \vve[0]^\pm(\param[\mminus,\mplus](s,t),t) = \vvE[0](s,t).
\end{equation}

To order $0$, \eqref{eq:PFMmassbal} gives
\begin{equation} \label{eq:inH5}
 \p_z \vvE[1] \cdot \n[\mminus,\mplus] + \p_s \vvE[0]\cdot \tang[\mminus,\mplus] = 0.
\end{equation}
We continue with \eqref{eq:PFMphi1} which, to order $-1$, becomes
\[ 
 \big{(} -\u[\mminus,\mplus] + \vvE[0] \cdot \n[\mminus,\mplus] \big{)} \p_z \phE[k,0] = \n[\mminus,\mplus] \cdot \p_z \JpE[k,0].
\] 
Integrating from $-\infty$ to $\infty$, using \eqref{eq:inH2b}, the matching condition \eqref{MC5}, and \eqref{eq:outH2} the right-hand side vanishes and we see that the interface normal velocity is given by the fluid velocity in normal direction, 
\begin{equation} \label{eq:inR1}
 \u[\mminus,\mplus] = \vvE[0] \cdot \n[\mminus,\mplus].
\end{equation}
Proceeding now as for $\JpE[k,-1]$ and using \eqref{MC5} and \eqref{eq:outH2} we see that 
\begin{equation} \label{eq:inH5b}
 \JpE[k,0] = 0.
\end{equation}
Recalling that $\JqE[-2] = 0$ from \eqref{eq:inH0} and that $\p_z \qqE[0] = 0$, equations \eqref{eq:PFMsag1}, \eqref{eq:PFMsag2} to order $-2$ read
\begin{align}
 0 &= - \n[\mminus,\mplus] \cdot \p_z \JqE[-1], \nonumber \\
 \JqE[-1] & = - \Mcc[\mminus,\mplus] \big{(} \a[\mminus,\mplus] + \w[\mminus,\mplus] \big{)} \big{(} \p_z \qqE[1] \n[\mminus,\mplus] + \p_s \qqE[0] \tang[\mminus,\mplus] \big{)}. \label{eq:inH9}
\end{align}
We can proceed as for $\qqE[0]$ to conclude that $\p_z \qqE[1] = 0$. For this purpose, we integrate with respect to $z$ and use $\n[\mminus,\mplus] \cdot \tang[\mminus,\mplus] =0$, the matching condition \eqref{MC5}, and the fact that $\Jqe[-1] = 0$, see \eqref{eq:outH4}. Regarding \eqref{eq:PFMmu} to order $0$ we obtain that (here, we dropped terms with $\p_s \phaE[0]$ as $\phaE[0]$ depends on $z$ only)
\begin{align}
 \mbE[l,0] = & \, \lamt[l](\qqE[0]) \xi_l'(\phE[l,0]) \nonumber \\
 & + \sigt[\mminus,\mplus]'(\qqE[0]) \qqE[1] \Big{(} \p_{\ph[l]} \a[\mminus,\mplus] + \p_{\ph[l]} \w[\mminus,\mplus] - \n[\mminus,\mplus] \cdot \frac{d}{dz} \big{(} \p_{\nabla \ph[l]} \a[\mminus,\mplus] \big{)} \Big{)} \nonumber \\
 & + \sigt[\mminus,\mplus](\qqE[0]) \sum_k \Big{(} \p_{\ph[k]} (\p_{\ph[l]} \a[\mminus,\mplus]) \phE[k,1] + \p_{\nabla \ph[k]} (\p_{\ph[l]} \a[\mminus,\mplus]) \n[\mminus,\mplus] \p_z \phE[k,1] \Big{)} \nonumber \\
 & + \sigt[\mminus,\mplus](\qqE[0]) \sum_k \Big{(} \p_{\ph[k]} (\p_{\ph[l]} \w[\mminus,\mplus]) \phE[k,1] \Big{)} \nonumber \\
 & - \sigt[\mminus,\mplus](\qqE[0]) \n[\mminus,\mplus] \cdot \sum_k \frac{d}{dz} \Big{(} \p_{\ph[k]} (\p_{\nabla \ph[l]} \a[\mminus,\mplus]) \phE[k,1] + \p_{\nabla \ph[k]} (\p_{\nabla \ph[l]} \a[\mminus,\mplus]) \n[\mminus,\mplus] \p_z \phE[k,1] \Big{)} \nonumber \\
 & - \tang[\mminus,\mplus] \cdot \frac{d}{ds} \big{(} \sigt[\mminus,\mplus](\qqE[0]) \p_{\nabla \ph[l]} \a[\mminus,\mplus] \big{)}. \label{eq:inH6}
\end{align}
Note that the terms involving $\qqE[1]$ are zero thanks to \eqref{eq:inH2}. We multiply with $\p_z \phE[l,0]$, sum up over $l$, and integrate with respect to $z$ from $-\infty$ to $\infty$. Recalling \eqref{eq:inH3c}, \eqref{eq:inprofile}, and noting that $\p_z \phaE[0] \in T\Sigma^M$ we obtain for the left-hand side thanks to the matching condition \eqref{MC0} for $\mu$ that
\[
 \int_{-\infty}^{\infty} \sum_l \mbE[l,0] \p_z \phE[l,0] = \mbE[\mplus,0] - \mbE[\mminus,0] = \big{[} \mbe[\cdot,0] \big{]}_{\mminus}^{\mplus}.
\]
Similarly, the first term on the right-hand side of \eqref{eq:inH6} becomes
\[
 \int_{-\infty}^{\infty} \sum_l \lamt[l](\qqE[0]) \xi_l'(\phE[l,0]) \p_z \phE[l,0] = \sum_l \lamt[l](\qqe[0]) \int_{-\infty}^{\infty} \p_z \xi_l(\phE[l,0]) = \lamt[\mplus](\qqe[0]) - \lamt[\mminus](\qqe[0]) = \big{[} \lamt[(\cdot)](\qqe[0]) \big{]}_{\mminus}^{\mplus}.
\]
Thanks to the isotropy of the surface energies $\p_{\nabla \ph[l]} \a[i,j](\phaE[0], \p_z \phaE[0] \otimes \n[\mminus,\mplus])$ points in the normal direction, 
\begin{equation} \label{eq:propanormal}
 \p_{\nabla \ph[k]} \a[i,j](\phaE[0], \p_z \phaE[0] \otimes \n[\mminus,\mplus]) = \big{(} \p_{\nabla \ph[k]} \a[i,j](\phaE[0], \p_z \phaE[0] \otimes \n[\mminus,\mplus]) \cdot \n[\mminus,\mplus] \big{)} \, \n[\mminus,\mplus],
\end{equation}
so that $\p_{\nabla \ph[l]} \a[i,j] \cdot \tang[\mminus,\mplus] = 0$. Using that $\p_s \tang[\mminus,\mplus] = \scurv[\mminus,\mplus] \n[\mminus,\mplus]$, the two-homogeneity of $\a[\mminus,\mplus]$ in the second argument, and \eqref{eq:inH3} with \eqref{eq:equi} we obtain for the last term of \eqref{eq:inH6} that
\begin{align*}
 & - \int_{-\infty}^{\infty} \sum_l \frac{d}{ds} \big{(} \sigt[\mminus,\mplus](\qqE[0]) \p_{\nabla \ph[l]} \a[\mminus,\mplus] \big{)} \cdot \tang[\mminus,\mplus] \p_z \phE[l,0] \\
 = & \int_{-\infty}^{\infty} \sum_l \sigt[\mminus,\mplus](\qqE[0]) \p_{\nabla \ph[l]} \a[\mminus,\mplus] \cdot \scurv[\mminus,\mplus] \n[\mminus,\mplus] \p_z \phE[l,0] \\
 = & \int_{-\infty}^{\infty} \sigt[\mminus,\mplus](\qqE[0]) \p_{\npha} \a[\mminus,\mplus](\phaE[0], \p_z \phaE[0] \otimes \n[\mminus,\mplus]) : \p_z \phaE[0] \otimes \n[\mminus,\mplus] \, \scurv[\mminus,\mplus] \\
 = & \, \sigt[\mminus,\mplus](\qqe[0]) \scurv[\mminus,\mplus] \int_{-\infty}^{\infty} 2 \a[\mminus,\mplus](\phaE[0], \p_z \phaE[0] \otimes \n[\mminus,\mplus]) \\
 = & \, \sigt[\mminus,\mplus](\qqe[0]) \scurv[\mminus,\mplus].
\end{align*}
Using integration by parts, the matching conditions \eqref{MC0}--\eqref{MC2}, and \eqref{eq:inH2} again one can show that the terms involving $\phaE[1]$ in \eqref{eq:inH6} sum up to zero, see \cite{GarNesSto98} for the details. In fact, one can consider \eqref{eq:inH6} as a differential equation for $\phaE[1]$ where the operator has a non-trivial kernel containing the span of $\p_z \phaE[0]$. The remaining terms of the mentioned operations then yield a solvability condition for $\phaE[1]$, which reads
\begin{equation} \label{eq:inH7}
 \big{[} \mbe[\cdot,0] \big{]}_{\mminus}^{\mplus} = \big{[} \lamt[(\cdot)](\qqe[0]) \big{]}_{\mminus}^{\mplus} + \sigt[\mminus,\mplus](\qqe[0]) \scurv[\mminus,\mplus].
\end{equation}

Now consider the surfactant equation \eqref{eq:PFMsag1} to order $-1$. Using \eqref{eq:inR1} and that $\p_z \qqE[0] = 0$ it reduces to
\begin{multline} \label{eq:inH8}
 \Big{(}\frac{d}{dt} - \p_t \param[\mminus,\mplus] \frac{d}{ds} + \vvE[0] \cdot \tang[\mminus,\mplus] \frac{d}{ds} + \vvE[1] \cdot \n[\mminus,\mplus] \frac{d}{dz} \Big{)} \big{(} (\a[\mminus,\mplus] + \w[\mminus,\mplus]) \cc[\mminus,\mplus](\qqE[0]) \big{)} \\
 = - \n[\mminus,\mplus] \cdot \p_z \JqE[0] - \tang[\mminus,\mplus] \cdot \p_s \JqE[-1].
\end{multline}
We integrate with respect to $z$ from $-\infty$ to $\infty$. Thanks to the equipartition of energy \eqref{eq:equi}, that $\phaE[0]$ is independent of $(s,t)$, the normalisation \eqref{eq:inH3}, the identity \eqref{eq:identnd}, and the matching condition \eqref{MC0} for $\qq$, the first and the second term on the left-hand side yield
\begin{align*}
  & \, \int_{-\infty}^{\infty} \Big{(} \frac{d}{dt} - \p_t \param[\mminus,\mplus] \frac{d}{ds} \Big{)} \big{(} (\a[\mminus,\mplus] + \w[\mminus,\mplus]) \cc[\mminus,\mplus](\qqE[0]) \big{)} \\
= & \, \int_{-\infty}^{\infty} \Big{(} \frac{d}{dt} - \p_t \param[\mminus,\mplus] \frac{d}{ds} \Big{)} \big{(} (2 \w[\mminus,\mplus]) \cc[\mminus,\mplus](\qqE[0]) \big{)} \\
= & \, \int_{-\infty}^{\infty} (2 \w[\mminus,\mplus]) dz \, \Big{(} \frac{d}{dt} - \p_t \param[\mminus,\mplus] \frac{d}{ds} \Big{)} \cc[\mminus,\mplus](\qqe[0]) \\
= & \, \ndu[\mminus,\mplus] \cc[\mminus,\mplus](\qqe[0]).
\end{align*}
For the third term on the left-hand side of \eqref{eq:inH8} we see that
\begin{multline*}
 \int_{-\infty}^{\infty} \vvE[0] \cdot \tang[\mminus,\mplus] \p_s \big{(} (\a[\mminus,\mplus] + \w[\mminus,\mplus]) \cc[\mminus,\mplus](\qqE[0]) \big{)} \\
 = \int_{-\infty}^{\infty} (2 \w[\mminus,\mplus]) dz \, \vvE[0] \cdot \tang[\mminus,\mplus] \p_s \cc[\mminus,\mplus](\qqe[0]) = \vvE[0] \cdot \tang[\mminus,\mplus] \p_s \cc[\mminus,\mplus](\qqe[0]).
\end{multline*}
For the forth term on the left-hand side of \eqref{eq:inH8} we also use \eqref{eq:inH5} and the matching conditions \eqref{MC0}, \eqref{MC1} for $\pha$ so that the boundary terms in the following integration by parts vanish:
\begin{align*}
  & \, \int_{-\infty}^{\infty} \vvE[1] \cdot \n[\mminus,\mplus] \p_z \big{(} (\a[\mminus,\mplus] + \w[\mminus,\mplus]) \cc[\mminus,\mplus](\qqE[0]) \big{)} \\
= & \, \int_{-\infty}^{\infty} - \p_z \vvE[1] \cdot \n[\mminus,\mplus] \big{(} (\a[\mminus,\mplus] + \w[\mminus,\mplus]) \cc[\mminus,\mplus](\qqE[0]) \big{)} \\
= & \, \int_{-\infty}^{\infty} \p_s \vvE[0] \cdot \tang[\mminus,\mplus] (\a[\mminus,\mplus] + \w[\mminus,\mplus]) \, \cc[\mminus,\mplus](\qqE[0]) \\
= & \, (\tang[\mminus,\mplus] \cdot \p_s \vve[0]) \cc[\mminus,\mplus](\qqe[0]). 
\end{align*}
With the matching condition \eqref{MC6} and the identities \eqref{eq:outH4}, \eqref{eq:outH5} the first term on the right-hand side of \eqref{eq:inH8} gives
\[
 \int_{-\infty}^{\infty} -\n[\mminus,\mplus] \cdot \p_z \JqE[0] = \n[\mminus,\mplus] \cdot \big{[} \Mc[\cdot] \nabla \qqe[0] \big{]}_{\mminus}^{\mplus},
\]
whilst for the second term on the right-hand side of \eqref{eq:inH8} we obtain with \eqref{eq:inH9} and \eqref{eq:Frenet} that
\begin{align*}
  & \, \int_{-\infty}^{\infty} - \tang[\mminus,\mplus] \cdot \p_s \JqE[-1] \\
= & \, \int_{-\infty}^{\infty} \tang[\mminus,\mplus] \cdot \frac{d}{ds} \big{(} \Mcc[\mminus,\mplus] 2 \w[\mminus,\mplus] \tang[\mminus,\mplus] \p_s \qqE[0] \big{)} \\
= & \, \tang[\mminus,\mplus] \cdot \frac{d}{ds} \big{(} \Mcc[\mminus,\mplus] \p_s \qqe[0] \tang[\mminus,\mplus] \big{)} \\
\end{align*}
Altogether, \eqref{eq:inH8} yields 
\begin{multline} \label{eq:inR2}
 \ndu[\mminus,\mplus] \cc[\mminus,\mplus](\qqe[0]) + \tang[\mminus,\mplus] \cdot \p_s \big{(} \vve[0] \, \cc[\mminus,\mplus](\qqe[0]) \big{)} \\
 = \n[\mminus,\mplus] \cdot \big{[} \Mc[\cdot] \nabla \qqe[0] \big{]}_{\mminus}^{\mplus} + \tang[\mminus,\mplus] \cdot \frac{d}{ds} \big{(} \Mcc[\mminus,\mplus] \p_s \qqe[0] \tang[\mminus,\mplus] \big{)},
\end{multline}
which is the surface surfactant equation in the form \eqref{eq:sagsurf2B}.

Finally, the left-hand side in equation \eqref{eq:PFMmombal} to order $-1$ is zero thanks to \eqref{eq:inR1}. Using also that $\JpE[k,-1] = \JpE[k,0] = 0$ for all $k \in \{ 1, \dots, M \}$ (see \eqref{eq:inH3b}, \eqref{eq:inH5b}) it reads
\begin{align}
 0 = & \, - \p_z \ppE[0] \n[\mminus,\mplus] + 2 \n[\mminus,\mplus] \cdot \frac{d}{dz} \big{(} \eta(\phaE[0]) (D(\p_z \vvE[1] \otimes \n[\mminus,\mplus]) + D(\p_s \vvE[0] \otimes \tang[\mminus,\mplus]) ) \big{)} \nonumber \\
 & \, + \tang[\mminus,\mplus] \cdot \frac{d}{ds} \Big{(} \sigt[\mminus,\mplus](\qqE[0]) (\beta_{(\mminus,\mplus)})_0 \Big{)} \nonumber \\
 & \, + \n[\mminus,\mplus] \cdot \frac{d}{dz} \Big{(} \sigt[\mminus,\mplus]'(\qqE[0]) \qqE[1] (\beta_{(\mminus,\mplus)})_0 + \sigt[\mminus,\mplus](\qqE[0]) (\beta_{(\mminus,\mplus)})_1 \Big{)}, \label{eq:inH10}
\end{align}
where 
\begin{align*}
 (\beta_{(\mminus,\mplus)})_0 \, &= (\a[\mminus,\mplus] + \w[\mminus,\mplus]) \bbb[I] - \sum_k \p_z \phE[k,0] \n[\mminus,\mplus] \otimes \p_{\nabla \ph[k]} \a[\mminus,\mplus], \\
 (\beta_{(\mminus,\mplus)})_1 \, &= (\p_{\pha} \a[\mminus,\mplus] \cdot \phaE[1] + \p_{\nabla \pha} \a[\mminus,\mplus] : (\p_z \phaE[1] \otimes \n[\mminus,\mplus] + \p_s \phaE[0] \otimes \tang[\mminus,\mplus]) + \p_{\pha} \w[\mminus,\mplus] \cdot \phaE[1]) \bbb[I] \\
 & \quad + \sum_k \big{(} \p_z \phE[k,1] \n[\mminus,\mplus] \otimes \p_{\nabla \ph[k]} \a[\mminus,\mplus] + \p_z \phE[k,0] \n[\mminus,\mplus] \otimes \p_{\pha} (\p_{\nabla \ph[k]} \a[\mminus,\mplus]) \cdot \phaE[1] \big{)} \\
 & \quad + \sum_k \big{(} \p_z \phE[k,0] \n[\mminus,\mplus] \otimes \p_{\nabla \pha} (\p_{\nabla \ph[k]} \a[\mminus,\mplus]) : (\p_z \phaE[1] \otimes \n[\mminus,\mplus] + \p_s \phaE[0] \otimes \tang[\mminus,\mplus]) \big{)}.
\end{align*}
Note that thanks to \eqref{eq:equi}
\begin{multline*}
 \n[\mminus,\mplus] \cdot (\beta_{(\mminus,\mplus)})_0 = (\a[\mminus,\mplus] + \w[\mminus,\mplus]) \n[\mminus,\mplus] - (\p_{\nabla \ph[k]} \a[\mminus,\mplus] : \p_z \phaE[0] \otimes \n[\mminus,\mplus]) \, \n[\mminus,\mplus] \\
 = (\a[\mminus,\mplus] + \w[\mminus,\mplus] - 2 \a[\mminus,\mplus]) \n[\mminus,\mplus] = 0,
\end{multline*}
so that the $\qqE[1]$ term vanishes in \eqref{eq:inH10}. We integrate that equation with respect to $z$ from $-\infty$ to $\infty$. Thanks to the matching conditions \eqref{MC0}, \eqref{MC2} for $\vv$ we see that, as $z \to \pm \infty$, 
\[
 \p_z \vvE[1] \otimes \n[\mminus,\mplus] + \p_s \vvE[0] \otimes \tang[\mminus,\mplus] \sim \nabla \vve[0].
\]
Therefore, after integrating with respect to $z$ from $-\infty$ to $\infty$, the first line of \eqref{eq:inH10} yields the jump of the normal stresses $[ - \ppe[0] \bbb[I] + 2 \et[\cdot] D(\vve[0]) ]_{\mminus}^{\mplus} \n[\mminus,\mplus]$.
With the matching conditions \eqref{MC0}--\eqref{MC2} applied to $\pha$ we see that $(\beta_{(\mminus,\mplus)})_1 \to 0$ as $z \to \pm \infty$ so that also the last term in \eqref{eq:inH10} vanishes after integrating. With respect to the second line we first recall \eqref{eq:propanormal} ($\p_{\nabla \ph[k]} \a[i,j](\phaE[0], \p_z \phaE[0] \otimes \n[\mminus,\mplus])$ points in the normal direction). Hence, with \eqref{eq:Frenet} we conclude that
\begin{align*}
 & \, \int_{-\infty}^{\infty} \tang[\mminus,\mplus] \cdot \frac{d}{ds} \Big{(} \sigt[\mminus,\mplus](\qqE[0]) (\beta_{(\mminus,\mplus)})_0 \Big{)} \nonumber \\
 = & \, \int_{-\infty}^{\infty} \big{(} \p_s \sigt[\mminus,\mplus](\qqE[0]) \big{)} \tang[\mminus,\mplus] (\a[\mminus,\mplus] + \w[\mminus,\mplus]) \nonumber \\
  & - \int_{-\infty}^{\infty} \sigt[\mminus,\mplus](\qqE[0]) \sum_k \p_z \phE[k,0] \n[\mminus,\mplus] \big{(} \p_{\nabla \ph[k]} \a[\mminus,\mplus] \cdot \n[\mminus,\mplus] \big{)} \, \p_s \n[\mminus,\mplus] \cdot \tang[\mminus,\mplus] \nonumber \displaybreak[0] \\
 = & \, \p_s \sigt[\mminus,\mplus](\qqe[0]) \tang[\mminus,\mplus] \int_{-\infty}^{\infty} (\a[\mminus,\mplus] + \w[\mminus,\mplus]) \nonumber \\
  & + \sigt[\mminus,\mplus](\qqe[0]) \int_{-\infty}^{\infty} \n[\mminus,\mplus] \big{(} \p_{\nabla \pha} \a[\mminus,\mplus] : \p_z \phaE[0] \otimes \n[\mminus,\mplus] \big{)} \, \scurv[\mminus,\mplus] \nonumber \displaybreak[0] \\
 = & \, \tang[\mminus,\mplus] \p_s \sigt[\mminus,\mplus](\qqe[0]) + \sigt[\mminus,\mplus](\qqe[0]) \curv[\mminus,\mplus],
\end{align*}
where we used the two-homogeneity of the $\a[i,j]$ and the equipartition of energy \eqref{eq:equi} for the last identity. Altogether, we obtain from \eqref{eq:inH10} that
\begin{equation} \label{eq:inR3}
 0 = \big{[} - \ppe[0] \bbb[I] + 2 \et[\cdot] D(\vve[0]) \big{]}_{\mminus}^{\mplus} \n[\mminus,\mplus] + \tang[\mminus,\mplus] \p_s \sigt[\mminus,\mplus](\qqe[0]) + \sigt[\mminus,\mplus](\qqe[0]) \curv[\mminus,\mplus].
\end{equation}
Together with \eqref{eq:inR1}, \eqref{eq:inR2} the recovery of the interface equations \eqref{eq:intvelcont2}--\eqref{eq:sagsurf2} of the sharp interface model is thus completed.

\subsection{Triple point expansions and matching conditions}
\label{sec:tri_expmat}

Still in the case $d=2$ we now consider a triple junction between three phases where $\phae[0] \approx \be[\mIn], \be[\mIp], \be[\mIm]$, respectively, with pairwise different indices $\mIn, \mIp, \mIm \in \{ 1, \dots, M \}$. These are separated by layers that converge to curves belonging to $\G[\mIn,\mIp]$, $\G[\mIp,\mIm]$, or $\G[\mIm,\mIn]$ as $\e \to 0$ and then form a triple point denoted by $\triple[\mIn,\mIp,\mIm](t) \in \T[\mIn,\mIp,\mIm]$. The thickness of the layers scaling with $\e$ motivates to introduce a local rescaled coordinate,
\begin{equation} \label{eq:tri_def_coord}
 y := \frac{x - \triple[\mIn,\mIp,\mIm](t)}{\e}, \quad x \mbox{ close to } \triple[\mIn,\mIp,\mIm](t).
\end{equation}
Differential operators then transform as follows: For a function $\zeta(x,t) = \hat{Z}(y,t)$
\begin{align*}
 \p_t \zeta(x,t) &= - \e^{-1} \nabla_y \hat{Z}(y,t) \cdot \uu[\mIn,\mIp,\mIm] + \p_t \hat{Z}(y,t), \\
 \nabla \zeta(x,t) &= \e^{-1} \nabla_y \hat{Z}(y,t),
\end{align*}
where the velocity $\uu[\mIn,\mIp,\mIm]$ is evaluated in $\triple[\mIn,\mIp,\mIm]$. Expansions of $\e$-solution fields to be plugged into the equations now read:
\begin{align}
 \zeta_\e(x,t) & = \hat{Z}_0(y,t) + \e \hat{Z}_1(y,t) + \e^2 \hat{Z}_2(y,t) + \dots, \label{eq:triexp} \\
 \Jp[k](x,t) & = \e^{-1} \JpT[k,-1](y,t) + \e^{0} \JpT[k,0](y,t) + \e^{1} \JpT[k,1](y,t) + \dots, \nonumber \\
 \Jq(x,t) & = \e^{-2} \JqT[-2](y,t) + \e^{-1} \JqT[-1](y,t) + \e^{0} \JqT[0](y,t) + \dots. \nonumber 
\end{align}

\begin{figure}[ht]
 \centering
 \includegraphics[height = 5cm]{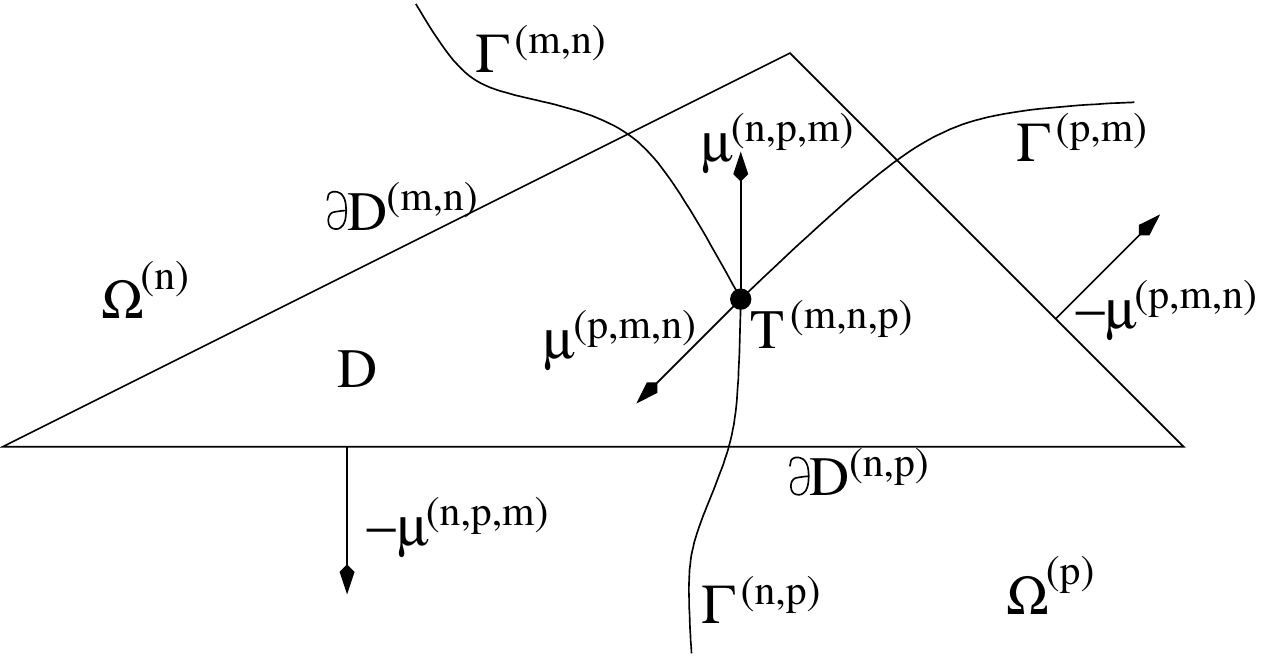}
 \caption{\em Illustration of the triangle for the asymptotic analysis around a triple point.} \label{fig:triangleD}
\end{figure}

As in \cite{BroRei93,GarNesSto98}, around $\triple[\mIn,\mIp,\mIm]$ we consider a triangle $D$ with three edges $\p D^{(\mIn,\mIp)}$, $\p D^{(\mIp,\mIm)}$, $\p D^{(\mIm,\mIn)}$ that have length scaling with $\e^{1/2}$. In addition, we assume that the exterior unit normal on $\p D^{(\mIn,\mIp)}$ coincides with $-\m[\mIn,\mIp,\mIm]$ and similarly for the other two edges. See Figure \ref{fig:triangleD} for an illustration. We remark that this additional assumption is for convenience only as the matching conditions below in \eqref{MCH1}--\eqref{MCH5} are easier to see, but it is not required: As in \cite{BroRei93,GarNesSto98} a more general triangle could be considered.

The parametrisation $\param[\mIn,\mIp]$ of $\G[\mIn,\mIp]$ is assumed to start in $s=0$ for simplicity, i.e., $\param[\mIn,\mIp](0,t) = \triple[\mIn,\mIp,\mIm](t)$. Requiring the local triple point expansion to match with the expansions along the interface layers (defined in Section \ref{sec:in_expmat}) also leads to conditions: As $\e \to 0$ for $y \in \p D^{(\mIn,\mIp)}$  
\begin{align}
 \hat{Z}_0(y,t) & \sim Z_0(0,y \cdot \n[\mIn,\mIp],t), \label{MCH1} \\
 \nabla_y \hat{Z}_0(y,t) & \sim \p_z Z_0(0,y \cdot \n[\mIn,\mIp],t) \n[\mIn,\mIp]. \label{MCH2}
\end{align}
For the fluxes of $\qq$ we have thanks to \eqref{eq:inH3}, \eqref{eq:inH9} and $\p_z \qqE[0] = \p_z \qqE[1] = 0$ that
\begin{align}
 \JqT[-2](y,t) & \sim \JqE[-2](0,y \cdot \n[\mIn,\mIp],t) = 0, \label{MCH3} \\
 \JqT[-1](y,t) & \sim \JqE[-1](0,y \cdot \n[\mIn,\mIp],t) = \Mcc[\mIn,\mIp] \big{(} \a[\mIn,\mIp] + \w[\mIn,\mIp] \big{)} \p_s \qqE[0](0,t) \m[\mIn,\mIp,\mIm], \label{MCH4}
\end{align}
where $\a[\mIn,\mIp]$ is evaluated in $(\phaE[0](z), \p_z \phaE[0](z) \otimes \n[\mIn,\mIp])$ and $\w[\mIn,\mIp]$ in $\phaE[0](z)$ with $z = y \cdot \n[\mIn,\mIp]$. Similarly, for the fluxes of the $\ph[k]$ we have
\begin{equation} \label{MCH5}
 \JpT[k,-1](y,t) \sim \JpE[k,-1](0,y \cdot \n[\mIn,\mIp],t) = 0, \quad \JpT[k,0](y,t) \sim \JpE[k,0](0,y \cdot \n[\mIn,\mIp],t) = 0,
\end{equation}
where we used \eqref{eq:inH3b} and \eqref{eq:inH5b}.

\subsection{Triple point solutions}
\label{sec:tri_sol}

For the surfactant equation \eqref{eq:PFMsag1}, \eqref{eq:PFMsag2} we obtain to leading order
\begin{equation} \label{eq:triH0}
\begin{split}
 0 &= - \nabla_y \cdot \JqT[-2], \\
 \JqT[-2] &= - \sum_{i<j} \Mcc[i,j] \big{(} \a[i,j](\phaT[0], \nabla_y \phaT[0]) + \w[i,j](\phaT[0]) \big{)} \nabla_y \qqT[0].
\end{split}
\end{equation}
The matching condition \eqref{MCH3} motivates to consider it as a PDE for $\qqT[0]$ with a no-flux boundary condition which, using the matching conditions \eqref{MCH1}, \eqref{MCH2} for $\pha$, reads (on $\p D^{(\mIn,\mIp)}$, similarly on the other two edges)
\begin{equation} \label{eq:triH0b} 
 0 = \Mcc[\mIn,\mIp] \big{(} \a[\mIn,\mIp](\phaE[0], \p_z \phaE[0] \otimes \n[\mIn,\mIp]) + \w[\mIn,\mIp](\phaE[0]) \big{)} \nabla_y \qqT[0] \cdot \m[\mIn,\mIp,\mIm].
\end{equation}
By multiplying \eqref{eq:triH0} with $\qqT[0]$, integrating by parts, and using this boundary condition we see that $\nabla_y \qqT[0] = 0$ provided that the prefactor doesn't degenerate, i.e., 
\begin{equation} \label{eq:Q0const}
 \sum_{i<j} \Mcc[i,j] \big{(} \a[i,j](\phaT[0], \nabla_y \phaT[0]) + \w[i,j](\phaT[0]) \big{)} \neq 0. 
\end{equation}
That this indeed is the case is the meaning of Assumption \ref{ass:AA3}. With \eqref{MCH1}, $\qqT[0]$ being constant implies that the limits of the values of the $\qqE[0]$ along the interfaces match up in $\triple[\mIn,\mIp,\mIm](t)$. As $\qqe[0]$ is continuous across the interfaces (see \eqref{eq:inQ0}) we obtain that $\qqe[0]$ is continuous at the triple junctions and, thus, in the whole domain. 

Using that also the $\sigt[i,j](\qqT[0])$ are constant, \eqref{eq:PFMmu} reads to order $-1$
\begin{equation} \label{eq:triH1}
 0 = \sum_{i<j} \sigt[i,j](\qqT[0]) \Big{(} \p_{\ph[l]} \big{(} \a[i,j](\phaT[0], \nabla_y \phaT[0]) + \w[i,j](\phaT[0]) \big{)} - \nabla_y \cdot \sum_l \p_{\nabla \ph[l]} \a[i,j](\phaT[0], \nabla_y \phaT[0]) \Big{)}.
\end{equation}
One can now proceed exactly as described in \cite{GarNesSto98}, Section 5.5, (see also \cite{BroRei93} for the techniques) and deduce the following identity, which is a solvability condition:
\begin{equation}  \label{eq:triR1}
 \sigt[\mIn,\mIp](\qqT[0]) \m[\mIn,\mIp,\mIm] + \sigt[\mIp,\mIm](\qqT[0]) \m[\mIp,\mIm,\mIn] + \sigt[\mIm,\mIn](\qqT[0]) \m[\mIm,\mIn,\mIp] = 0.
\end{equation}

Using that $\qqT[0]$ is constant, the surfactant equation \eqref{eq:PFMsag1}, \eqref{eq:PFMsag2} to the next order is
\begin{align}
 0 &= - \nabla_y \cdot \JqT[-1], \nonumber \\
 \JqT[-1] &= - \sum_{i<j} \Mcc[i,j] \big{(} \a[i,j](\phaT[0], \nabla_y \phaT[0]) + \w[i,j](\phaT[0]) \big{)} \nabla_y \qqT[1]. \label{eq:triH1b}
\end{align}
We integrate over $D$, use the divergence theorem, the matching condition \eqref{MCH4} and the normalisation \eqref{eq:inH3}:
\begin{align}
 0 \, &= -\int_D \nabla_y \cdot \JqT[-1] \nonumber \\
 &= \int_{ \p D^{(\mIn,\mIp)} } \JqT[-1] \cdot \m[\mIn,\mIp,\mIm] + \int_{ \p D^{(\mIp,\mIm)} } \JqT[-1] \cdot \m[\mIp,\mIm,\mIn] + \int_{ \p D^{(\mIm,\mIn)} } \JqT[-1] \cdot \m[\mIm,\mIn,\mIp] \nonumber \\
 &\to \sum_{ (k,l) \in \{ (\mIn,\mIp), (\mIp,\mIm), (\mIm,\mIn) \} } \Mcc[k,l] \int_{\R} \big{(} \a[k,l](\phaE[0], \p_z \phaE[0] \otimes \n[k,l]) + \w[k,l](\phaE[0]) \big{)} \p_s \qqE[0] \nonumber \\
 &= \Mcc[\mIn,\mIp] \nabg[\mIn,\mIp] \qqe[0] \cdot \m[\mIn,\mIp,\mIm] + \Mcc[\mIp,\mIm] \nabg[\mIp,\mIm] \qqe[0] \cdot \m[\mIp,\mIm,\mIn] + \Mcc[\mIm,\mIn] \nabg[\mIm,\mIn] \qqe[0] \cdot \m[\mIm,\mIn,\mIp], \label{eq:triR2}
\end{align}
where the interface fields $\phaE[0]$ and $\qqE[0]$ in the third line are evaluated in the boundary point, i.e., the triple junction. 

To leading order, equations \eqref{eq:PFMphi1}, \eqref{eq:PFMphi2} are
\[
 0 = - \nabla_y \cdot \JpT[k,-1], \quad \JpT[k,-1] = - \sum_l \L[k,l]_0 (\phaT[0],\qqT[0]) \nabla_y \mbT[l,0].
\]
This system of equations for $\mbaT[0]$ in $D$ is closed with a no-flux boundary condition thanks to the matching condition \eqref{MCH5}. Thus, multiplying with $\mbT[k,0]$, integrating over $D$ and applying a Green's formula we see that 
\begin{equation} \label{eq:triH2}
 0 = - \int_D \sum_{k,l} \L[k,l]_0 (\phaT[0],\qqT[0]) \nabla_y \mbT[k,0] \cdot \nabla_y \mbT[l,0].
\end{equation}
Stating Assumption \ref{ass:AA3} more precisely, $\phaT[0](y)$ is assumed to be no corner of the Gibb's simplex in all points $y \in D$. Hence, by Assumption \ref{ass:AA4} on the kernel of $\{ \L[k,l]_0 (\phaT[0](y),\qqT[0](y)) \}_{k,l}$, $\p_{y_1} \mbaT[0](y)$ and $\p_{y_2} \mbaT[0](y)$ are multiples of $\bbb[1]$. Consequently, 
\[
 \JpT[k,-1] = 0.
\]

This fact simplifies the momentum balance \eqref{eq:PFMmombal} which, to leading order $-2$, reads
\begin{multline}
 0 = \nabla_y \cdot \Big{(} \eta \big{(} \nabla_y \vvT[0] + (\nabla_y \vvT[0])^\top \Big{)} \\
 + \sum_{i<j} \sigt[i,j](\qqT[0]) \Big{(} \nabla_y \big{(} \a[i,j] + \w[i,j] \big{)} - \sum_k \nabla_y \cdot \big{(} \nabla_y \phT[k,0] \otimes \p_{\nabla \ph[k]} \a[i,j] \big{)} \Big{)}, \label{eq:triH3}
\end{multline}
where the $\a[i,j]$ are evaluated in $(\phaT[0],\nabla_y \phaT[0])$ and $\eta$ and the $\w[i,j]$ in $\phaT[0]$. On noting that
\[
 \nabla_y \cdot \big{(} \nabla_y \phT[k,0] \otimes \p_{\nabla \ph[k]} \a[i,j] \big{)} = \nabla_y \phT[k,0] \nabla_y \cdot \p_{\nabla \ph[k]} \a[i,j] + \nabla_y^2 \phT[k,0] \p_{\nabla \ph[k]} \a[i,j]
\]
and that $\nabla_y \a[i,j] = \p_{\ph[k]} \a[i,j] \nabla_y \phT[k,0] + \nabla^2 \phT[k,0] \p_{\nabla \ph[k]} \a[i,j]$, the second line of \eqref{eq:triH3} yields
\begin{multline*}
 \sum_{i<j} \sigt[i,j](\qqT[0]) \Big{(} \nabla_y \big{(} \a[i,j] + \w[i,j] \big{)} - \sum_k \nabla_y \cdot \big{(} \nabla_y \phT[k,0] \otimes \p_{\nabla \ph[k]} \a[i,j] \big{)} \Big{)} \\
 = \sum_k \sum_{i<j} \sigt[i,j](\qqT[0]) \Big{(} \p_{\ph[k]} \a[i,j] + \p_{\ph[k]} \w[i,j] - \nabla_y \cdot \p_{\nabla \ph[k]} \a[i,j] \Big{)} \nabla_y \phT[k,0] = 0
\end{multline*}
thanks to \eqref{eq:triH1}. Applying the matching condition \eqref{MCH2} to $\vv$ and using that $\vvE[0]$ is constant across the interface layers (see \eqref{eq:inV0}) we see that $\nabla_y \vvT[0] + (\nabla_y \vvT[0])^\top$ vanishes on $\partial D$. Thus, multiplying \eqref{eq:triH3} with $\vvT[0]$ and integrating over $D$ we obtain that 
\[
 0 = - \int_D \eta | \nabla_y \vvT[0] + (\nabla_y \vvT[0])^\top |^2.
\]
Now, $\nabla_y \vvT[0] + (\nabla_y \vvT[0])^\top = 0$ implies that $\vvT[0]$ is a linear function of the form $\vvT[0](y) = A y + b$ with a skew symmetric matrix $A \in \R^{2 \times 2}$ and a vector $b \in \R^2$. However, as $\vvT[0]$ is constant along each of the three edges of $D$ implies that $A = 0$. Hence, $\vvT[0]$ is constant, and as for $\qqe[0]$ we can conclude that also the velocity 
\begin{equation} \label{eq:triR0}
\vve[0] \mbox{ is continuous in the whole domain.}
\end{equation}

Finally, considering \eqref{eq:PFMphi1} to order $-1$, 
\begin{equation} \label{eq:triH4}
 - \uu[\mIn,\mIp,\mIm] \cdot \nabla_y \phT[k,0] + \vvT[0] \cdot \nabla_y \phT[k,0] = - \nabla_y \cdot \JpT[k,0].
\end{equation}
Integrating over $D$, applying the divergence theorem, and applying the matching condition \eqref{MCH5} the right-hand side vanishes. But the vector $\nabla_y \phT[k,0]$ cannot vanish for all $k$ in $D$. Thus, necessarily $\vvT[0] = \uu[\mIn,\mIp,\mIm]$, that is \eqref{eq:triplevelcont} for $\vve[0]$ in two spatial dimensions.

Together with \eqref{eq:triR1} and \eqref{eq:triR2} the recovery of the triple junction equations \eqref{eq:triplevelcont2}--\eqref{eq:tripleforcebal2} of the sharp interface model is thus completed.

\subsection{Triple lines}
\label{sec:tri_3D}

Let us now discuss a triple line belonging to $\T[\mIn,\mIp,\mIm](t)$ in the three dimensional case. Proceeding similarly to the previous two-dimensional case we keep the arguments rather short and confine ourselves to highlight the differences. 

Given a point $\triple[\mIn,\mIp,\mIm](t) \in \T[\mIn,\mIp,\mIm](t)$ we denote the plane through $\triple[\mIn,\mIp,\mIm](t)$ orthogonal to $\T[\mIn,\mIp,\mIm](t)$ by $\Y[\mIn,\mIp,\mIm](t)$. The point is assumed to move in normal direction, i.e., 
\[
 \p_t \triple[\mIn,\mIp,\mIm](t) = \uu[\mIn,\mIp,\mIm](\triple[\mIn,\mIp,\mIm](t),t) \in \Y[\mIn,\mIp,\mIm](t).
\]
A local parametrisation or the triple line is denoted by $\paramtri[\mIn,\mIp,\mIm](s,t)$ with an arc-length parameter $s$ so that the unit tangential vector $\tangtri[\mIn,\mIp,\mIm](\paramtri[\mIn,\mIp,\mIm](s,t),t) := \p_s \paramtri[\mIn,\mIp,\mIm](s,t)$ is normal to $\Y[\mIn,\mIp,\mIm]$. In the following, we will write for the projection to $\Y[\mIn,\mIp,\mIm]$ 
\[
 \bbb[P]_{(\T[i,j,k])^\perp} =: \projY[\mIn,\mIp,\mIm] = \bbb[I] - \tangtri[\mIn,\mIp,\mIm] \otimes \tangtri[\mIn,\mIp,\mIm].
\]

For a point $x \in \Y[\mIn,\mIp,\mIm]$ close to $\triple[\mIn,\mIp,\mIm]$ we only rescale the coordinates in $\Y[\mIn,\mIp,\mIm]$ but not the coordinate $s$ along the triple line and define
\[
 y := \frac{x - \triple[\mIn,\mIp,\mIm](t)}{\e}, \quad x \in \Y[\mIn,\mIp,\mIm](t) \mbox{ close to the triple line}.
\]
We consider expansions of the form
\[
 \zeta_{\e}(x,t) = \hat{Z}_0(s,y,t) + \e \hat{Z}_1(s,y,t) + \e^2 \hat{Z}_2(s,y,t) + \dots
\]
and similarly for the fluxes $\Jp[k]$, $\Jq$, starting at the same order as before (see after \eqref{eq:triexp}). For the differential operators we note the transformation
 \begin{align*}
  \p_t \zeta(x,t) &= - \e^{-1} \nabla_y \hat{Z}(s,y,t) \cdot \uu[\mIn,\mIp,\mIm] + \mmm[O](\e^{0}), \\
  \nabla \zeta(x,t) &= \e^{-1} \nabla_y \hat{Z}(s,y,t) + \p_s \hat{Z}(s,y,t) \tangtri[\mIn,\mIp,\mIm] + \mmm[O](\e^{1}),
 \end{align*}
for a function $\zeta(x,t) = Z(s,y,t)$, where the triple line fields $\uu[\mIn,\mIp,\mIm]$ and $\tangtri[\mIn,\mIp,\mIm]$ are evaluated in $\triple[\mIn,\mIp,\mIm](t)$. Observe that the operator $\nabla_y$ now is the tangential gradient of the plane $\Y[\mIn,\mIp,\mIm](t)$, i.e., $$\nabla_y = \nabla_{\Y[\mIn,\mIp,\mIm](t)}.$$ 

The matching conditions \eqref{MCH1} and \eqref{MCH2} still are valid. However, the fluxes $\Jp[k]$ and $\Jq$ may exhibit contributions out of the plane $\Y[\mIn,\mIp,\mIm]$. Hence, \eqref{MCH3}--\eqref{MCH5} are only true for the tangential contributions, i.e., replacing $\JpE[k,\cdot]$ with $\projY[\mIn,\mIp,\mIm] \JpE[k,\cdot]$ and analogously for the other fluxes. 

The surfactant equation \eqref{eq:PFMsag1}, \eqref{eq:PFMsag2} to leading order yields $\nabla_y \qqT[0] = 0$ and, thus, the continuity of $\qqe[0]$. Also as before this leads to the triple junction condition \eqref{eq:triR1}, where we remark that all vector lie in $\Y[\mIn,\mIp,\mIm]$, so this force balance is indeed a condition intrinsic to the plane normal to the triple line. Subsequently, instead of \eqref{eq:triH1b} we now obtain
\[
 \JqT[-1] = - \sum_{i<j} \Mcc[i,j] ( \a[i,j] + \w[i,j] ) \big{(} \nabla_y \qqT[1] + \p_s \qqT[0] \tangtri[\mIn,\mIp,\mIm] \big{)},
\]
where the last term is normal to $\Y[\mIn,\mIp,\mIm]$. However, as the matching condition \eqref{MCH4} still is true for the in-plane contributions to the fluxes and vectors such as $\m[\mIn,\mIp,\mIm]$ are tangential to $\Y[\mIn,\mIp,\mIm]$ we can proceed as before and recover \eqref{eq:triR2}. 

The arguments around \eqref{eq:triH2} still are valid and lead to $\JpT[k,-1] = 0$. Also the conclusions after \eqref{eq:triH3} still are true whence the momentum equation \eqref{eq:PFMmombal} to leading order $-2$ becomes
\begin{align*}
 0 \, & = \nabla_y \cdot \Big{(} \eta \big{(} \nabla_y \vvT[0] + (\nabla_y \vvT[0])^\top \Big{)} \\
 & = \nabla_y \cdot \Big{(} \eta \big{(} \nabla_y (\projY[\mIn,\mIp,\mIm] \vvT[0]) + (\nabla_y (\projY[\mIn,\mIp,\mIm]\vvT[0]) )^\top \Big{)} \\
 & \quad + \nabla_y \cdot \Big{(} \eta \big{(} \tangtri[\mIn,\mIp,\mIm] \otimes \nabla_y (\vvT[0] \cdot \tangtri[\mIn,\mIp,\mIm]) + \nabla_y (\vvT[0] \cdot \tangtri[\mIn,\mIp,\mIm]) \otimes \tangtri[\mIn,\mIp,\mIm] \Big{)}, 
\end{align*}
where we split the field $\vvT[0]$ into contributions tangential and orthogonal to $\Y[\mIn,\mIp,\mIm]$. We now multiply the equation with $\projY[\mIn,\mIp,\mIm] \vvT[0]$ and integrate over $D$. Integrating by parts and applying the matching condition \eqref{MCH2} and \eqref{eq:inH4b} to $\vv$ to get rid of the boundary terms we obtain that
\begin{multline*}
 0 = \int_D \eta \big{|} \nabla_y (\projY[\mIn,\mIp,\mIm] \vvT[0]) + (\nabla_y (\projY[\mIn,\mIp,\mIm]\vvT[0]) )^\top \big{|}^2 \\
 + \eta \big{(} \tangtri[\mIn,\mIp,\mIm] \otimes \nabla_y (\vvT[0] \cdot \tangtri[\mIn,\mIp,\mIm]) \big{)} : \nabla_y (\projY[\mIn,\mIp,\mIm] \vvT[0]) \\
 + \eta \big{(} \nabla_y (\vvT[0] \cdot \tangtri[\mIn,\mIp,\mIm]) \otimes \tangtri[\mIn,\mIp,\mIm] \big{)} : \nabla_y (\projY[\mIn,\mIp,\mIm] \vvT[0]).
\end{multline*}
The last term vanishes as $\tangtri[\mIn,\mIp,\mIm]$ is normal and $\nabla_y$ yields a field tangential to $\Y[\mIn,\mIp,\mIm]$. Also the second last term vanishes as $\tangtri[\mIn,\mIp,\mIm]$ is independent of $y$ and orthogonal to $\projY[\mIn,\mIp,\mIm] \vvT[0]$:
\[
 \big{(} \tangtri[\mIn,\mIp,\mIm] \otimes \nabla_y (\vvT[0] \cdot \tangtri[\mIn,\mIp,\mIm]) \big{)} : \nabla_y (\projY[\mIn,\mIp,\mIm] \vvT[0]) = \nabla_y (\vvT[0] \cdot \tangtri[\mIn,\mIp,\mIm]) \cdot \nabla_y \big{(} (\projY[\mIn,\mIp,\mIm] \vvT[0]) \cdot \tangtri[\mIn,\mIp,\mIm] \big{)} = 0.
\]
Hence, proceeding as before \eqref{eq:triR0} we can conclude that $\projY[\mIn,\mIp,\mIm] \vvT[0]$ locally is constant. In equation \eqref{eq:triH4} we can replace $\vvT[0] \cdot \nabla_y \phT[k,0]$ with $\projY[\mIn,\mIp,\mIm] \vvT[0] \cdot \nabla_y \phT[k,0]$ and analogously show that
\[
 \uu[\mIn,\mIp,\mIm] = \projY[\mIn,\mIp,\mIm] \vvT[0].
\]
Together with \eqref{eq:triR1} and \eqref{eq:triR2} that, as discussed, both still are valid the triple junctions equations \eqref{eq:triplevelcont2}--\eqref{eq:tripleforcebal2} are recovered also in the case $d=3$.

\subsection{Quadruple points}
\label{sec:quad_expsol}

The asymptotic analysis around quadruple points is very similar to the analysis around triple points in the two-dimensional case (see Section \ref{sec:tri_sol}) whence we only provide a sketch. Around a quadruple point we consider a tetrahedron. The external unit-normals of its faces are tangential to the triple lines that form the quadruple point. A rescaled variable around the quadruple point as in \eqref{eq:tri_def_coord} is introduced and used for defining local expansions. The edge lengths of the tetrahedron scale with a power between $\frac{1}{2}$ and $1$ in $\e$. Thus, on the tetrahedron's faces the expansions match with those near the triple lines that are considered in the previous section \ref{sec:tri_3D}. 

The surfactant equation to leading order reads like \eqref{eq:triH0} in the tetrahedron around the quadruple point. It is closed with a no-flux boundary condition like \eqref{eq:triH0b} thanks to matching with the triple line solutions. Following the arguments after those two equations one can show again that the surfactant chemical potential $\qq$ is constant to leading order so that $\qqe[0]$ is continuous at the quadruple junction. This requires \eqref{eq:Q0const}, which is meant to be satisfied by Assumption \ref{ass:AA3}.

With $\qq$ being constant the surface tensions are constant to leading order, too. The geometry of the quadruple junction therefore is already fully determined by the force balances \eqref{eq:tripleforcebal2} (or \eqref{eq:triR1}) at the triple lines, see \cite{BroGarSto98}, Section 3, for a discussion. Expanding \eqref{eq:PFMmu}, which previously lead to solvability conditions at the interfaces and the triple junctions, therefore doesn't yield any more insight. As there is no surfactant mass flux along the triple lines there is also no need to discuss any higher order expansions of the surfactant equation.

It remains to show continuity of the velocity to leading order. But this can be done by expanding the momentum equation \eqref{eq:PFMmombal} and arguing as before equation \eqref{eq:triR0}. Moreover, the velocity coincides with the quadruple point velocity, which can be shown by considering the phase field equation \eqref{eq:PFMphi1} to leading order and proceeding as for a triple point (see around equation \eqref{eq:triH4}).

\subsection{Boundary conditions}
\label{sec:all_bc}

Expanding the boundary conditions \eqref{eq:PFMBCnoflux}, \eqref{eq:PFMBCnoforce}, and \eqref{eq:PFMBCqnoflux} of the phase field model by plugging in the outer expansions immediately yields the boundary conditions \eqref{eq:BCvel1}, \eqref{eq:BCvel2}, and \eqref{eq:BCsagbulk} for the velocity and bulk surfactant mass flux to leading order. 

Intersections of interfaces $\G[i,j]$ with $\p \Omega$ form triple junctions with the additional constraint that $\p \Omega$ is fixed but still can be analysed similarly, see \cite{OweRubSte90,GarNesSto98}. In analogy to \eqref{eq:triR1} the angle condition \eqref{eq:BCangle} is obtained as a local solvability condition by expanding \eqref{eq:PFMmu}. Also for the the surfactant equation \eqref{eq:PFMsag1}, \eqref{eq:PFMsag2} one can proceed as around \eqref{eq:triH1b} and \eqref{eq:triR2} to obtain the interface no-flux condition \eqref{eq:BCsagsurf}.

\section{Numerical simulation results}
\label{sec:numsim}

In this section we aim to support the results of the asymptotic analysis and showcase the capability of the new phase field model by some numerical computations. The results are on a qualitative level. We don't have any specific substances in mind but simply choose material parameters for convenience. The computational method is based on adaptive finite elements in space and a fractional step $\theta$--scheme in time, and we plan to provide details on it in future publications. Computations are carried out in the C++ finite element toolbox DUNE-FEM \cite{DedKloNolOhl10}, with use of the DUNE-Alugrid module \cite{AlkDedKloNol16} for construction of the adaptive parallel grids. For solving linear systems in parallel we made use of the flexible PETSc interface \cite{petsc-user-ref}, in particular the HYPRE BoomerAMG preconditioner \cite{FalYan02}. Visualisation is presented in the ParaView software \cite{book:Aya15}, graphs are constructed with GNUPlot. 1D computations were carried out in  MATLAB and credit is due to the MATLAB and Statistics Toolbox Release 2014a-2016b, The MathWorks, Inc., Natick, Massachusetts, United States.

\subsection{Surfactant diffusion through a triple junction}
\label{sec:econvq}

\begin{figure}[ht]
\vspace{-30pt}
        \centering
       \includegraphics[width=0.55\textwidth]{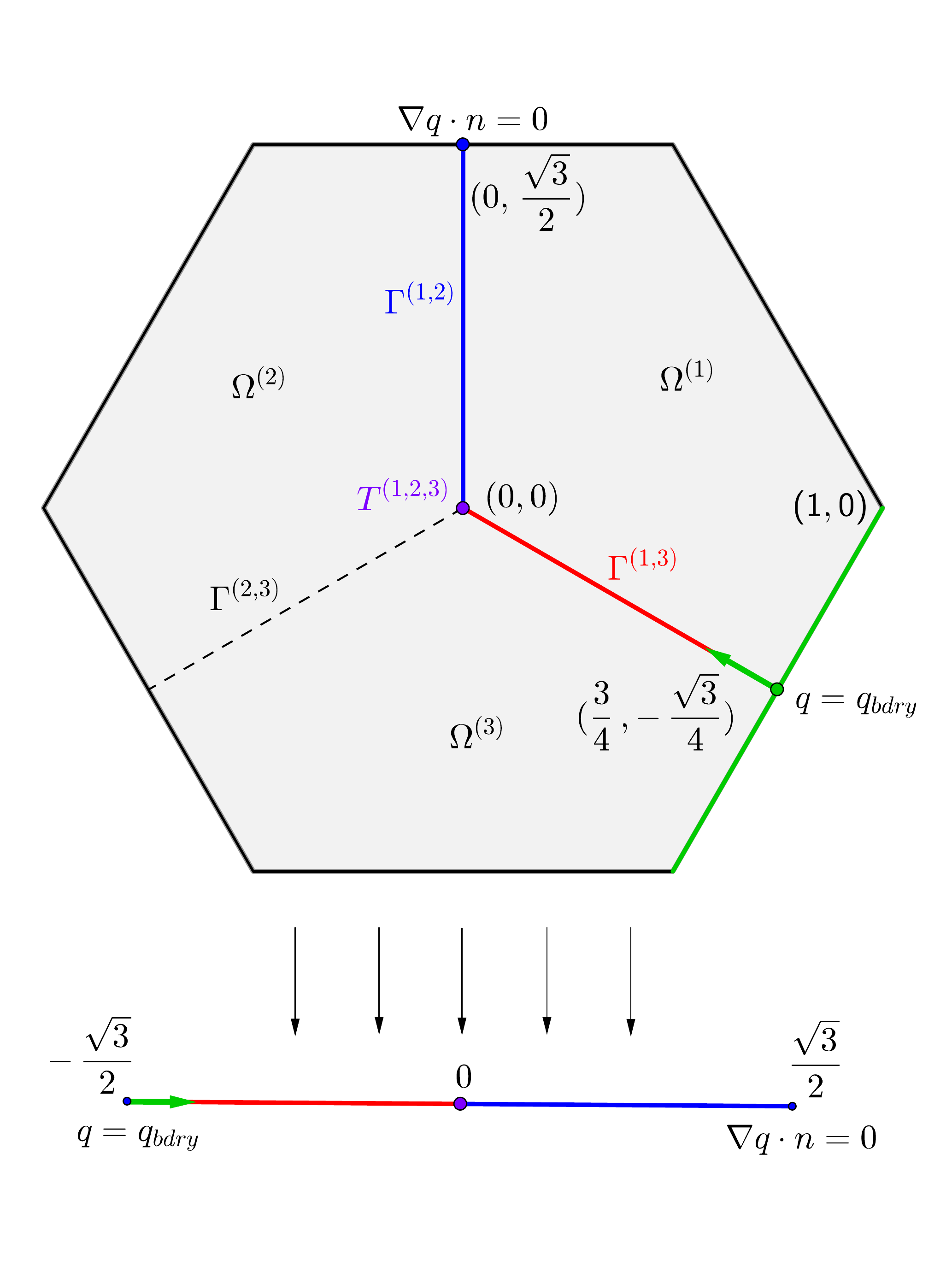}
       \vspace{-30pt}
        \caption{\em Setup for the $\e$-convergence test for the surfactant equation as considered in Section \ref{sec:econvq}.}
        \label{fig:Hex}
\end{figure}

We first discuss the convergence of the diffuse surfactant equation \eqref{eq:PFMsag1}-\eqref{eq:PFMsag2} to the sharp interface setting \eqref{eq:sagsurf2} in the case of a stationary network, i.e., no fluid flow and no interface motion. We consider three phases where the surfactant is present only at the interfaces $\G[1,2]$ and $\G[1,3]$ (we set $\cc[2,3](q)=0$ and $\c[i](q)=0, \ \forall i=1,2,3$). The idea is to supply surfactant on one of the boundaries and observe its diffusion along the interfaces and through a triple junction. The specific configuration is illustrated in Figure \ref{fig:Hex}. The domain is a regular hexagon of side length 1. It is comprised of 3 subregions $\O[i]$, $i=1,2,3$ separated by fixed straight interfaces $\G[i,j]$, $(i,j) = (1,2),(1,3),(2,3)$, which meet at a triple junction $\T[1,2,3]$ at the origin. 

We assume the following free energies for $(i,j) \in \{ (1,2), (1,3) \}$ that are quadratic in $q$ and thus lead to linear dependencies of the surfactant densities on the potential $q$ (see \eqref{q}):
\begin{equation}
 \gg[i,j](\cc[i,j](q))-\sigma_0 = \frac{1}{2}\beta_{i,j}(\cc[i,j])^2=\frac{1}{2}\frac{q^2}{\beta_{i,j}}, \label{eq:toyint}
\end{equation}
with some constants $\sigma_0,\beta_{i,j}>0$, and then $\c[i,j](q) = q/\beta_{i,j}$, and $\sigt[i,j](q) = \sigma_0 - q^2/(2\beta_{i,j})$.
For the test series, we choose the model parameters $\beta_{1,2}=4,\ \beta_{1,3}=1$. Moreover, we assume constant mobilities $\Mcc[1,2]=25,\ \Mcc[1,3]=100$. Initially, the surfactant is absent ($q(x,0) = 0$ $\forall x \in \Omega$). It is then supplied at the intersection point of $\G[1,3]$ with $\partial \Omega$ by imposing a Dirichlet condition with
\[
 q_{bdry}(t) = 
 \begin{cases} 
  5000t, & \quad \mbox{if } 0 \leq t \leq 10^{-4}, \\
  0.5, & \quad \mbox{if } t \geq 10^{-4}.
 \end{cases}
\]
Otherwise, we impose the natural (homogeneous Neumann) boundary condition. 

\begin{figure}[ht]
        \centering
        \includegraphics[width=0.8\textwidth]{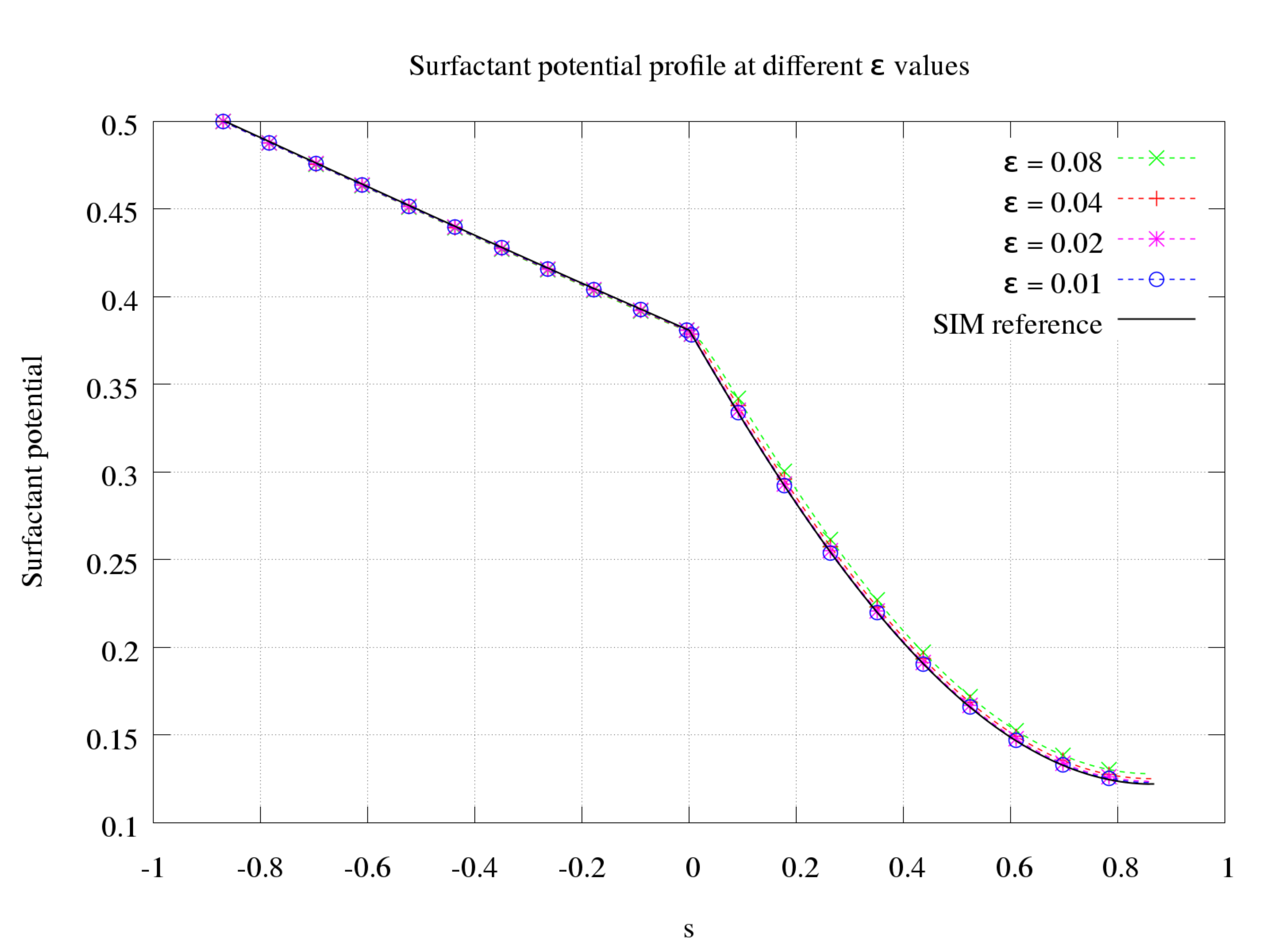}
        \caption{\em Results for the problem in Section \ref{sec:econvq}: Profiles of $\qq$ for different $\e$ values and of the solution $q$ to \eqref{eq:surf1Dtest} at time $T$. The fields $\qq$ were sampled along the interfaces and mapped onto the interval $\big(-\frac{\sqrt{3}}{2},\frac{\sqrt{3}}{2}\big)$ as illustrated in Figure \ref{fig:Hex}.} 
        \label{fig:econvhex_mb}
\end{figure}

\begin{figure}[ht]        
        \centering
        \includegraphics[width=0.8\textwidth]{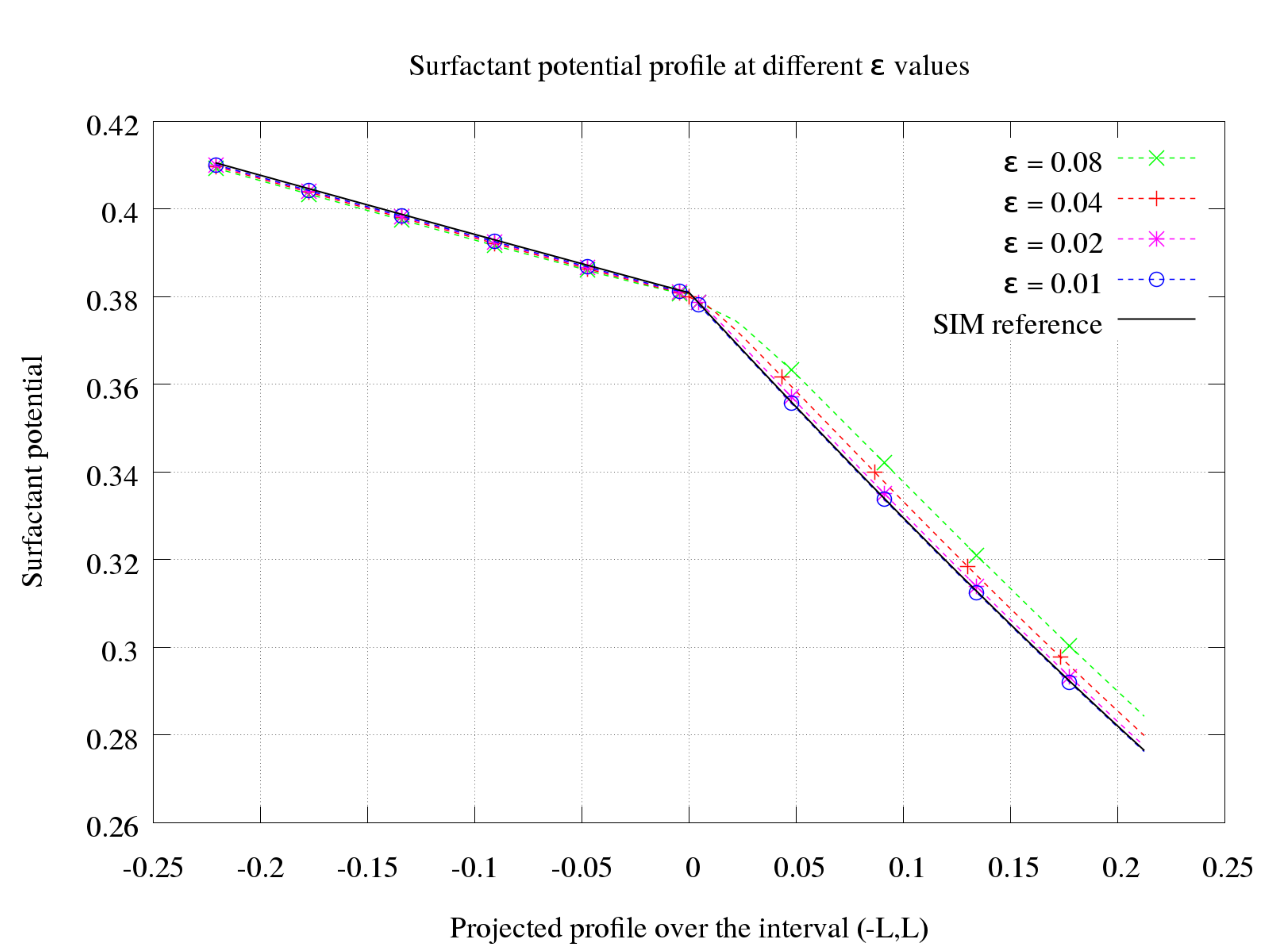}
        \caption{\em A zoom into Figure \ref{fig:econvhex_mb} around the triple junction, which is located at $s=0$.}
        \label{fig:econvhex_mb250}
\end{figure}

For $L=\sqrt{3}/2$, we map the equations for $\cc[1,2](q)$ on $\G[1,2]$ onto the interval $(0,L)$ and $\cc[1,3](q)$ on $\G[1,3]$ onto the interval $(-L,0)$ as sketched in Figure \ref{fig:Hex}. The triple junction $\T[i,j,k]$ is thus mapped to $0$ and boundary intersections to $\pm L$. Application of the above transformation, conditions and assumptions to the equations \eqref{eq:sagsurf2} and \eqref{eq:sagtripmasscond2}, result in the following problem:
\begin{align}
 \p_t \cc[1,3](q(s,t)) - \Mc[1,3] \p_{ss}q(s,t) &= 0, \qquad \forall(s,t) \in (-L,0) \times (0,T), \nonumber \\
 \p_t \cc[1,2](q(s,t)) - \Mc[1,2] \p_{ss}q(s,t) &= 0, \qquad \forall(s,t) \in (0,L) \times (0,T),\nonumber \\
  q(-L,t) = q_{\text{bdry}},\hspace{15pt} \p_sq(L,t) &=0, \qquad \forall t \in (0,T), \label{eq:surf1Dtest} \\
   [\Mc[\cdot,\cdot]q']^+_- (0,t) = 0, \hspace{13pt} [q]^{+}_{-}(0,t) &= 0, \qquad \forall t \in (0,T), \nonumber \\
  q(s,0) &= 0,\qquad \forall s\in[-L,L]. \nonumber
\end{align}
With standard finite difference techniques a sufficiently accurate approximation to the solution is obtained. The errors are evaluated at time $T=0.01$, for which significant gradients still are present. Figure \ref{fig:econvhex_mb} gives an impression of the solution at time $T$. A closer view of the triple junction region is displayed in Figure \ref{fig:econvhex_mb250}. 

The approximation of the solution to the diffuse interface model involves two stages. First, we initialise the phase fields by smoothed characteristic functions of the three domains. The initial data are relaxed to a diffuse triple junction in equilibrium by solving the Cahn-Hilliard system \eqref{eq:PFMphi1}-\eqref{eq:PFMmombal} ($\mathcal{L}^{i,j}$, $\a[i,j]$, $\w[i,j]$ as in Section \ref{sec:BLM}) in the absence of fluid flow ($\vv = 0$) and surfactant ($\qq = 0$, $\sigt[i,j]=1$) until a stationary state is reached. The second stage involves fixing the final computed phase field solutions, and then substituting them into equations \eqref{eq:PFMsag1}-\eqref{eq:PFMsag2}. These are then solved with conditions that approximate the setting of the sharp interface model test. We set $\xi_i=0$ and $\d[2,3]=0$, and the surfactant is supplied at the boundary $\p\Omega_{\text{in}}\coloneqq\big\{(r,\sqrt{3}/(2(r-1))\big| r\in(0,1)\big\}$, using the Dirichlet data $q_{\text{bdry}}(t)$. On all other boundaries a (natural) homogeneous Neumann boundary condition is imposed. 

We replace the distributions $\d[1,2],\d[1,3]$ with a regularisation 
\begin{equation}
 \tilde{\delta}_{i,j}(\pha,\npha)=\begin{cases}
                       \d[i,j](\pha,\npha), & \hspace{15pt} \text{if }|\d[i,j](\pha,\npha)|>C\e^2,\\
                       C\e^2, & \hspace{15pt} \text{otherwise}.
                      \end{cases}\label{eq:regdelta}
\end{equation}
As the $\d[i,j]$ decay exponentially fast outside of the interfacial regions, this bulk degeneracy would otherwise cause numerical instability. We found $C=0.001$ sufficient for the comparison of different values of $\e$. A similar technique of regularisation in the case of degeneracy in the bulk is presented in \cite{RatRibVoi06}. 

The discretisation parameters in space and time are chosen to ensure sufficient accuracy to correctly observe the $\e$-convergence. At the final time $t=T$, we sample $\qq$ at $N=400$ equi-distributed points along the straight segments representing $\G[1,3] \cong (-L,0)$ and $\G[1,2] \cong (0,L)$ for the comparison with $q$. For a sample $y_\e=(y_k)_{k=1}^N$, $y_k = (\qq)_k - q_k$, we compute the errors 
\begin{equation}\label{eq:lpdef} 
 \| y_\e \|_{l^\infty} = \max_{k=1,\dots,N}(y_k),\qquad \| y_\e \|_{l^2} = \Big(\sum_{k=1}^N y_k^2\Big)^{\frac{1}{2}},
\end{equation}
and for a series of samples $\{ y_\e \}_{\e}$ we estimate the order of convergence in a norm $\|\cdot\|_*$ by
\begin{equation}\label{eq:EOCdef}
 *\text{-EOC}(\e_1,\e_2) = \log(\| y_{\e_1} \| / \| y_{\e_2} \|) / \log(\e_1 / \e_2). 
\end{equation}

\begin{table}[ht]
\centering
\begin{tabular}{|c|c|c|c|c|c|}\hline
$\e$ &$q_{\e }(L)$& $\|y_{\e }\|_{l^\infty}$ & $l^\infty$-EOC & $\|y_{\e }\|_{l^2}$ & $l^2$-EOC \\ \hline 
0.08 & 0.127939 & 0.005887 & 0.910448 & 0.005596 & 0.865238   \\ \hline
0.04 & 0.125184 & 0.003132 & 0.927244 & 0.003072 & 0.882334 \\ \hline
0.02 & 0.123699 & 0.001617 & 0.919096 & 0.001667 & 0.852642 \\ \hline
0.01 & 0.122923 & 0.000871 & --        & 0.000923 & -- \\ \hline
ref & 0.122052  &  --     & --&  --     & --\\ \hline
\end{tabular}
\caption{\em Errors and convergence rates for the problem in Section \ref{sec:econvq}.}
\label{tab:econvhex_mbdiff}
\end{table}

The profiles of the surfactant potential $\qq$ at time $T$ for different values of $\e$ and the sharp interface model solution are displayed in Figure \ref{fig:econvhex_mb} with a zoom into the triple junction in Figure \ref{fig:econvhex_mb250}. In Table \ref{tab:econvhex_mbdiff} we present the errors and EOCs, where we remark that the maximum error occurs at points $s=L$ furthest from the source (Dirichlet) boundary condition. 

We observe good agreement of the solution across all $\e$ values ($4.82\%$ in $l^\infty$ for $\e=0.08$ and $0.71\%$ at $\e=0.01$). The estimated orders of convergence at around $0.9$ in $l^\infty$ and $0.85$ in $l^2$ clearly indicate convergence, which seems slightly sub-linear. Figure \ref{fig:econvhex_mb250} demonstrates that the greatest source of inaccuracy of the model arises from the jump at $s=0$ in the gradient $\p_s q$.

\subsection{Angles at a triple junction}\label{sec:angles}

\begin{figure}[ht]
        \centering
        \includegraphics[width=0.3\textwidth]{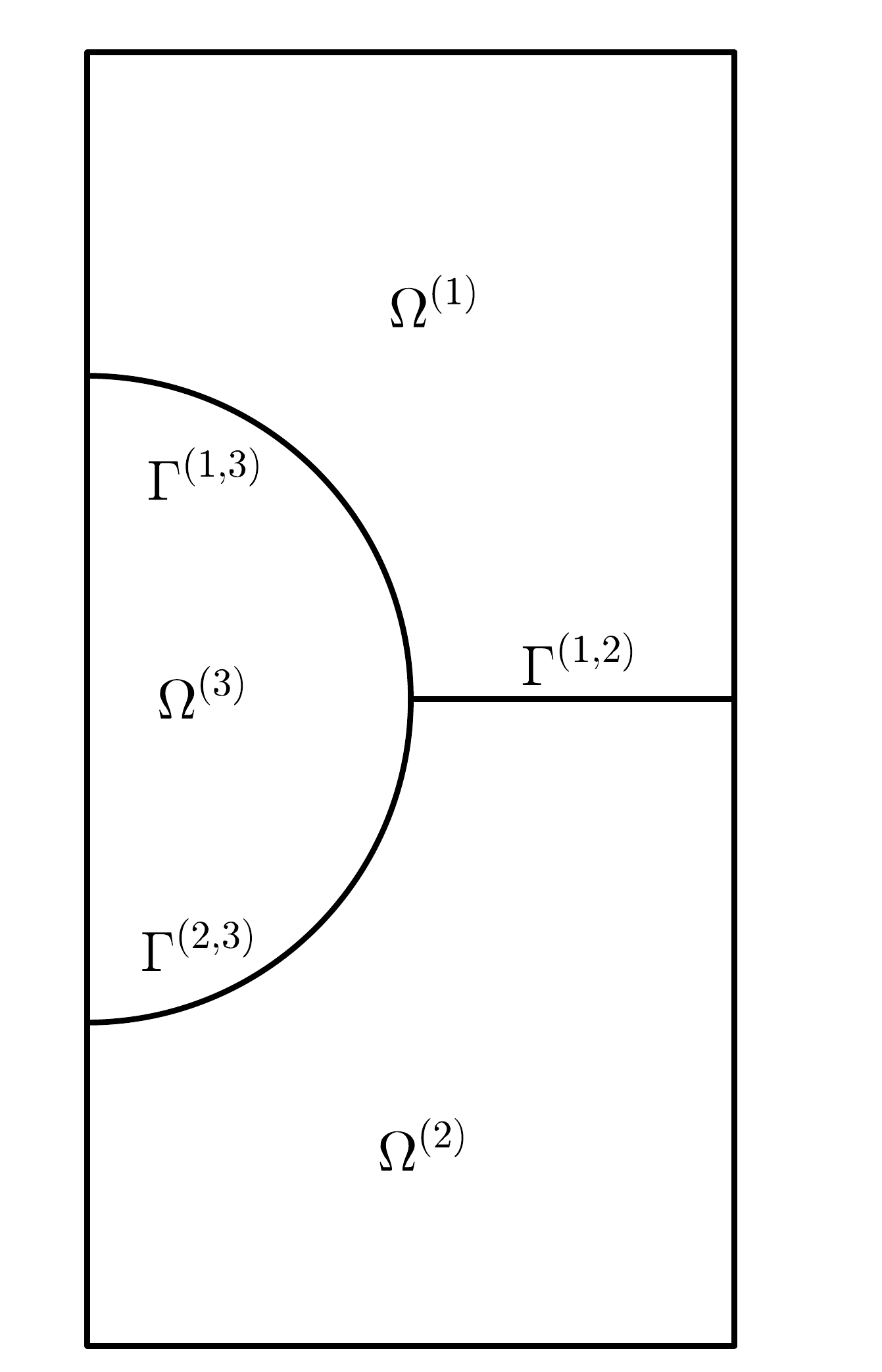}
        \hfill 
        \includegraphics[width=0.53\textwidth]{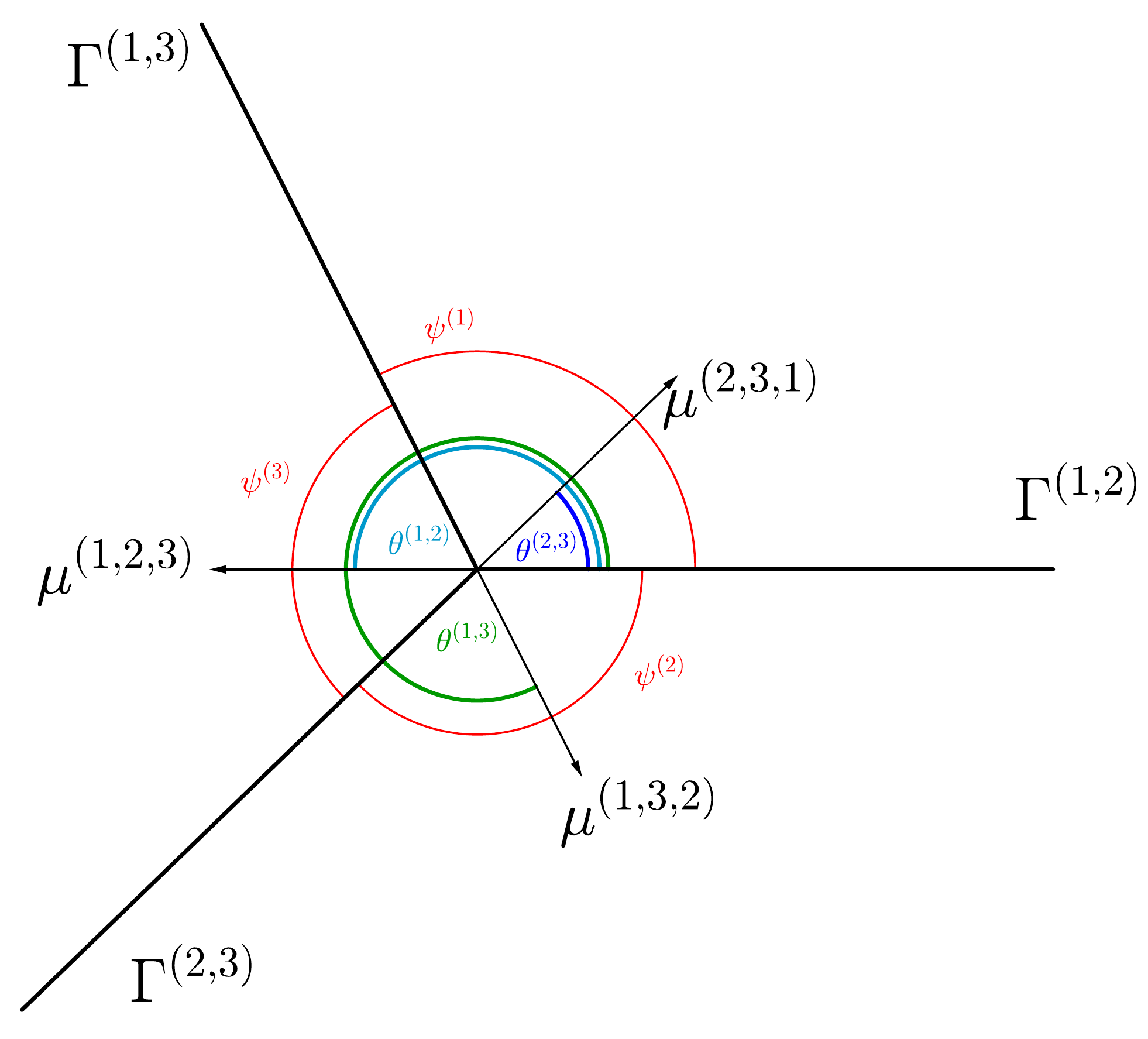}
        
        \caption{\em Left: Initial configuration for the problem in Section \ref{sec:angles}. Right: Diagram of a triple junction.
        The angles $\psi^{(k)}$ are the angles between the hypersurfaces $\G[j,k],\G[k,i]$, and the angles $\t^{(i,j)}$ can be used to express the mechanical force balance, see \eqref{eq:angletrig}. 
        }
        \label{fig:LensInitHalf} \label{fig:tripleangles}
\end{figure}

The angles at a triple junction are determined by the mechanical equilibrium \eqref{eq:tripleforcebal} of the surface tensions, which depend on the surfactant densities and thus can change over time. We here start with a three phase half-lens setting as illustrated in Figure \ref{fig:LensInitHalf} (left). It can be viewed as a fluid droplet trapped between two fluids with a line of symmetry in the centre of the droplet. We then relax the configuration whilst supplying surfactant by a Dirichlet boundary condition. The material parameters are such that, without surfactant, angles  of $2\pi/3$ (120 degrees) form at the triple junction, while at full surfactant saturation the equilibrium angles are $\pi/2,\ 2\pi/3$, and $5\pi/6$ (90, 120, and 150 degrees) in the sharp interface model. 

The domain is given by $\Omega=(-0,2)\times(-2,2)$, and time $t\in(0,T)$ with end time $T=10$. Initially, $\G[1,2](0) = \{ (z,0) | z \in (r,2) \}$, and $\G[1,3](0)$ (and $\G[2,3](0)$) is given by the open upper (resp. lower) right quarter-circle of radius $r$ centred at the origin. The triple junction $\T[1,2,3](0)$ is located at $(r,0)$. In this test we take $r=\sqrt{3/\pi}$, giving $|\O[3]|=3/2$.

We may pick coordinates so that the triple junction is located at the origin and the $\G[1,2]$ interface is along the positive $x$ axis. The angle formed in phase $\O[k]$ is denoted by $\phi^{(k)}$. We can also define $\theta^{(i,j)}$ to be the angles anticlockwise from $\G[1,2]$ to the co-normals $\m[i,j,k]$, see  Figure \ref{fig:tripleangles} (right; note that always $\theta^{(1,2)} = \pi$). Then \eqref{eq:tripleforcebal} can be written as 
\begin{equation} \label{eq:angletrig}
 \begin{cases}
   \sigt[1,2](q) \cos(\theta^{(1,2)}) + \sigt[1,3](q) \cos(\theta^{(1,3)}) + \sigt[2,3](q) \cos(\theta^{(2,3)}) &= 0, \\
   \sigt[1,2](q) \sin(\theta^{(1,2)}) + \sigt[1,3](q) \sin(\theta^{(1,3)}) + \sigt[2,3](q) \sin(\theta^{(2,3)}) &= 0. 
 \end{cases}
\end{equation}
We choose surface free energies as in \eqref{eq:toyint} (for all interfaces) and recall that then $\sigt[i,j](q) = \sigma_0 - q^2/(2\beta_{i,j})$. We take $\sigma_0 = 4$, $\beta_{1,2}=1/24$, $\beta_{1,3} = 1/(8(4-\sqrt{3}))$, and $\beta_{2,3} = 1/16$. If $q = 0$ then the surface tensions are the same, $\sigt[i,j](0) = 4$, and the equilibrium angles then are $\psi^{(k)} = 2\pi/3$ or, equivalently, $\theta^{(1,2)} = \pi$, $\theta^{(1,3)} = 5\pi/3$, $\theta^{(2,3)} = \pi/3$. If $q = 0.5$ then $\sigt[1,2](0.5) = 1$, $\sigt[1,3] = \sqrt{3}$, and $\sigt[2,3](0.5) = 2$, which means that $\theta^{(1,2)} = \pi$, $\theta^{(1,3)} = 3\pi/2$, $\theta^{(2,3)} = \pi/3$.  Hence, the equilibrium angles are
\begin{equation} \label{eq:angles_target}
 \psi^{(1)} = \pi/2, \, \psi^{(2)} = 2 \pi /3, \, \psi^{(3)} = 5 \pi/6 \qquad \mbox{if } q = 0.5.
\end{equation}
We choose the bulk free energies as 
\begin{equation} \label{eq:toybulk}
 \g[i](\c[i](q)) = \beta_i(\c[i](q))^2 / 2 = q^2 / (2\beta_i) 
\end{equation}
with $\beta_i =1$. We also take constant mobilities $\Mc[i]=100$ and $\Mcc[i,j]=100/\beta_{i,j}$.

Initially, there is no surfactant present in the domain ($q(x,0) = 0$ $\forall x \in \Omega$). We impose the following boundary condition on $z_R \in \{ (2,y) | y \in[-2,2] \}$:
\[
q(z_R,t)=
\begin{cases}
0, & \text{for } t\in (0,T_0),\\
\frac{t-T_0}{2T_q}, & \text{for } t\in (T_0,T_0+T_q),\\
0.5, & \text{for } t\in (T_q,T),
\end{cases}
\]
where $T_0 = 1$ and $T_q = 1$. The surfactant diffuses into the whole domain and approaches $q=0.5$, which is achieved before the final time $T$. 

For the dynamics of the geometry, we consider the model in Section \ref{sec:sumSIM} but neglect the fluid flow, so that the equations \eqref{eq:massbal2}, \eqref{eq:mombal2}, \eqref{eq:intvelcont2}, \eqref{eq:triplevelcont2} disappear and $\v = 0$. Regarding the diffuse interface approximation, the phase fields are initialised by using the leading order profiles from the asymptotic analysis along two phase interfaces: Given $\e$, for $z=(x,y)\in\Omega$ we set
\[
\ph[1](x,y,t=0) = \begin{cases}
                        \frac{1}{2}(1+\tanh(\frac{2y}{\e})), & \quad\text{for}\quad (x,y)\in(r,2]\times [-2,2],\\
                        \frac{1}{2}(1+\tanh(\frac{2(r-1)}{\e})), & \quad\text{for}\quad (x,y) \in [-r,r]\times(0,2],\\
                        0, & \quad\text{otherwise},
                     \end{cases}
\]
\[
\ph[2](x,y,t=0) = \begin{cases}
                        \frac{1}{2}(1+\tanh(\frac{-2y}{\e})), & \quad\text{for}\quad (x,y)\in(r,2]\times[-2,2],\\
                        \frac{1}{2}(1+\tanh(\frac{2(r-1)}{\e})), & \quad\text{for}\quad (x,y) \in [-r,r]\times[-2,0],\\
                        0, & \quad\text{otherwise},
                     \end{cases}
\]
and $\ph[3] = 1-\ph[1]-\ph[2]$ to guarantee \eqref{eq:phiconstraint}. Due to the absence of the flow we only have to solve the system \eqref{eq:PFMphi1}--\eqref{eq:PFMsag2}. We take the phase field model presented in Section \ref{sec:BLM} for the Cahn-Hilliard potentials with a constant mobility parameter $M_c=0.1$ and a  regularisation parameter $\Lambda=0.1$. Recall that this regularisation is required to avoid third phase contributions along the interfacial layers, but let us remark that choosing $\Lambda$ as small as possible is desirable as high values are detrimental to the recovery of the angles at the triple junction. 

\begin{figure}[ht]
        \centering
        \includegraphics[height=7cm]{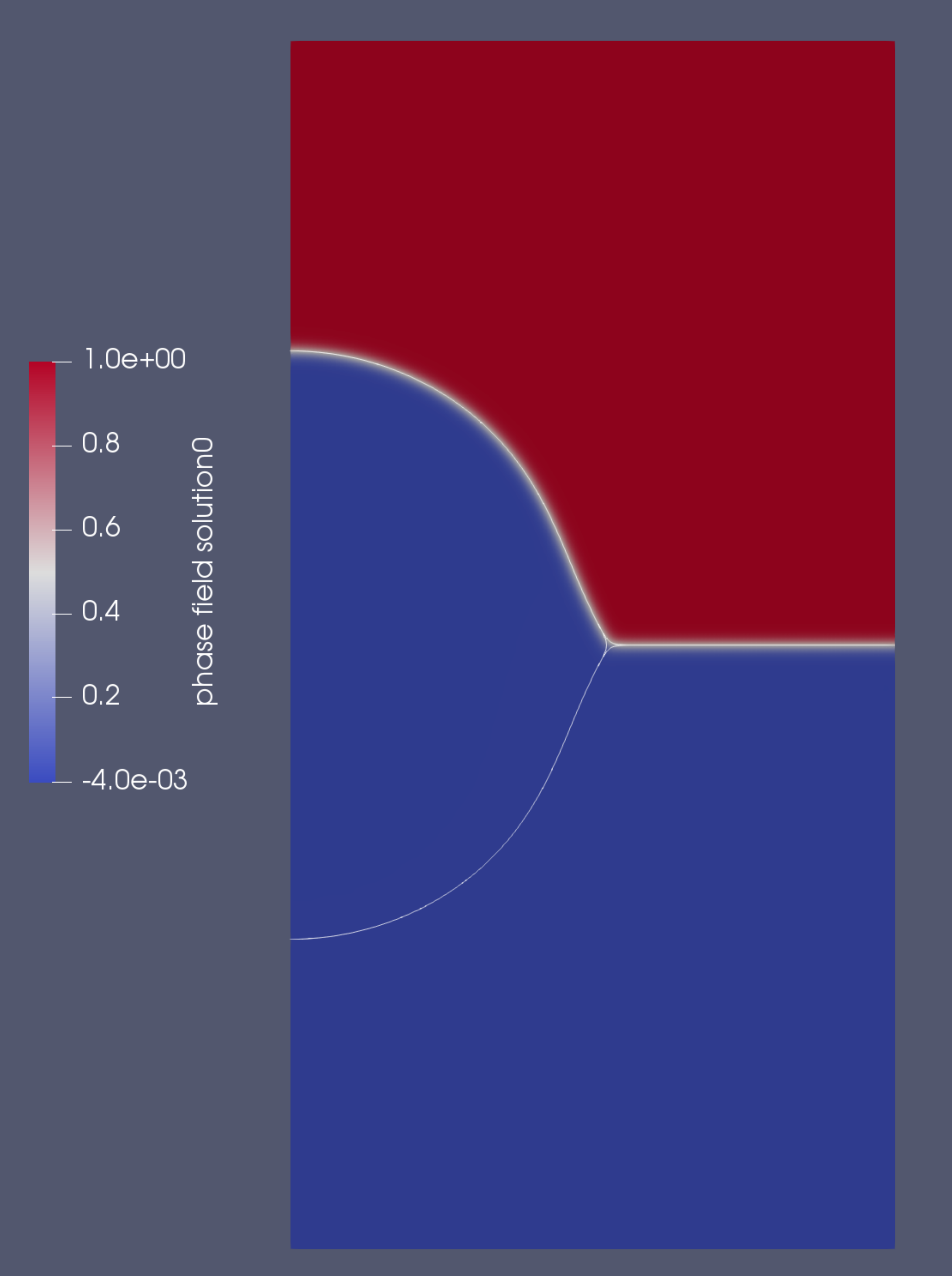}
        \hfill
        \includegraphics[height=7cm]{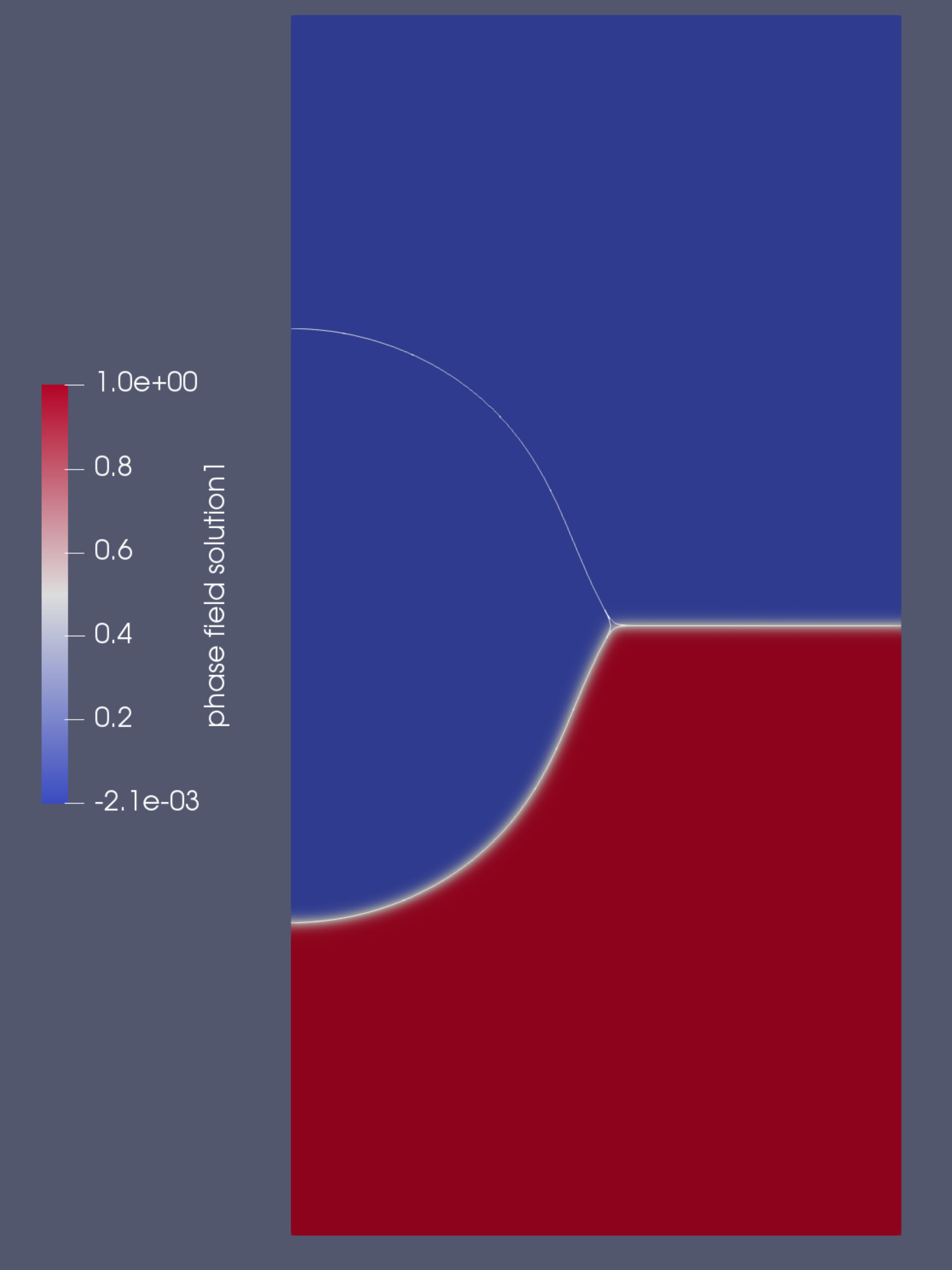}
        \hfill
        \includegraphics[height=7cm]{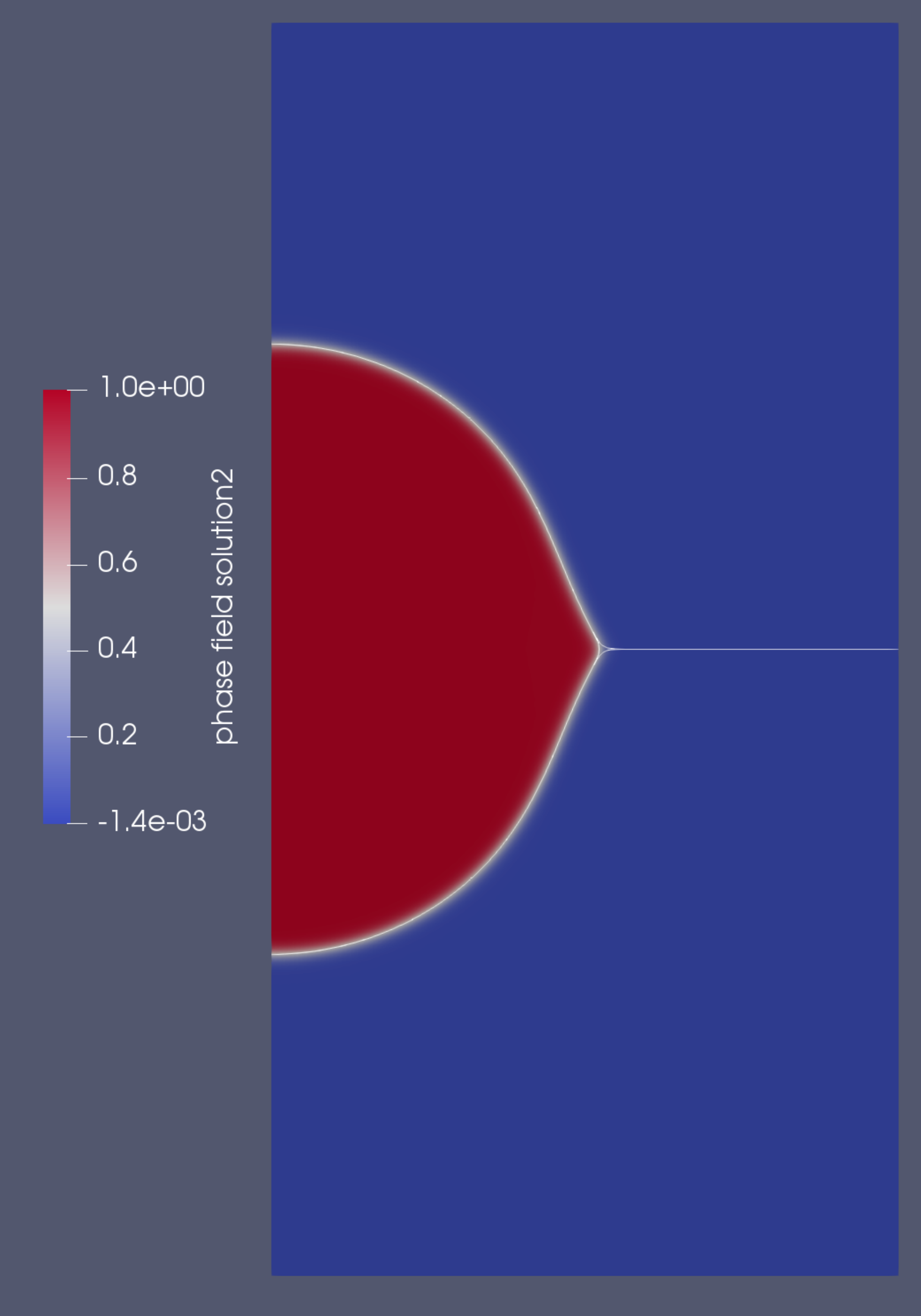}
        \caption{\em Phase fields at time $T_0=1$ (left $\ph[1]$, centre $\ph[2]$, right $\ph[3]$) with all $\ph[k]=1/2$ level sets in each plot for $\e=0.05$ prior to the introduction of surfactant.}
        \label{fig:LensPreSurf}
\end{figure}

\begin{figure}[ht]
        \centering
        \includegraphics[width=0.47\textwidth]{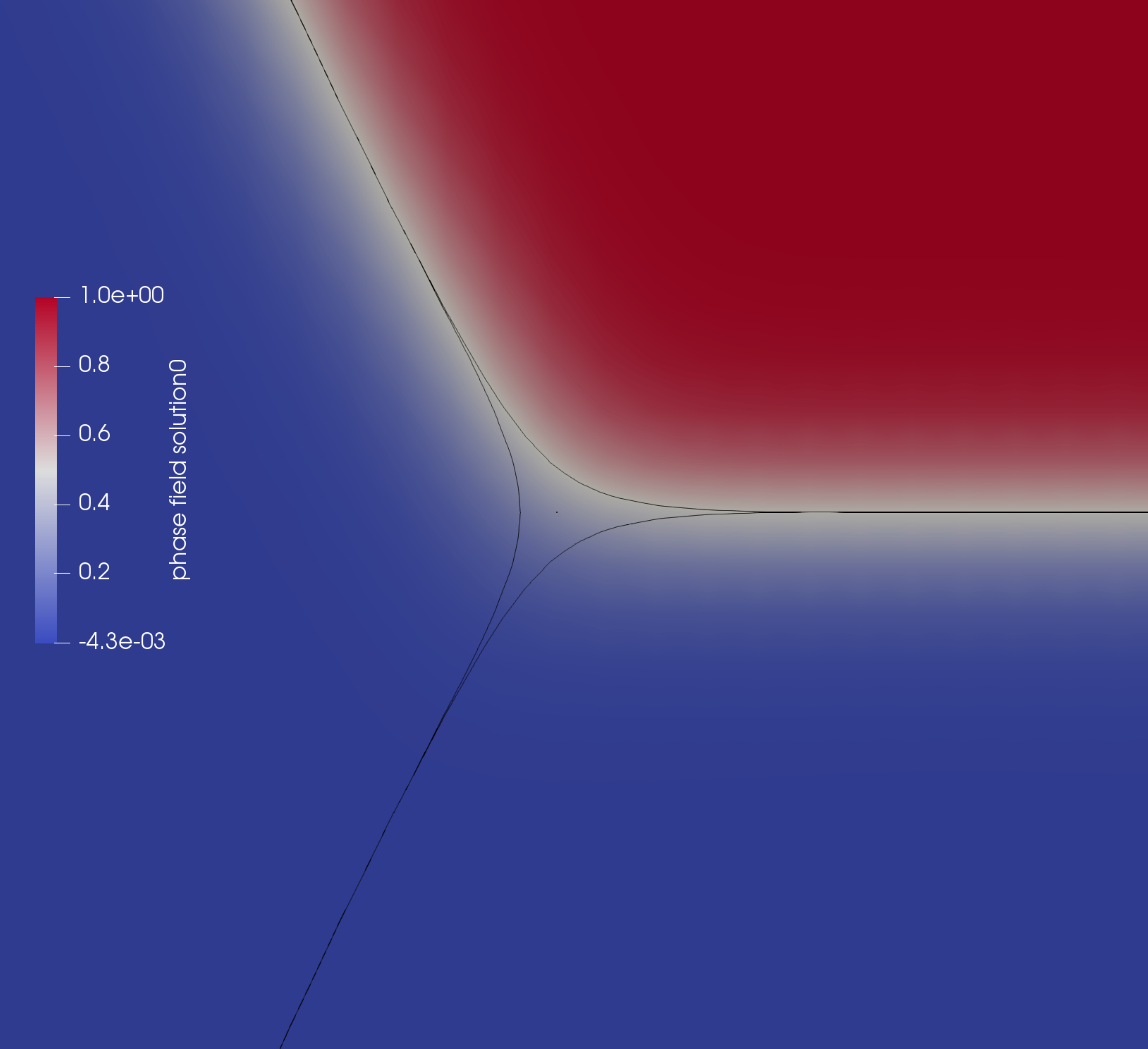}
        \hfill
        \includegraphics[width=0.47\textwidth]{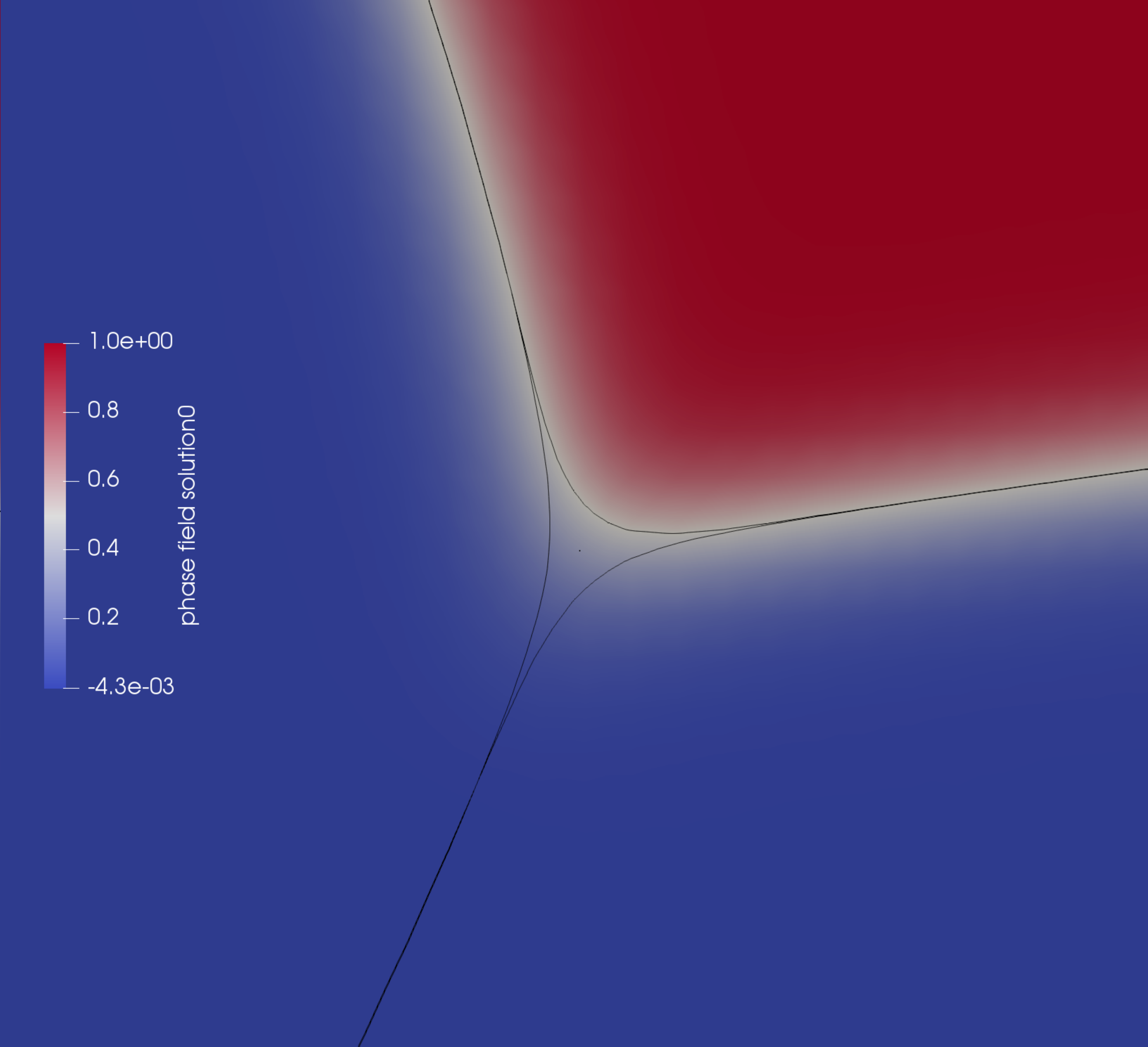}
        \hfill
        \caption{\em Close-up of the triple junction. The colour indicates the phase $\ph[1]$. Left: at time $T_0=1$ prior to the introduction of surfactant. Right: at time $T=10$, when $q$ is almost constant and the junction relaxed. We include the $\ph[k]=1/2$ level sets on each plot.}
        \label{fig:BeforeAfter}
\end{figure}

Until time $T_0$ the angles relax to nearly $2\pi/3$, of which Figure \ref{fig:LensPreSurf} gives an impression. At that time the surfactant is supplied on the boundary and starts to diffuse in. Subsequently, the angles at the triple junction change. We display a closeup of the final angle at time $T$ side by side with a closeup of the angle at time $T_0$ for comparison in Figure \ref{fig:BeforeAfter}. 

The angles $\psi^{(\cdot)}$ are measured at the junction as follows:
\begin{enumerate}[1.]
 \item First, we find the triple junction, where $\ph[k] = 1/3$ for all $k$, by looping over the elements and monitoring sign changes of the values $\ph[k] - 1/3$ in the vertices. 
 \item For every triple $(i,j,k) \in \{(1,2,3),(1,3,2),(2,3,1)\}$ and 20 evenly spaced values $\eta$ in $(9/20,1/2)$ we proceed similarly to find points, called $\eta$-junctions below, where $\ph[i] = \ph[j] = \eta$ and $\ph[k]=1-2\eta$. 
 \item For every triple $(i,j,k) \in \{(1,2,3),(1,3,2),(2,3,1)\}$ we perform two different linear regressions between the $\eta$-junctions.
 For the first we restrict the regression line to pass through the triple junction ({\em anchored}) $\ph[1]=\ph[2]=\ph[3]=1/3$, for the second method we do not use the triple junction to enforce any constraint ({\em unanchored}).
 \item We use the directions of these lines to compute the approximate angles between each pair of interfaces for both the anchored (labelled $\psi^{(\cdot)}_a$) and the unanchored (labelled $\psi^{(\cdot)}_u$) regression.
\end{enumerate}

\begin{figure}[ht]
        \centering
        \includegraphics[width=0.93\textwidth]{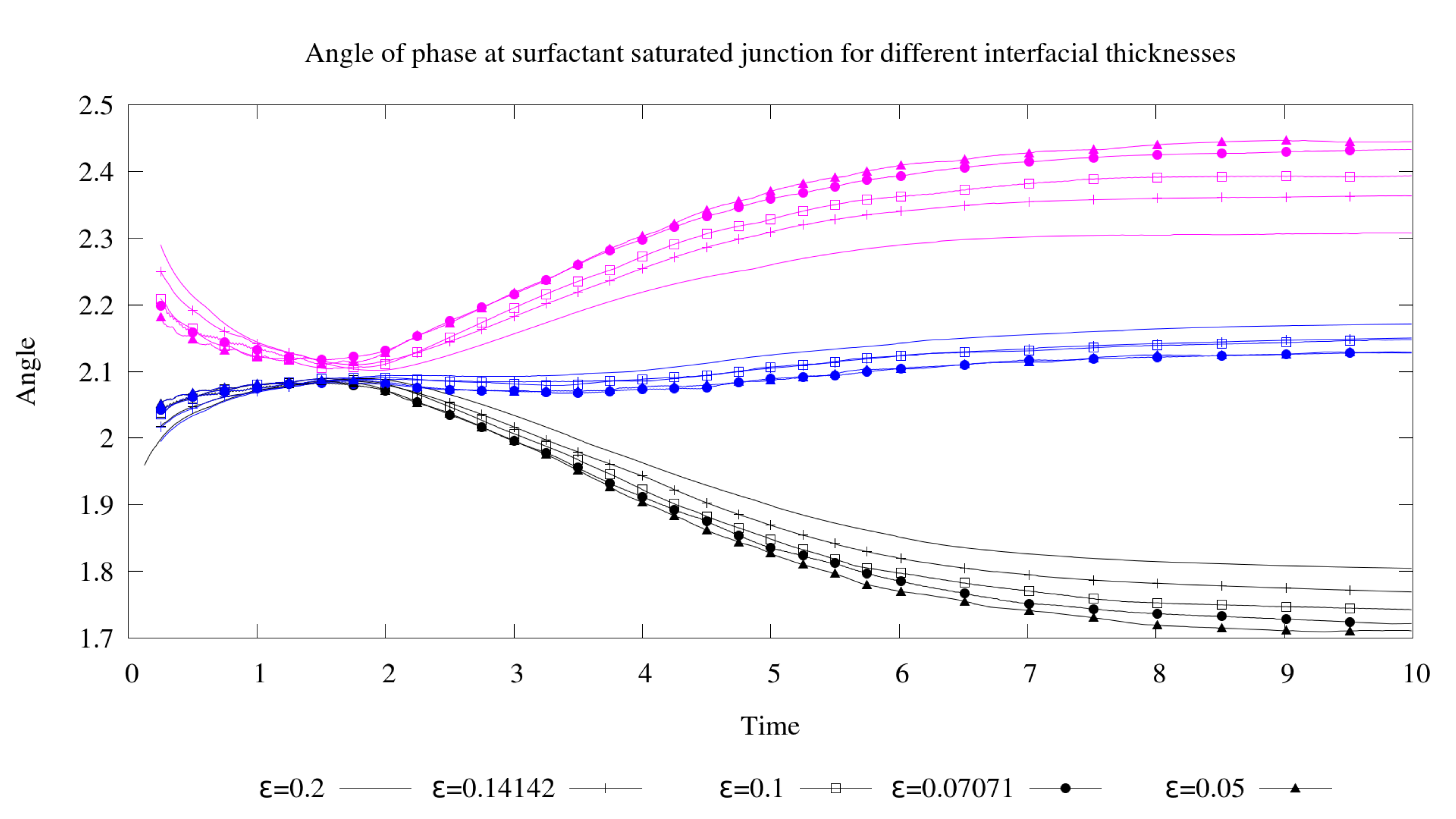}
        \caption{\em Evolution of measured angles $\psi^{(\cdot)}_a$ (anchored regression) through triple junction (black: $\ph[1]$, blue: $\ph[2]$, pink: $\ph[3]$) for different values of $\e$.}
        \label{fig:TJplot}
\end{figure}

\begin{figure}[ht]
        \centering
        \includegraphics[width=0.93\textwidth]{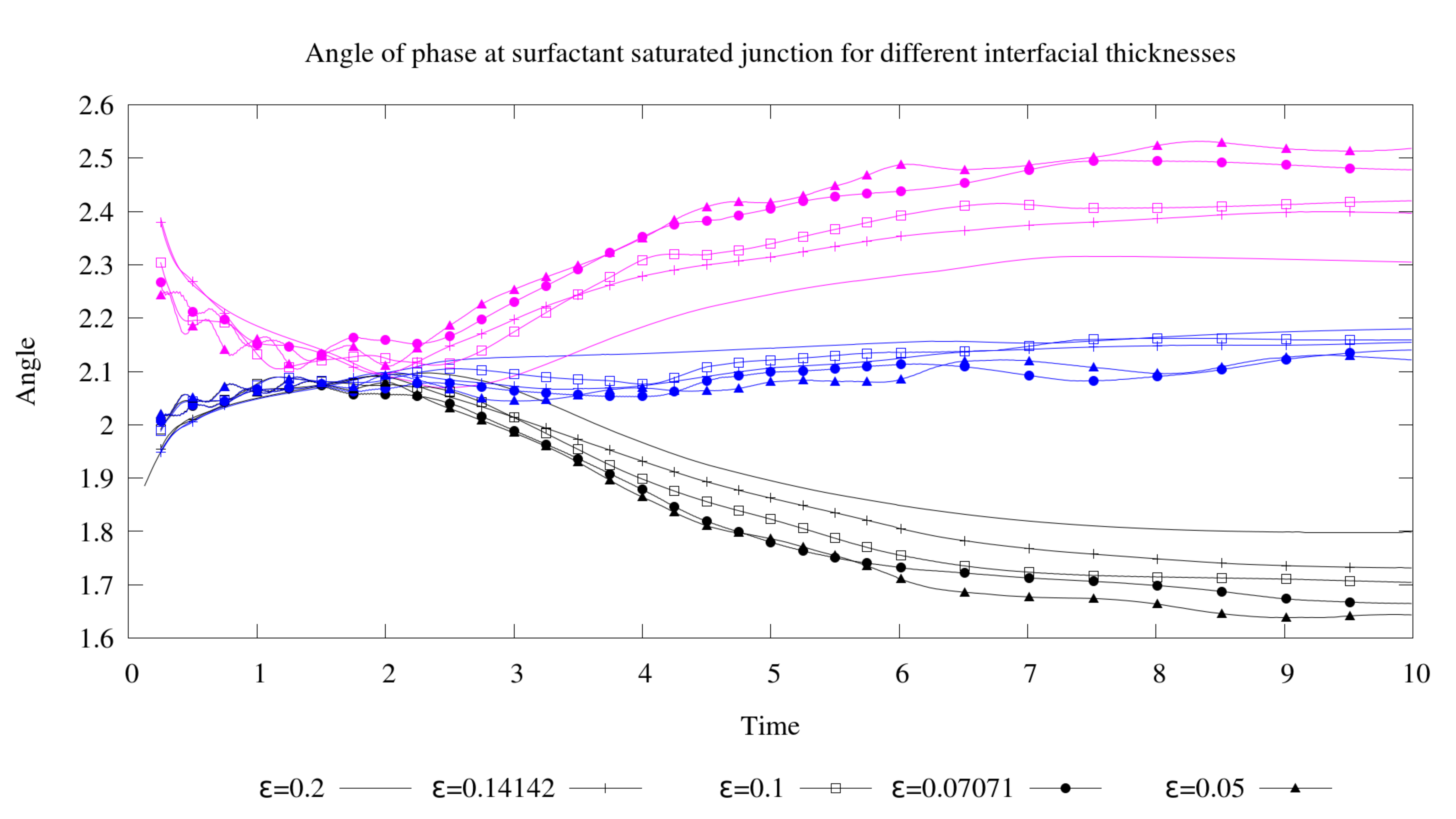}
        \caption{\em Evolution of measured angles $\psi^{(\cdot)}_u$ (unanchored regression) through triple junction (black: $\ph[1]$, blue: $\ph[2]$, pink: $\ph[3]$) for different values of $\e$.}
        \label{fig:TJplot2}
\end{figure}

\begin{table}[ht]
 \centering
\begin{tabular}{|c||c|c|c|c|c|}\hline
$\e$ &  $\psi_a^{(1)}(T)$& $\psi_a^{(2)}(T)$ & $\psi_a^{(3)}(T)$ & $\|\psi_a(T)-\text{ref}\|_2$ & EOC\\ \hline 
0.2 &		1.81433	&2.16857&	2.30029&0.40712 & --\\ \hline
0.14142 &	1.7732&	2.14781&	2.36217&0.33056 & 0.60113\\ \hline
0.1 &  		1.7446&	2.14421	&2.39438&  0.28756 &0.40202\\ \hline
0.070711&	1.72379&	2.12526&	2.43414& 0.24117 &0.50768\\ \hline
0.05 & 		1.70546&	2.12579	&2.44993  &0.21764&0.29624\\ \hline \hline
ref&		1.57079& 2.09439& 2.61799& -- & --\\ \hline
\end{tabular}
\caption{\em Measured angles $\psi^{(\cdot)}_a$ (anchored regression) for the problem in Section \ref{sec:angles} at the final time $T=10$, with reference values corresponding to the sharp interface model (see \eqref{eq:angles_target}.}
\label{tab:anchanglessurf}
\end{table}
\begin{table}[ht]
 \centering
\begin{tabular}{|c||c|c|c|c|c|}\hline
$\e$ &  $\psi_u^{(1)}(T)$& $\psi_u^{(2)}(T)$ & $\psi_u^{(3)}(T)$& $(\|\psi_u(T)-\text{ref}\|_2$ & EOC\\ \hline 
0.2 &		1.79865	&2.17966&2.30488  & 0.39652 & --\\ \hline
0.14142 &	1.74152	&2.15455&2.37711 & 0.30132 &0.79224\\ \hline
0.1 &  		1.70451	&2.15894&2.41973 & 0.2477&  0.56537\\ \hline
0.070711&	1.66479	&2.14011&2.47828 & 0.17449& 1.01096\\ \hline
0.05 & 		1.64349	&2.12166&2.51804 & 0.12657& 0.92643\\ \hline \hline
ref&		1.57079& 2.09439& 2.61799& -- & --  \\ \hline
\end{tabular}
\caption{\em Measured angles $\psi^{(\cdot)}_u$ (unanchored regression) for the problem in Section \ref{sec:angles} at the final time $T=10$, with reference values corresponding to the sharp interface model (see \eqref{eq:angles_target}.}
\label{tab:unanchanglessurf}
\end{table}

We run the simulation for $\e=0.2/\sqrt{2^k}$, $k=0,1,2,3,4$. The measurements of the angles over time displayed in Figure \ref{fig:TJplot} for the anchored case and in Figure \ref{fig:TJplot2} for the unanchored case. We observe convergence of the angles as $\e \to 0$ in both cases. Tables \ref{tab:anchanglessurf} and \ref{tab:unanchanglessurf} list the angles at the final time $T$, along with the target angles, the differences, and EOCs. In the unanchored case we notice a good approximation rate of a bit less than 1. The anchored convergence rate is poorer but, comparing Figures \ref{fig:TJplot} and \ref{fig:TJplot2}, we notice that their measurement is a bit less volatile over time.

\subsection{Marangoni effect}
\label{sec:maran}

We now demonstrate the capability of the new model to describe the effects of Marangoni forces. For each $\e$ we relax a liquid lens to equilibrium, and then introduce surfactant into the domain via the boundary. We retain a persistent surfactant gradient with boundary conditions, and this induces a Marangoni force along some of the interfaces and causes the lens to move. At a final time $T=10$ we compare the relative positions of the lenses for different values of $\e$.

The initial configuration similar to Figure \ref{fig:LensInitHalf} (left) after reflection about the left boundary, $\O[3]$ forms a disc-shaped lens between $\O[1]$ and $\O[2]$. The domain is $\Omega = (-2,4)\times(-2,2)$ and we choose the disc centred at the origin, with radius $r=1$.
The viscosities and densities of each subregion are fixed in time and matched across the different phases. We set $\et[i]=0.01$, $\r[i]=0.1$ for all $i=1,2,3$, which corresponds to a Reynolds number of $10$, and leads to constant functions $\eta(\pha)=0.01$ and $\rho(\pha)=0.1$. The initial velocity is zero, and we impose a homogeneous Dirichlet boundary condition. For the surfactant we choose free energies as in \eqref{eq:toyint} and \eqref{eq:toybulk} with $\sigma_0=1$, $\beta_{i}=1$ and $\beta_{i,j}=0.2$ for all $i=1,2,3$ and $(i,j)=(1,2),(1,3),(2,3)$. The surfactant mobilities we choose as $\Mc[i]=10$ and $\Mcc[i,j]=50$. Initially, there is no surfactant present within the domain. At time $T_0=2.4$ we introduce it with a source at the left hand boundary and a sink along the right hand boundary. We realise this boundary condition linearly over a time $T_q=0.05$: for $z_L \in \{ (-2,y) | y \in [-2,2] \}$ and $z_R \in \{ (4,y) | y \in [-2,2] \}$
\begin{equation*}
 q(z_L,t) = \begin{cases}
          0,& \text{ for } t \leq T_0,\\
	  10(t-T_0), & \text{ for }  T_0 \leq t \leq T_q,\\
	  0.5, & \text{ for } t \geq T_q, \\
          \end{cases}
\qquad \text{and} \qquad q(z_R,t) = 0, \text{ for } t \geq 0.
\end{equation*}
On the other parts of the boundary, a homogeneous Neumann condition is imposed. 

To initialise the phase field variables we proceed similarly as in Section \ref{sec:angles} by using the leading order profile to smooth the initial hypersurfaces $\G[i,j](0)$. The velocity $\vv$ is initialised with zero. The phase field potentials are given as in Section \ref{sec:BLM}, where we set $M_c=0.005$ and $\Lambda=0.1$. 

\begin{figure}[ht]
        \centering 
        \includegraphics[width=0.48\textwidth]{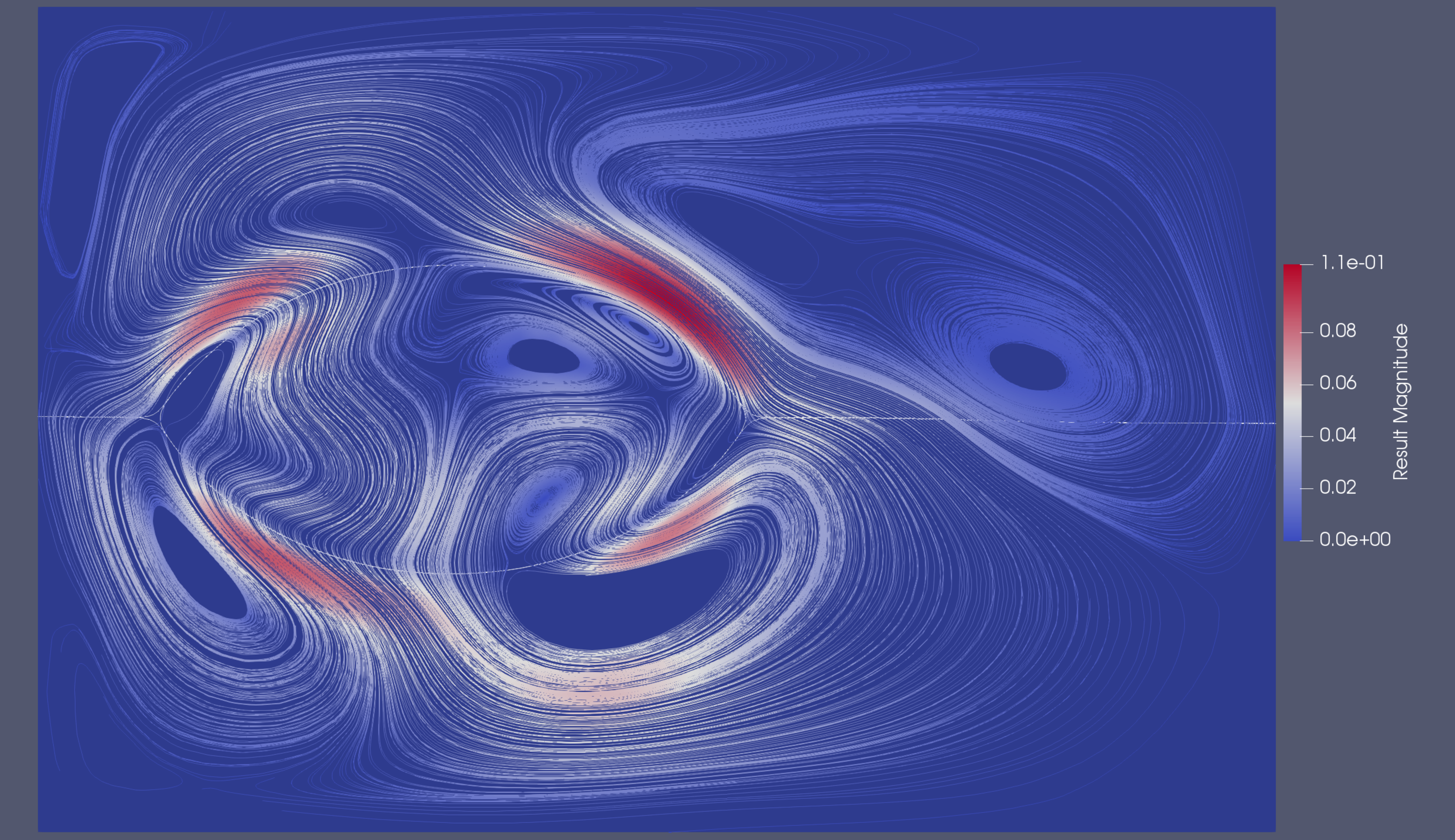} \hfill 
        \includegraphics[width=0.48\textwidth]{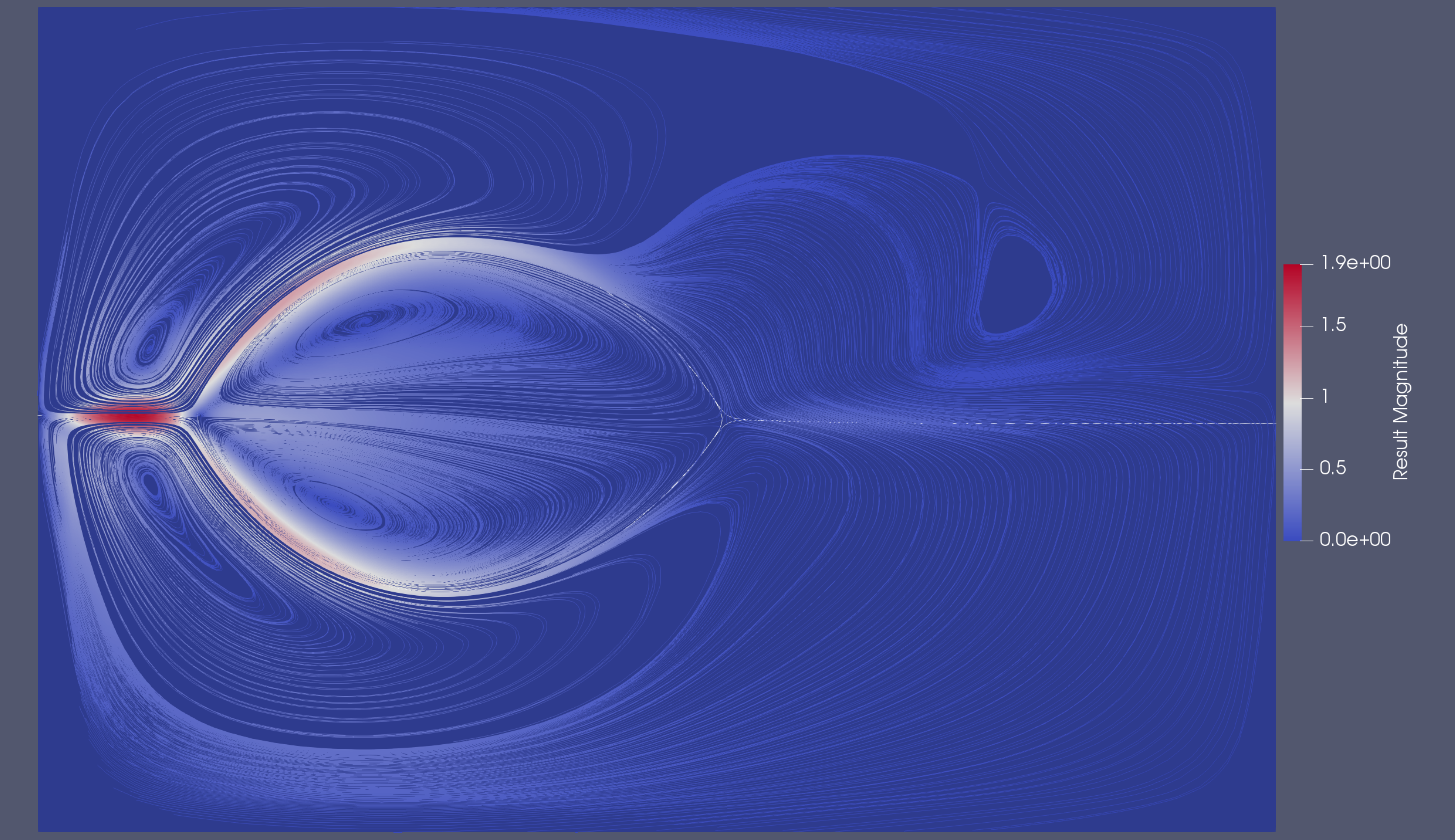} \\[8pt] 
        \includegraphics[width=0.48\textwidth]{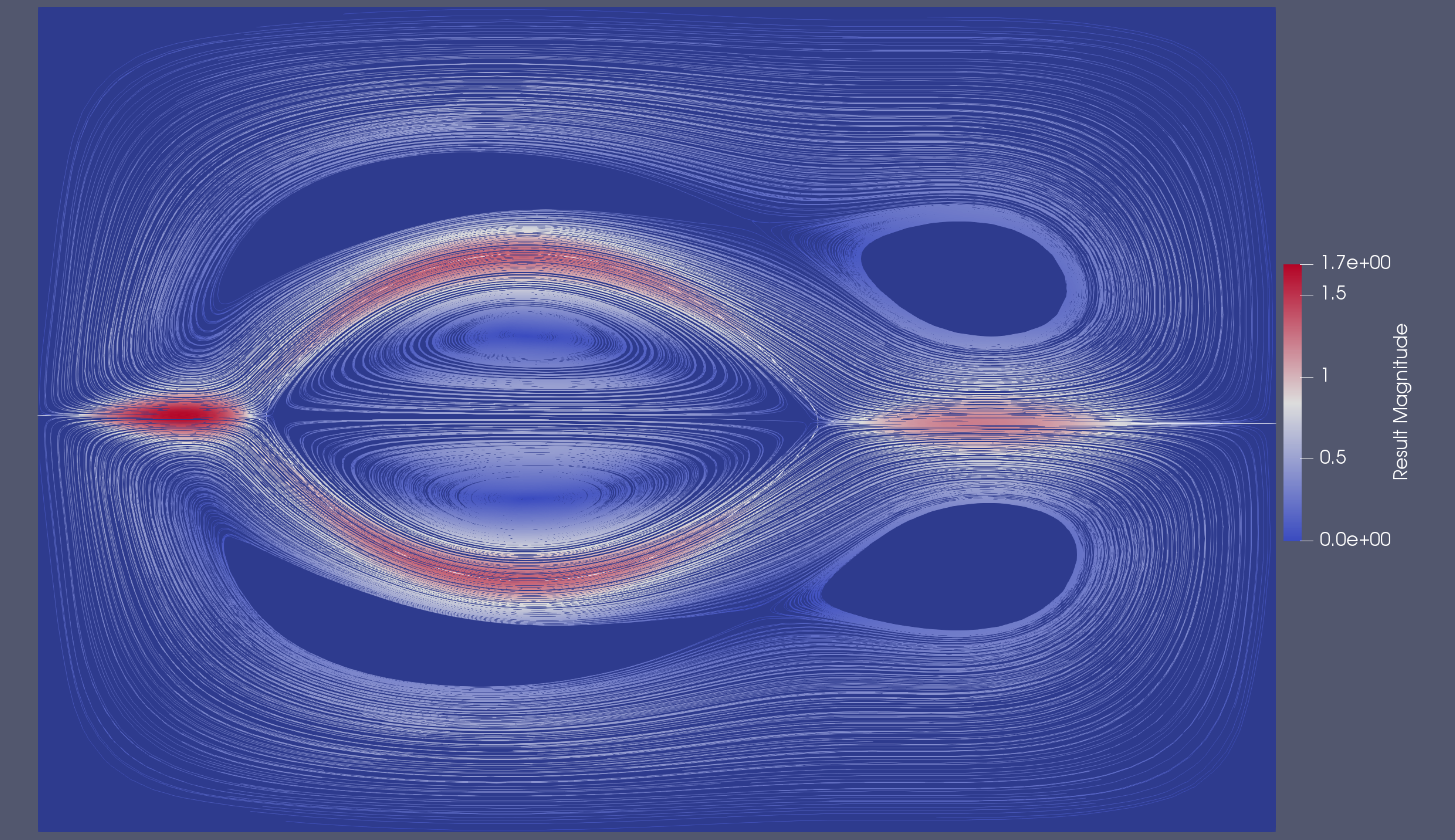}
        \caption{\em Snapshots of a simulation for $\e=0.1$. The $\ph[k]=1/2$ level sets in white represent the phase interfaces. Fluid streamlines are coloured by the fluid velocity magnitude. Top left: $t=2.2<T_0$, no surfactant present. Top right: $t=3.0$ soon after the introduction of surfactant. Bottom: $t=T=10.0$ at the final time.}
        \label{fig:maranflow}
\end{figure}

\begin{figure}[ht]
        \centering 
        \includegraphics[width=0.70\textwidth]{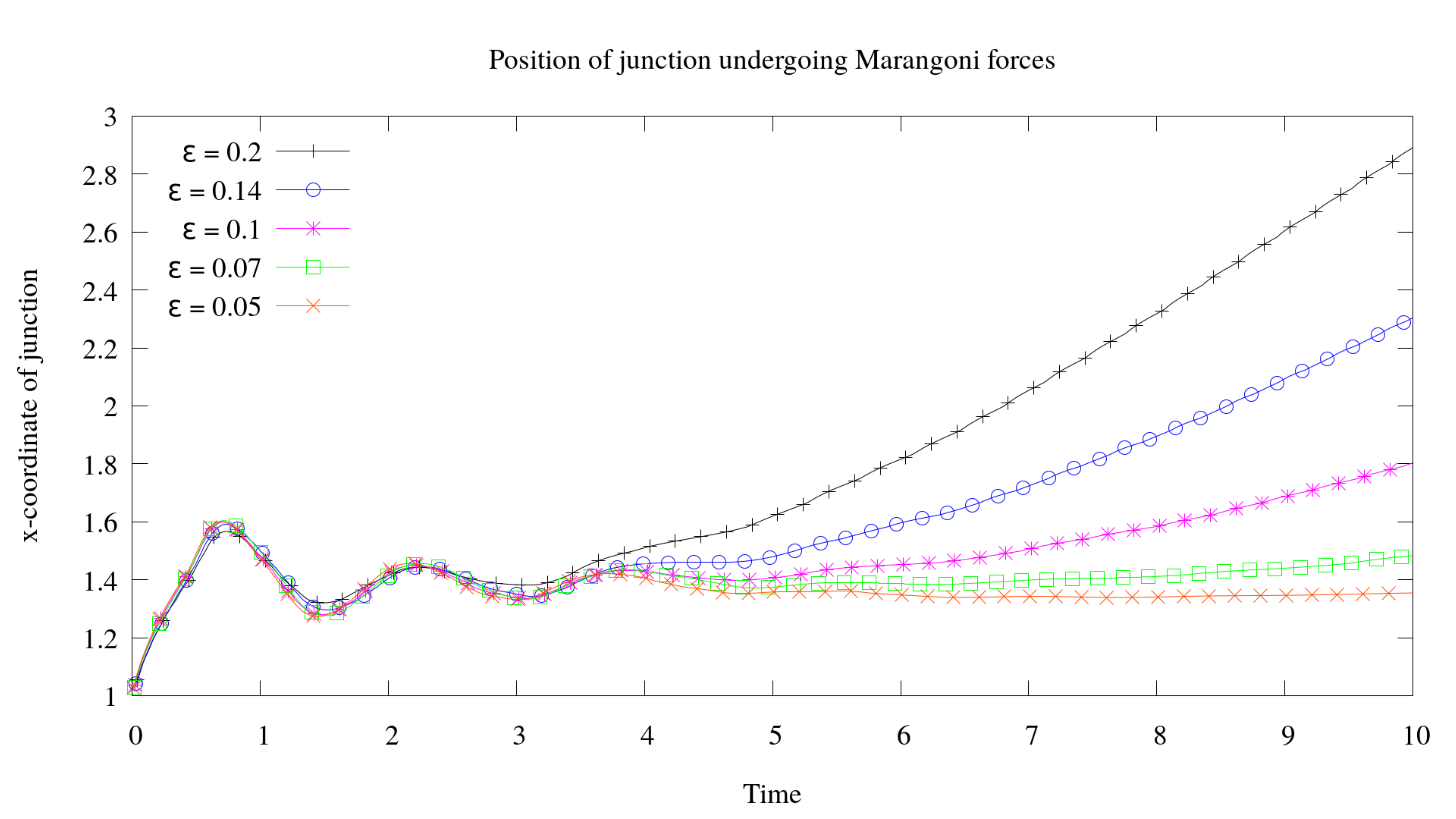} \\
         \includegraphics[width=0.70\textwidth]{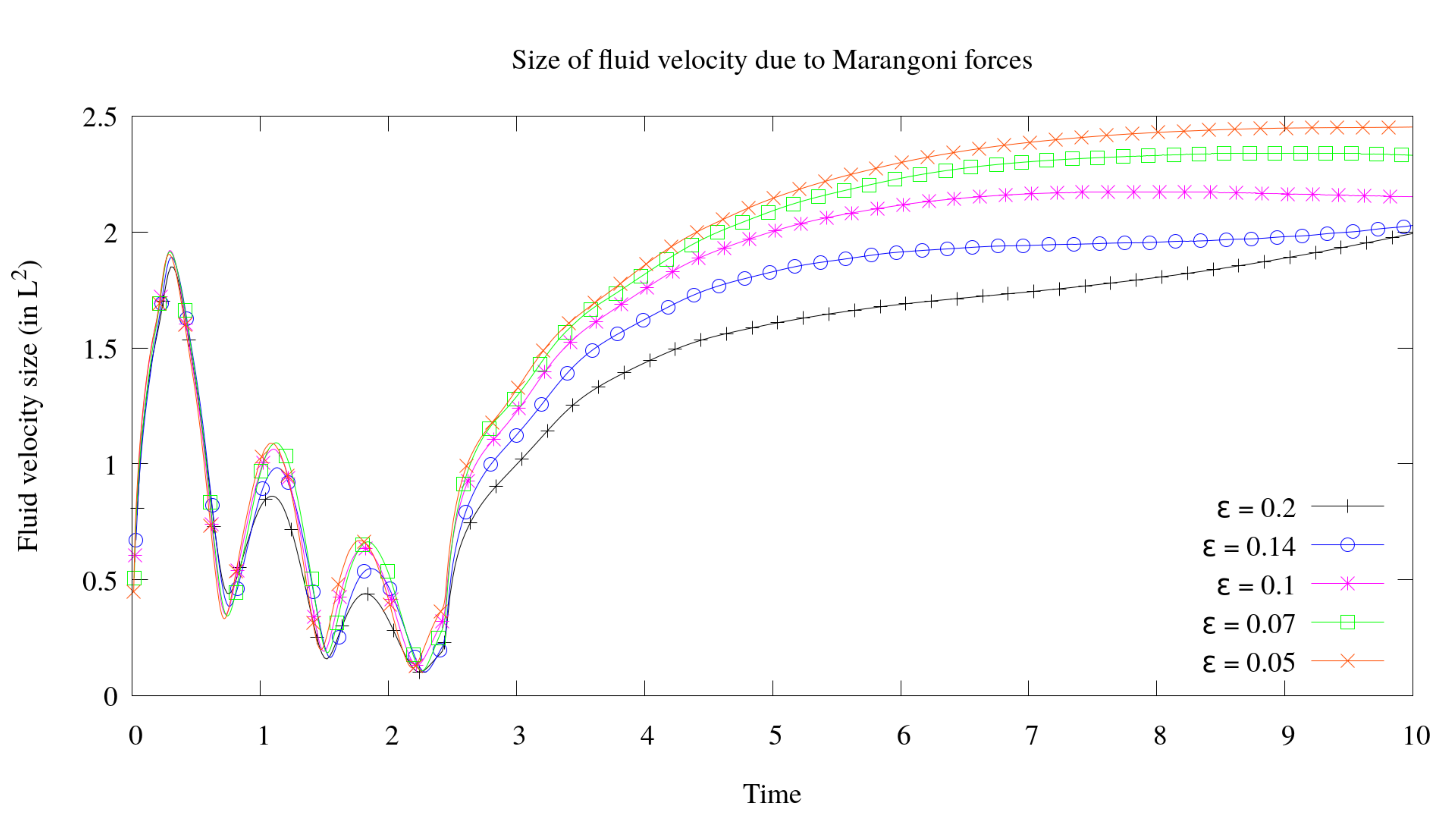}
        \caption{\em Graphs showing the system behaviour over time for different $\e$ ($\e=0.2/\sqrt{2^k}$ $k=0,1,2,3,4$). Top: $x$--coordinate of the left triple junction (see also Figure \ref{fig:marandiffe}). Bottom: $L^2$ norm of the fluid velocity.}
        \label{fig:marantripx}
\end{figure}

\begin{figure}[ht]
        \centering 
        \includegraphics[width=0.6\textwidth]{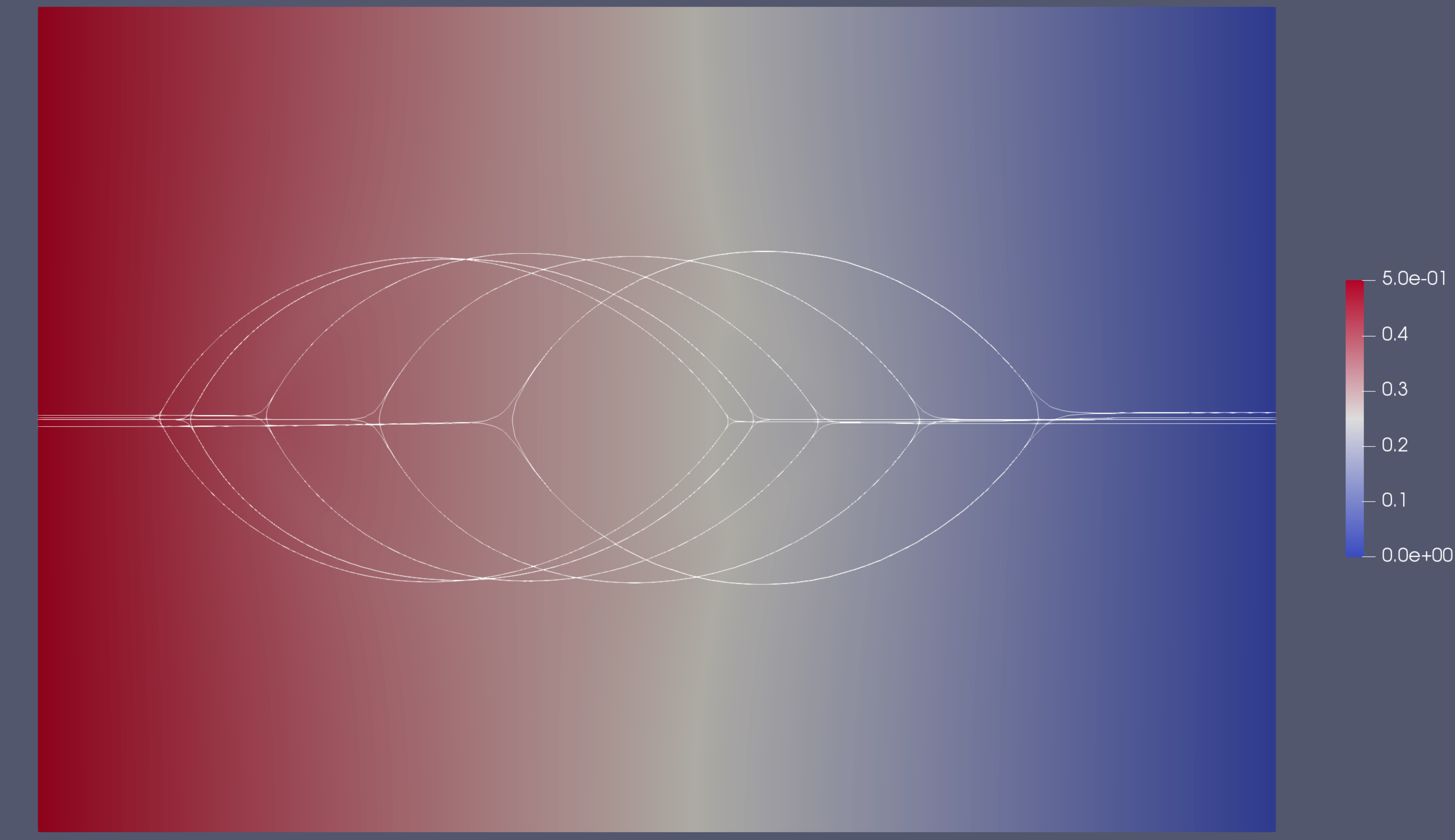}
        \caption{\em Positions of the lenses for different values of $\e$ ($\e=0.2/\sqrt{2^k}$, $k=4,3,2,1,0$ from left to right) at the end of the simulations at time $t=T=10.0$. The $\ph[k]=0.5$ level sets are in white to show the phase interfaces. The colour is determined by the concentration gradient of surfactant.}
        \label{fig:marandiffe}
\end{figure}

We show typical progression of a run with $\e =0.1$ in Figure \ref{fig:maranflow}. The effects of the Marangoni forces for different values of $\e$ are displayed in Figure \ref{fig:marantripx} (where we show the position of the left triple junction and $L^2$ norm of the velocity) and Figure \ref{fig:marandiffe} (giving an impression of the position of the lens). 

Regarding Figure \ref{fig:marantripx} there appears to be convergence of the solutions as $\e$ is decreased. The larger $\e$ solutions display less oscillations until the surfactant is introduced, which can be explained by the (generally known) smoothing effect of the diffuse interface on the fluid flow. At least on the time scale of our observation the impact of the surfactant introduced thereafter is also more pronounced for larger $\e$, indicating that the interfacial forces are dominating inertial forces more strongly. However, we also note in the second plot of Figure \ref{fig:marantripx} that the $L^2$ norm of the velocity increases as $\e$ decreases once the surfactant effect is present. A possible explanation for these observations is that, when the interfacial layer is thicker, there is a greater volume of the fluid for the Marangoni forces to take effect. Vice versa, the bulk droplet volume relative to the interfacial layer volume is reduced for large $\e$, whence there is a smaller region for bulk viscous forces to produce inertia to the motion. This leads to higher dominance of interfacial forces (i.e., Marangoni forces) that move the droplet faster to the right in the large $\e$ regime. We finally remark that, as $t$ approaches $T=10$, there are some boundary effects visible for $\e=0.2,\ 0.1\sqrt{2}$ ($L^2$ norms picking up) as the droplets then approach the right domain boundary.

\section{Conclusion}
\label{sec:conclusion}

We have derived a general moving boundary problem for multi-phase flow of immiscible, incompressible fluids with surfactant, the governing equations of which are presented in Section \ref{sec:sumSIM}. The surfactant is subject to advection-diffusion equations in the bulk and on the interfaces and impacts on the flow via the capillary term and the Marangoni force. A general phase field model is then derived and summarised in Section \ref{sec:sumDIM} following the same procedures. A detailed asymptotic analysis has been performed, which links the two models in the sense that the sharp interface limit of the phase field model is the moving boundary problem. Some numerical simulations of surfactant diffusion through a stationary triple junction, of changing triple junction angles due to changes in the surfactant densities, and of a Marangoni effect showcase the capability of the model and support the results of the asymptotic analysis. 

We have restricted our considerations to chemical equilibrium of the surfactant at the interfaces (instantaneous adsorption) and in the triple junctions. A generalisation to multiple surfactants seems relatively straightforward as $\cc[i,j]$ and $\c[i]$ could be vector-valued fields. We also plan to address the non-instantaneous case in forthcoming work, as well as the details of the numerical method that was used for the simulations in Section \ref{sec:numsim}.

\subsection*{Acknowledgments}

The research was supported by Engineering and Physical Sciences Research Council (EPSRC) grant EP/H023364/1. The third author would also like to thank the Isaac Newton Institute for Mathematical Sciences, Cambridge, for its hospitality during the programme Coupling Geometric PDEs with Physics for Cell Morphology, Motility and Pattern Formation supported by EPSRC grant EP/K032208/1.

\appendix
\section{Appendix}

We state some useful calculus identities on and with moving surfaces (for instance, see \cite{Bet86} and \cite{CerFriGur05}).

\noindent{\em Reynold's transport identity}: For a time dependent domain $\Omega(t) \subset \R^d$ with exterior unit normal $\nu$ and with associated velocity field $\bbb[v]$ (here note necessarily divergence free) and for a field $f(t) : \Omega(t) \to \R$ we have that 
\begin{equation} \label{eq:Rey1}
 \frac{\dd }{\dd t } \int_{\Omega(\cdot)} f(\cdot) \Big|_t = \int_{\Omega(t)} \p_t f(t) + \int_{\p \Omega(t)} f(t) \, \bbb[v](t) \cdot \nu(t) = \int_{\Omega(t)} \md[\v] f(t) + f(t) \, \nabla \cdot \bbb[v](t).
\end{equation}
For a time dependent hypersurface $\Gamma(t)$ with velocity $\bbb[v]$ and for a field $f(t) : \Gamma(t) \to \R$ we have that 
\begin{equation} \label{eq:Rey2}
 \frac{\dd}{\dd t} \int_{\Gamma(\cdot)} f(\cdot) \Big|_t = \int_{\Gamma(t)} \md[\v] f(t) + f(t) \nabla_{\Gamma(t)} \cdot \bbb[v](t).
\end{equation}

\noindent{\em Gauss-Green Formula}: For an orientable hypersurface $\Gamma$ with unit normal $\bbb[\nu]$ and with outward unit conormal $\bbb[\mu]$ on $\p \Gamma$ and for any differentiable scalar function $f : \Gamma \to \R$ we have that
\begin{equation}\label{GG}
 \int_\Gamma \nabla_\Gamma f = -\int_\Gamma f \bbb[\kappa] + \int_{\p \Gamma} f \bbb[\mu],
\end{equation}
with the 
\[
\mbox{curvature vector } \bbb[\kappa] = \nabla_\Gamma \cdot \bbb[\nu] \, \bbb[\nu].
\]
Equivalently, for any differentiable vector field $\bbb[w] : \Gamma \to \R^d$
\begin{equation} \label{div}
 \int_{\Gamma} \nabla_\Gamma \cdot \bbb[w] = -\int_{\Gamma} \bbb[w] \cdot \bbb[\kappa] + \int_{\p \Gamma} \bbb[w] \cdot \bbb[\mu].
\end{equation}

\end{document}